\documentclass[12pt]{amsart}
\DeclareFontFamily{OML}{script}{}
\DeclareFontShape{OML}{script}{m}{it}
{ <5-20> rsfs10 }{}
\DeclareMathAlphabet{\mathscript}{OML}{script}{m}{it}

\renewcommand{\mathcal}[1]{{\mathscript #1}\hspace{0.2ex}}
\usepackage{color}
\ifx\red\undefined
\newcommand{\red}{\color{red}}

\fi
\usepackage[ansinew]{inputenc}
\usepackage[encapsulated]{CJK}

\usepackage{graphicx,graphics}
\usepackage{cite}
\usepackage{pifont}
\usepackage{amsthm}
\usepackage{subfigure}
\usepackage{amscd}
\usepackage{amsmath}
\usepackage{latexsym}
\usepackage{amsfonts}
\usepackage{amssymb}
\usepackage{color}
\usepackage{multicol}
\usepackage{hyperref}
\hypersetup{hypertex=true,
            colorlinks=true,
            linkcolor=blue,
            anchorcolor=blue,
            citecolor=blue}
\usepackage{bm}
\allowdisplaybreaks[4]

\newcommand{\Rmnum}[1]{\uppercase\expandafter{\romannumeral #1}}

\def\s{\sigma}
\renewcommand{\epsilon}{\varepsilon}

\ifx\text\undefined
\newcommand{\text}{\mbox}
\fi
\ifx\operatorname\undefined
\newcommand{\operatorname}{\mathop}
\fi

\newcommand\be{\begin{equation}}
\newcommand\ee{\end{equation}}
\newcommand\bea{\begin{eqnarray}}
\newcommand\eea{\end{eqnarray}}
\newcommand\beaa{\begin{eqnarray*}}
\newcommand\eeaa{\end{eqnarray*}}

\newenvironment{eqa}{\begin{equation}%
  \begin{array}{rcl}}{\end{array}\end{equation}}

\newcommand\beqa{\begin{eqa}}
\newcommand\eeqa{\end{eqa}}

\numberwithin{equation}{section}

\renewcommand{\tilde}{\widetilde}
\renewcommand{\hat}{\widehat}

\newtheorem{thm}{Theorem}[section]

\newtheorem{lem}[thm]{Lemma}

\newtheorem{rem}{Remark}[section]

\newcommand{\void}[1]{}

\numberwithin{equation}{section}

\begin{document}\begin{CJK}{UTF8}{gkai}

\title[Non-equilibrium-diffusion limit]{Non-equilibrium-diffusion limit of the compressible Euler radiation model}
\author{Qiangchang Ju, Lei Li and Zhengce Zhang}
\date{\today}
\address[Qiangchang Ju]{Institute of Applied Physics and Computational Mathematics, Beijing, 100094, P. R. China}
\email{ju\_qiangchang@iapcm.ac.cn}
\address[Lei Li]{School of Mathematics and Statistics, Xi'an Jiaotong University, Xi'an, 710049, P. R. China}
\email{lileimathmatica@outlook.com}
\address[Zhengce Zhang]{School of Mathematics and Statistics, Xi'an Jiaotong University, Xi'an, 710049, P. R. China}
\email{zhangzc@mail.xjtu.edu.cn}
\thanks{Corresponding author: Zhengce Zhang}

\maketitle
\begin{abstract}
 We justify rigorously the non-equilibrium-diffusion limit of the compressible Euler model coupled with a radiative transfer equation arising in radiation hydrodynamics. For general initial data, we establish the uniform existence of the solution to the coupled model in $\mathbb{T}^{3}$ and prove the convergence of the solutions to the limiting system in the nonequilibrium-diffusion regime. Moreover, the initial layer for the radiative density is constructed to get the strong convergence in $L^\infty$ norm.
\end{abstract}
\section{Introduction}
The Euler system coupled with a radiative transfer equation (see \cite{bu-8, bu-9, bu-10}) can describe the motion of a compressible inviscid fluid under the presence of the radiative field as follows
\begin{equation}\label{research equations-x-1}\left\{
\begin{split}
&\partial_{t}\rho+\mathrm{div}_{x}(\rho\vec{u})=0,\\
&\partial_{t}(\rho\vec{u})+\mathrm{div}_{x}(\rho\vec{u}\otimes\vec{u})+\nabla_{x}P
=-\vec{S}_{F},\\
&\partial_{t}(\rho e)+\mathrm{div}_{x}(\rho e\vec{u})+P\mathrm{div}_{x}\vec{u}
=-S_{E}+\vec{S}_{F}\cdot\vec u,\\
&\frac{1}{c}\partial_{t}f+\vec w\cdot\nabla_{x}f=S.\\
\end{split}\right.
\end{equation}
Here the unknowns $\rho, \vec u=(u^{1}, u^{2}, u^{3})$ and $e$ denote the density, velocity, and specific internal energy, respectively. $f=f(t, \vec x, \nu, \vec w)$ is the radiative intensity which depends on time $t\geq 0$, the special variable $\vec x\in\Omega:=\mathbb{T}^{3}=(\mathbb{R}\setminus(2\pi\mathbb{Z}))^3$, the frequency $\nu\geq 0$ and direction $\vec w\in\mathbb{S}^{2}$ of photons with $\mathbb{S}^{2}\subset\mathbb{R}^{3}$ being the unit sphere. The pressure $P=P(\rho, \theta)$ and the internal energy $e=e(\rho, \theta)$ are smooth functions of $\rho$ and the temperature $\theta$.

We suppose that the radiative source term $S$ in $(\ref{research equations-x-1})_{4}$ is given by
\begin{equation}\nonumber
S=\sigma_{a}[B(\nu, \theta)-f(t, \vec x, \nu, \vec w)]+\sigma_{s}\Big(\frac{1}{4\pi}\int_{\mathbb{S}^{2}}f(t, \vec x, \nu, \vec w)\mathrm{d}\vec w-f(t, \vec x, \nu, \vec w)\Big),
\end{equation}
where $\sigma_{a}=\sigma_{a}(\nu, \theta)$ and $\sigma_{s}=\sigma_{s}(\nu, \theta)$ are (given) nonnegative functions and represent the absorption and the scattering coefficient, respectively. Moreover, the radiative equilibrium function $B(\nu, \theta)$ is given by
\begin{equation}\nonumber
B(\nu, \theta)=2h\nu^{3}c^{-2}\Big(e^{\frac{h\nu}{k_{B}\theta}}-1\Big)^{-1},
\end{equation}
where $c, k_{B}$ and  $h$ are the speed of light, the Boltzmann and Planck constants, respectively. Furthermore, the radiative flux $\vec{S}_{F}$ and the radiative energy source $S_{E}$ are defined by
\begin{equation}\nonumber
\vec{S}_{F}=\frac{1}{c}\int_{0}^{\infty}\int_{\mathbb{S}^{2}}\vec{w}S\mathrm{d}\vec w\mathrm{d}\nu
\end{equation}
and
\begin{equation}\nonumber
S_{E}=\int_{0}^{\infty}\int_{\mathbb{S}^{2}}S\mathrm{d}\vec w\mathrm{d}\nu.
\end{equation}

The system (\ref{research equations-x-1}) is a simplified model in radiation hydrodynamics \cite{bu-10}. More realistic systems have been proposed by
Lowrie, Morel and Hittinger \cite{bu-9} and Buet and Despr\'es \cite{bu-8}.

In order to identify the appropriate limit regime, we nondimensionalize the system \eqref{research equations-x-1} as in \cite{DN1, bu-6} by setting
\begin{equation}\label{equation}\nonumber\left\{
\begin{split}
&\hat{t}=\frac{t}{T_{\infty}},~~~ \hat{\vec{x}}=\frac{\vec x}{L_{\infty}},~~~ \hat{\vec u}=\frac{\vec u}{U_{\infty}},~~~ \hat\rho=\frac{\rho}{\rho_{\infty}},~~~ \hat\theta=\frac{\theta}{\theta_{\infty}},~~~ \hat P=\frac{P}{P_{\infty}},\\
&\hat e=\frac{e}{e_{\infty}},~~~ \hat{\sigma_{a}}=\frac{\sigma_{a}}{\sigma_{a, \infty}},~~~ \hat{\sigma_{s}}=\frac{\sigma_{s}}{\sigma_{s, \infty}},~~~ \hat{f}=\frac{f}{f_{\infty}},~~~ \hat{\nu}=\frac{\nu}{\nu_{\infty}}.\\
\end{split}\right.
\end{equation}
Here we use $T_{\infty}, L_{\infty}, U_{\infty}, \rho_{\infty}, \theta_{\infty}, P_{\infty}$ and $e_{\infty}$ to denote the reference hydrodynamical quantities (time, length, velocity, density, temperature, pressure, energy). Furthermore, we use $f_{\infty}, \nu_{\infty}, \sigma_{a, \infty},$ and $\sigma_{s, \infty}$ to denote the reference radiative quantities (radiative intensity, frequency, absorption and scattering coefficients). Besides, we use the carets to symbolize nondimensional quantities. Moreover, we also assume the compatibility conditions $P_{\infty}=\rho_{\infty}e_{\infty}, \nu_{\infty}=\frac{k_{B}\theta_{\infty}}{h}, f_{\infty}=\frac{2h\nu_{\infty}^{2}}{c^{2}}$ and we denote by
\begin{equation}\nonumber
S_{r}:=\frac{L_{\infty}}{T_{\infty}U_{\infty}},~~~\; Ma:=\frac{U_{\infty}}{\sqrt{P_{\infty}/\rho_{\infty}}},~~~\; \mathcal{C}:=\frac{c}{U_{\infty}}
\end{equation}
the Strouhal, Mach and ``infrarelativistic" numbers corresponding to hydrodynamics, and by
\begin{equation}\nonumber
\mathcal{L}:=L_{\infty}\sigma_{a, \infty},~~~ \;\mathcal{L}_{s}:=\frac{\sigma_{s, \infty}}{\sigma_{a, \infty}},~~~\; \mathcal{P}:=\frac{2k_{B}^{4}\theta_{\infty}^{4}}{h^{3}c^{3}\rho_{\infty}e_{\infty}},
\end{equation}
the various dimensionless numbers corresponding to radiation. Using above scalings and omitting the carets, we can write (\ref{research equations-x-1}) as follows
\begin{equation}\label{research equations-x-2}\left\{
\begin{split}
&S_{r}\partial_{t}\rho+\mathrm{div}_{x}(\rho\vec{u})=0,\\
&S_{r}\partial_{t}(\rho\vec{u})+\mathrm{div}_{x}(\rho\vec{u}\otimes\vec{u})+\frac{1}{Ma^{2}}\nabla_{x}P
=-\frac{\mathcal{P}}{Ma^{2}}\vec{S}_{F},\\
&S_{r}\partial_{t}(\rho e)+\mathrm{div}_{x}(\rho e\vec{u})+P\mathrm{div}_{x}\vec{u}
=-\mathcal{P}\mathcal{C}S_{E}+\mathcal{P}\vec{S}_{F}\cdot\vec u,\\
&\frac{S_{r}}{\mathcal{C}}\partial_{t}f+\vec{w}\cdot\nabla_{x}f=S,\\
\end{split}\right.
\end{equation}
where
\begin{equation}\nonumber\left\{
\begin{split}
&S_{F}=\int_{0}^{\infty}\int_{\mathbb{S}^{2}}\vec w S\mathrm{d}\vec w \mathrm{d}\nu,~~~ \;S_{E}=\int_{0}^{\infty}\int_{\mathbb{S}^{2}}S\mathrm{d}\vec w \mathrm{d}\nu,\\
&S=\mathcal{L}\sigma_{a}(B(\nu, \theta)-f)+\mathcal{L}\mathcal{L}_{s}\sigma_{s}(\overline f-f),\\
&B(\nu, \theta)=\frac{\nu^{3}}{e^{\frac{\nu}{\theta}}-1},~~~\; \overline f=\frac{1}{4\pi}\int_{\mathbb{S}^{2}}f\mathrm{d}\vec w.
\end{split}\right.
\end{equation}

In the present paper, we consider the grey hypothesis such that the transport coefficients $\sigma_{a}$ and $\sigma_{s}$ are independent of the frequency $\nu$. For simplicity, we assume that $\sigma_{a}$ and $\sigma_{s}$ are two positive constants in the following derivation. It is worthy to note that the general case $\sigma_{a}=\sigma_{a}(\rho, \theta)$ and $\sigma_{s}=\sigma_{s}(\rho, \theta)$ can be dealt with similarly. Integrating the radiative transfer equation in (\ref{research equations-x-2}) with respect to the frequency $\nu$, and noticing that
\begin{equation}\nonumber
\int_{0}^{\infty}B(\nu, \theta)\mathrm{d}\nu=\int_{0}^{\infty}\frac{\nu^{3}}{e^{\frac{\nu}{\theta}}-1}\mathrm{d}\nu=C\theta^{4}:=B(\theta)
\end{equation}
for some positive constant $C$, we get
\begin{equation}\nonumber
\begin{aligned}
\frac{Sr}{\mathcal{C}}\partial_{t}\int_{0}^{\infty}f\mathrm{d}\nu+\vec w\cdot\nabla_{x}\int_{0}^{\infty}f \mathrm{d}\nu=&\mathcal{L}\sigma_{a}\Big(B(\theta)-\int_{0}^{\infty}f\mathrm{d}\nu\Big)
\\&+\mathcal{L}\mathcal{L}_{s}\sigma_{s}\Big(\frac{1}{4\pi}\int_{\mathbb{S}^{2}}\int_{0}^{\infty}f\mathrm{d}\nu \mathrm{d}\vec{w}-\int_{0}^{\infty}f\mathrm{d}\nu\Big).
\end{aligned}
\end{equation}
Moreover, if we denote $\int_{0}^{\infty}f \mathrm{d}\nu$ by $f$, then the above equation can be written as
\begin{equation}\nonumber
\begin{aligned}
\frac{Sr}{\mathcal{C}}\partial_{t}f+\vec w\cdot\nabla_{x}f=\mathcal{L}\sigma_{a}(B(\theta)-f)
+\mathcal{L}\mathcal{L}_{s}\sigma_{s}\Big(\frac{1}{4\pi}\int_{\mathbb{S}^{2}}f\mathrm{d}\vec{w}-f\Big).
\end{aligned}
\end{equation}
Defining $\langle f\rangle=\overline{f}$, we can rewrite (\ref{research equations-x-2}) after careful calculations as follows
\begin{equation}\label{research equations-x-5}\left\{
\begin{split}
&S_{r}\partial_{t}\rho+\mathrm{div}_{x}(\rho\vec{u})=0,\\
&S_{r}\partial_{t}(\rho\vec{u})+\mathrm{div}_{x}(\rho\vec{u}\otimes\vec{u})+\frac{1}{Ma^{2}}\nabla_{x}P
=\frac{\mathcal{P}}{Ma^{2}}(\mathcal{L}\sigma_{a}+\mathcal{L}\mathcal{L}_{s}\sigma_{s})\langle\vec{w}(f-\overline f)\rangle,\\
&S_{r}\partial_{t}(\rho e)+\mathrm{div}_{x}(\rho e\vec{u})+P\mathrm{div}_{x}\vec{u}
=-\mathcal{P}\mathcal{C}\mathcal{L}\sigma_{a}(B(\theta)-\overline f)\\&\hspace{6.0cm}-\mathcal{P}(\mathcal{L}\sigma_{a}+\mathcal{L}\mathcal{L}_{s}\sigma_{s})\langle\vec{w}(f-\overline f)\rangle\cdot\vec u,\\
&\frac{Sr}{\mathcal{C}}\partial_{t}f+\vec w\cdot\nabla_{x}f=\mathcal{L}\sigma_{a}(B(\theta)-f)
+\mathcal{L}\mathcal{L}_{s}\sigma_{s}(\overline f-f).\\
\end{split}\right.
\end{equation}

The main purpose of the present paper is to identify the non-equilibrium diffusion regime for the system (\ref{research equations-x-2}) in a three-dimensional torus. So, we put $\mathcal{L}_{s}=1/\epsilon^{2}$ and $\mathcal{L}=\epsilon$. For the strong under-relativistic flow, we put $\mathcal{C}=1/\epsilon$. We also assume that a certain amount of radiation is present, i.e., $\mathcal{P}=1$, and the Mach number is fixed, i.e. $Ma=1$. Moreover, we put $S_{r}=1$ in the previous system. We shall focus on the fluids obeying the ideal polytropic gas relations
\begin{equation}
e=C_{V}\theta, P=R\rho\theta,
\end{equation}
where the specific gas constant $R$ and the specific heat at constant volume $C_{V}$ are positive constants. For simplicity, we take $\sigma_{a}, \sigma_{s}, C_{V}$ and $R$ in (\ref{research equations-x-5}) to be one. As a result, it is equivalent to proving the limit of the solution $(f^{\epsilon}, \rho^{\epsilon}, \vec{u}^{\epsilon}, \theta^{\epsilon})$ to the initial value problem of the following system
\begin{equation}\label{research equations}\left\{
\begin{split}
&\partial_{t}\rho^{\epsilon}+\mathrm{div}_{x}(\rho^{\epsilon}\vec{u}^{\epsilon})=0,~~\mathrm{in} \ \ (0,T)\times \Omega,\\
&\partial_{t}(\rho^{\epsilon}\vec{u}^{\epsilon})+\mathrm{div}_{x}(\rho^{\epsilon}\vec{u}^{\epsilon}\otimes\vec{u}^{\epsilon}+\rho^{\epsilon}\theta^{\epsilon}\mathbb{I}_{3\times 3})
=\Big(\frac{1}{\epsilon}+\epsilon\Big)\langle\vec{w}(f^{\epsilon}-\overline f^{\epsilon})\rangle,~~\mathrm{in} \ \ (0,T)\times \Omega,\\
&\partial_{t}(\rho^{\epsilon}\theta^{\epsilon})+\mathrm{div}_{x}(\rho^{\epsilon}\theta^{\epsilon}\vec{u}^{\epsilon})+\rho^{\epsilon}\theta^{\epsilon}\mathrm{div}_{x}\vec{u}^{\epsilon}
=\overline f^{\epsilon}-B(\theta^{\epsilon})\\&\hspace{6.9cm}-\Big(\frac{1}{\epsilon}+\epsilon\Big)\langle\vec{w}(f^{\epsilon}-\overline f^{\epsilon})\rangle\cdot\vec u^{\epsilon},~~\mathrm{in} \ \ (0,T)\times \Omega, \\
&{\epsilon}^{2}\partial_{t}f^{\epsilon}+\epsilon \vec w\cdot \nabla_{x}f^{\epsilon}+f^{\epsilon}-\overline f^{\epsilon}+{\epsilon}^{2}f^{\epsilon}={\epsilon}^{2}B(\theta^{\epsilon}),~~\mathrm{in} \ \ (0,T)\times \Omega\times \mathbb{S}^{2},\\
&\rho^{\epsilon}(0, \vec{x})=\rho^{0}(\vec x),\;\;\vec{u}^{\epsilon}(0, \vec{x})=\vec{u}^{0}(\vec x)~~\;\mathrm{in}  \ \ \Omega,\\
&\theta^{\epsilon}(0, \vec{x})=\theta^{0}(\vec{x}),\;\;\;f^{\epsilon}(0, \vec{x}, \vec{w})=h(\vec{x}, \vec{w})~~\;\mathrm{in}  \ \ \Omega\\
\end{split}\right.
\end{equation}
as $\epsilon$ tends to zero.

The first theorem in the paper is as follows.
\begin{thm}\label{definition weak solution-0}
Assume that $h(\vec{x}, \vec{w})\in L^{\infty}(\mathbb{S}^{2}; W_{2}^{6}(\Omega))$ and $\rho^{0}(\vec x), \vec{u}^{0}(\vec x), \theta^{0}(\vec{x})\in W_{2}^{6}(\Omega)$ with $h(\vec{x}, \vec{w})>0, \rho^{0}(\vec x), \theta^{0}(\vec x)\geq a$ where $a$ is a positive constant.  There exists a $T>0$ independent of $\epsilon$, such that the problem (\ref{research equations}) has a unique solution $(f^{\epsilon}, \rho^{\epsilon}, \vec{u}^{\epsilon}, \theta^{\epsilon})$ satisfying $(f^{\epsilon}, \rho^{\epsilon}, \vec{u}^{\epsilon}, \theta^{\epsilon})\in G=\{(f,\rho, \vec u, \theta), \|f\|_{L^{\infty}(0, T; L^{2}(\mathbb{S}^{2}; W_{2}^{6}(\Omega)))}+\|(\rho,\vec u,\theta)\|_{L^{\infty}(0, T; W_{2}^{6}(\Omega))}\leq r\}$ for some positive constant $r$. Furthermore, the unique solution $(f^{\epsilon}, \rho^{\epsilon}, \vec{u}^{\epsilon}, \theta^{\epsilon})$ satisfies
\begin{equation}\label{jielunbudengshi 1}
\lim_{\epsilon\rightarrow 0}\|f^{\epsilon}-\overline f\|_{L^{2}(0, T; L^{2}(\mathbb{S}^{2}; W_{2}^{5}(\Omega)))}=0,
\end{equation} and
\begin{equation}\label{jielunbudengshi 1.5}
\lim_{\epsilon\rightarrow 0}\|(\overline f^{\epsilon}-\overline f,\rho^{\epsilon}-\rho, \vec u^{\epsilon}-\vec u, \theta^{\epsilon}-\theta) \|_{C^{0}([0, T]; W_{2}^{5}(\Omega))}=0,
\end{equation}
where $(\overline f, \rho, \vec u, \theta)\in (L^{\infty}(0, T; W_{2}^{6}(\Omega))\cap L^{2}(0, T; W_{2}^{7}(\Omega))\times L^{\infty}(0, T; W_{2}^{6}(\Omega))^3)$ satisfies
\begin{equation}\label{equations about limit of overline f epsilon}\left\{
\begin{split}
&\partial_{t}\rho+\mathrm{div}_{x}(\rho\vec{u})=0,~~\mathrm{in} \ \ (0,T)\times \Omega,\\
&\partial_{t}(\rho\vec{u})+\mathrm{div}_{x}(\rho\vec{u}\otimes\vec{u}+\rho\theta\mathbb{I}_{3\times 3})
=-\frac{1}{3}\nabla_{x}\overline f,~~\mathrm{in} \ \ (0,T)\times \Omega,\\
&\partial_{t}(\rho\theta)+\mathrm{div}_{x}(\rho\theta\vec{u})+\rho\theta\mathrm{div}_{x}\vec{u}
=\overline f-B(\theta)+\frac{1}{3}\nabla_{x}\overline f\cdot\vec{u},~~\mathrm{in} \ \ (0,T)\times \Omega, \\
&\partial_{t}\overline f-\frac{1}{3}\Delta_{x}\overline f+\overline f=B(\theta)\ \ \mathrm{in}\ \ \ (0,T]\times \Omega,\\
&\rho(0, \vec{x})=\rho^{0}(\vec x),\;\;\;\vec{u}(0, \vec{x})=\vec{u}^{0}(\vec x)~~\;\mathrm{in}  \ \ \Omega,\\
&\theta(0, \vec{x})=\theta^{0}(\vec{x}),\;\;\;\overline f(0, \vec{x})=\overline h(\vec{x})~~\;\mathrm{in}  \ \ \Omega.
\end{split}\right.
\end{equation}
\end{thm}

\begin{rem}
In Theorem \ref{definition weak solution-0}, it suffices to assume that $\rho^{0}(\vec x), \vec{u}^{0}(\vec x), \theta^{0}(\vec{x})\in W_{2}^{3}(\Omega)$ and $h(\vec{x}, \vec w)\in L^{2}(\mathbb{S}^{2}; W_{2}^{3}(\Omega))$ to get the local existence and the similar convergences. Here we have assumed more regularity on the data so as to get the uniform convergence of $f^{\epsilon}$  by using the Hilbert expansion in Section \ref{wu} combined with the uniqueness of the solution to the system (\ref{equations about limit of overline f epsilon}).
\end{rem}

One of the difficulties of getting the uniform existence of the system (\ref{research equations}) is how to handle the singular term $(f^{\epsilon}-\overline f^{\epsilon})/\epsilon$. Inspired by \cite{local-existence-1}, we write the equations satisfied by the density, temperature and velocity in the form of (\ref{research equations-3-2}). Then, we consider the linearized system of (\ref{research equations}) in the form of (\ref{research equations-4}) and get the solution sequences $(f_{k+1}, V_{k+1})$ according to Lemma \ref{main result-7-bu} and the similar version of Theorem I in \cite{cauchy-problem}. By defining the energy radius as $r_{k}=\|f_{k}\|_{L^{\infty}(0, T; L^{2}(\mathbb{S}^{2}; W_{2}^{6}(\Omega)))}+\|V_{k}\|_{L^{\infty}(0, T; W_{2}^{6}(\Omega))}$, one can control the singular term $(f_{k+1}-\overline f_{k+1})/\epsilon$ as $\|(f_{k+1}-\overline f_{k+1})/\epsilon\|_{L^{2}(0, T; L^{2}(\mathbb{S}^{2}; W_{2}^{6}(\Omega)))}^{2}
\leq Cr_{k}^{8}T+C\|h\|_{L^{2}(\mathbb{S}^{2}; W_{2}^{6}(\Omega))}^{2}$ (see (\ref{gronwall's inequality 1-1})). The uniform estimates about $k$ and $\epsilon$ and the existence time $T$ can be further obtained by a iterative method. Thus the local existence of the solution to the system (\ref{research equations}) is proved by the Banach contraction mapping principle.

In the next theorem, we will get the uniform convergence of $f^{\epsilon}$ in $L^{\infty}((0, T)\times\Omega\times\mathbb{S}^{2})$ by resorting to the Hilbert expansion.
\begin{thm}\label{definition weak solution-1}
For $f^{\epsilon}(t, \vec{x}, \vec{w})$ and $\overline f$ obtained in Theorem \ref{definition weak solution-0}, it holds that
\begin{equation}\label{jielunbudengshi 1-jia}
\lim_{\epsilon\rightarrow 0}\|f^{\epsilon}-\overline f-f_{I, 0}\|_{L^{\infty}((0, T)\times\Omega\times \mathbb{S}^{2})}=0
\end{equation}
where $f_{I, 0}$ satisfies
\begin{equation}
f_{I, 0}=e^{-\tau}\Big(h(\vec{x}, \vec{w})-\frac{1}{4\pi}\int_{\mathbb{S}^{2}}h(\vec{x}, \vec{w})\mathrm{d}\vec{w}\Big)
\end{equation}
for $\tau=\frac{t}{\epsilon^{2}}$.
\end{thm}

There are many classical results \cite{lionssidawenzhang, geometric-CMP, unsteady-neutron-transport-2-d-unit-disk, convex-domain-123, 3-d-diffusive-limit-convex-domain, 3-d-zuinanqingxing} about the single neutron transport equation with physical boundary. Instead of the single neutron equation, the coupled model (\ref{research equations}) we consider is more complicated and  physical. For simplicity, we consider the coupled model in $\mathbb{T}^{3}$. For the transport equation in the system (\ref{research equations}), there are two extra terms $\epsilon^{2}f^{\epsilon}$ and $\epsilon^{2}B(\theta^{\epsilon})$ involved. The extra term $\epsilon^{2}B(\theta^{\epsilon})$ cause some difficulties for the remainder estimates of $f^{\epsilon}$. Even through the proof of the theorems are mainly motivated by the works on the diffusive limit of the single neutron transport equation (see \cite{unsteady-neutron-transport-2-d-unit-disk, 3-d-zuinanqingxing}), some new ideas and different arguments are needed here to overcome the difficulties.   One of the key steps in our analysis is to control the remainder
 $\mathcal{R}_{2}=f^{\epsilon}-\sum_{k=0}^{2}\epsilon^{k}f_{k}-\sum_{k=0}^{2}\epsilon^{k}f_{I, k}$ satisfying
\begin{equation}\nonumber\left\{
\begin{split}
&{\epsilon}^{2}\partial_{t}\mathcal{R}_{2}+\epsilon \vec w\cdot \nabla_{x}\mathcal{R}_{2}+\mathcal{R}_{2}-\overline {\mathcal{R}_{2}}+{\epsilon}^{2}\mathcal{R}_{2}=\mathcal{L}\mathcal{R}_{2}~~\mathrm{in} \ \ (0,T]\times \Omega\times \mathbb{S}^{2},\\
&\mathcal{R}_{2}(0, \vec{x}, \vec{w})=0~~\mathrm{in}~~~\Omega\times \mathbb{S}^{2},\\
\end{split}\right.
\end{equation}
where the operator is defined in (\ref{two kinds of equations about f epsilon}). For equation (\ref{research equation}), our idea is to get the following kind of $L^{\infty}$ estimate
\begin{equation}\nonumber
\|f\|_{L^{\infty}((0, \infty)\times \Omega\times \mathbb{S}^{2})}\leq \frac{C}{\epsilon^{2}}(\|F\|_{L^{\infty}((0, \infty)\times \Omega\times \mathbb{S}^{2})}+\|h\|_{L^{\infty}(\Omega\times \mathbb{S}^{2})})
\end{equation}
by using the extra term $\epsilon^{2}f$. Then we can control the remainder $\mathcal{R}_{2}$ as follows
\begin{equation}\nonumber
\begin{aligned}
\|\mathcal{R}_{2}\|_{L^{\infty}((0, T)\times\Omega\times \mathbb{S}^{2})}\leq &\frac{C}{\epsilon^{2}}\|\mathcal{L}\mathcal{R}_{2}\|_{L^{\infty}((0, T)\times\Omega\times S^{2})}\\\leq&C\|B(\theta^{\epsilon})-B(\theta_{0})\|_{L^{\infty}((0, T)\times\Omega)}.
\end{aligned}
\end{equation}

On the other hand, in order to obtain the convergence $\lim_{\epsilon\rightarrow 0}\|\mathcal{R}_{2}\|_{L^{\infty}((0, T)\times\Omega\times \mathbb{S}^{2})}=0$, we need to prove that $\lim_{\epsilon\rightarrow 0}\|B(\theta^{\epsilon})-B(\theta_{0})\|_{L^{\infty}((0, T)\times\Omega\times \mathbb{S}^{2})}=0$. First, we prove that
\begin{equation}\nonumber
\lim_{\epsilon\rightarrow 0}\|\theta^{\epsilon}-\theta\|_{C^{0}([0, T], W_{2}^{5}(\Omega))}=0
\end{equation}
where $(\overline f, \rho, \vec u, \theta)$ is the solution of equations (\ref{equations about limit of overline f epsilon}) by employing the compactness analysis. Then we conclude $\overline f_{0}=\overline f,$ and $ \theta_{0}=\theta$ through the uniqueness of equations (\ref{equations about limit of overline f epsilon}).
So we have
\begin{equation}\nonumber
\lim_{\epsilon\rightarrow 0}\|\mathcal{R}_{2}\|_{L^{\infty}((0, T)\times\Omega\times \mathbb{S}^{2})}\leq \lim_{\epsilon\rightarrow 0}C(r)\|\theta^{\epsilon}-\theta_{0}\|_{C^{0}([0, T], W_{2}^{5}(\Omega))}=0.
\end{equation}
Finally, we obtain that
\begin{equation}\nonumber
\begin{aligned}
\|f^{\epsilon}-\overline f_{0}-f_{I, 0}\|_{L^{\infty}((0, T)\times\Omega\times \mathbb{S}^{2})}&\leq \|\sum_{k=1}^{2}\epsilon^{k}(f_{k}+f_{I, k})\|_{L^{\infty}((0, T)\times\Omega\times \mathbb{S}^{2})}+\|\mathcal{R}_{2}\|_{L^{\infty}((0, T)\times\Omega\times \mathbb{S}^{2})}\\&=O(\epsilon)+\|\mathcal{R}_{2}\|_{L^{\infty}((0, T)\times\Omega\times \mathbb{S}^{2})}
\end{aligned}
\end{equation}
which further implies (\ref{jielunbudengshi 1-jia}) by the fact that $\overline f=\overline f_{0}$.

In addition, we mention that one can take the P1 hypothesis to get the compressible Euler/Navier-Stokes-P1 approximation model when the distribution of photons is almost isotropic (see \cite{bu-10}). Recently, there are some investigations on the singular limits of such approximation model. For example, Jiang, Li and Xie \cite{bu-7} studied the nonrelativistic limit problem for the Navier-Stokes-Fourier-P1 approximation model as the reciprocal of the light speed tends to zero. Jiang, Ju and Liao \cite{bu-6} showed the diffusive limit of the compressible Euler-P1 approximation model arising in radiation hydrodynamics as the Mach number tends to zero.

Throughout this paper, $C$ denotes a certain positive constant and $C(\cdot)$ and $C_{1}(\cdot)$ are the positive constants depending on the quantity $"\cdot"$. We use $\hookrightarrow$ and $\hookrightarrow\hookrightarrow$ to denote the continuous embedding and compact embedding, respectively. $D_{x}$ and $D_{t}$ represent the operator $\partial_{x_{i}}, i=1, 2, 3$ and $\partial_{t}$, respectively. $W_{q}^{k}(\Omega)=W^{k, q}(\Omega)$ $(1\leq q\leq\infty)$ denotes the usual Lebesgue space on $\Omega$ with norm $\|\cdot\|_{W^{k}_{q}(\Omega)}$. We denote by $L^{p}(0, T; W^{k}_{q}(\Omega))$ $(1\leq p\leq\infty)$ the space of $L^{p}$ functions on $(0, T)$ with values in $W^{k}_{q}(\Omega)$, and $C^{0}(I, W^{k}_{q}(\Omega))$ standards for the space of continuous functions on the interval $I$ with values in $W_{q}^{k}(\Omega)$.  For $k\geq 0$ and $k\in\mathbb{Z}$, $C^{2k}([0, \infty)\times\overline\Omega)$ denotes the set satisfying $D_{x}^{l}D_{t}^{r}f\in C^{0}([0, \infty)\times\overline\Omega)$ where $0\leq l+2r\leq 2k$ while $C^{2k+1}([0, \infty)\times\overline\Omega)$ represents the set satisfying $D_{x}^{l}D_{t}^{r}f\in C^{0}([0, \infty)\times\overline\Omega)$ where $0\leq l+2r\leq 2k$ and $D_{x}^{2k+1}f\in C^{0}([0, \infty)\times\overline\Omega)$.

The paper is organized as follows. In Section \ref{san}, we discuss the well-posedness of neutron transport equation and prove the first part of Theorem \ref{definition weak solution-0}. In Section \ref{si}, we discuss the equations satisfied by the moments of $f^{\epsilon}$ and give the proof of the second part of Theorem \ref{definition weak solution-0}. In Section \ref{wu}, we show the details of the construction for the Hilbert expansion of system (\ref{research equations}). At the end of this section, we will give the proof of Theorem \ref{definition weak solution-1}.

\section{Well-posedness of system (\ref{research equations})}\label{san}
In this section, we consider the well-posedness of the system (\ref{research equations}). We first consider the following transport equation:
\begin{equation}\label{research equation}\left\{
\begin{split}
&{\epsilon}^{2}\frac{\partial f}{\partial t}+\epsilon \vec{w}\cdot \nabla_{x}f+f-\overline f+{\epsilon}^{2}f=F(t, \vec{x}, \vec{w})~~\mathrm{in} \ \ (0, \infty)\times \Omega\times \mathbb{S}^{2},\\
&f(0, \vec{x}, \vec{w})=h(\vec{x}, \vec{w})~~\mathrm{in} \ \ \Omega\times \mathbb{S}^{2}.\\
\end{split}\right.
\end{equation}

Define the $L^{2}$, $L^{\infty}$ and $C^{0}$ norms in $\Omega\times \mathbb{S}^{2}$ by
\begin{equation}\nonumber
\|f\|_{L^{2}(\Omega\times \mathbb{S}^{2})}=\Big(\int_{\Omega}\int_{\mathbb{S}^{2}}|f(\vec{x}, \vec{w})|^{2}\mathrm{d}\vec{w}\mathrm{d}\vec{x}\Big)^{\frac{1}{2}},
\end{equation}

\begin{equation}\nonumber
\|f\|_{L^{\infty}(\Omega\times \mathbb{S}^{2})}=\sup_{(\vec{x}, \vec{w})\in \Omega\times \mathbb{S}^{2}}|f(\vec{x}, \vec{w})|.\;\;
\end{equation}
Similar notations also applies to the space $[0, T]\times\Omega\times\mathbb{S}^{2}.$

\subsection{$L^{\infty}$ well-posedness}
We will show the existence and $L^{\infty}$ estimate of the transport equation (\ref{research equation}) in the following.
\vskip 2mm

\begin{lem}\label{main result}
Assume $F(t, \vec{x}, \vec{w})\in L^{\infty}((0, \infty)\times \Omega\times \mathbb{S}^{2}), h(\vec{x}, \vec{w})\in L^{\infty}(\Omega\times \mathbb{S}^{2})$ and $F, h\geq 0$. Then there exists a unique nonnegative solution $f(t, \vec{x}, \vec{w})\in L^{\infty}((0, \infty)\times \Omega\times \mathbb{S}^{2})$ of the transport equation (\ref{research equation}) which satisfies
\begin{equation}\label{estimate result}
\|f\|_{L^{\infty}((0, \infty)\times \Omega\times \mathbb{S}^{2})}\leq \frac{C}{\epsilon^{2}}(\|h\|_{L^{\infty}(\Omega\times \mathbb{S}^{2})}+\|F\|_{L^{\infty}((0, \infty)\times \Omega\times \mathbb{S}^{2})}).
\end{equation}

\end{lem}
\noindent {\bf Proof.} Since $f, F$ and $h$ are periodic, we can consider (\ref{research equation}) hold on $(0, \infty)\times\mathbb{R}^{3}\times \mathbb{S}^{2}$.

We construct an approximating sequence $\{f_{k}\}_{k=0}^{\infty}$ which are also periodic, where $f_{0}=0$ and

\begin{equation}\label{approximating sequence equation}\left\{
\begin{split}
&{\epsilon}^{2}\frac{\partial f_{k}}{\partial t}+\epsilon \vec w\cdot \nabla_{x}f_{k}+(1+\epsilon^{2})f_{k}
-\overline f_{k-1}=F(t, \vec{x}, \vec{w})~~\mathrm{in} \ \ (0, \infty)\times \mathbb{R}^{3}\times \mathbb{S}^{2},\\
&f_{k}(0, \vec{x}, \vec{w})=h(\vec{x}, \vec{w})~~\mathrm{in} \ \ \mathbb{R}^{3}\times \mathbb{S}^{2}.\\
\end{split}\right.
\end{equation}

The characteristics $(T(s), X(s), W(s))$ of equation (\ref{approximating sequence equation}) which goes through $(t, \vec{x}, \vec{w})$ is defined by
\begin{equation}\label{equation}\left\{
\begin{split}
&(T(0), X(0), W(0))=(t, \vec{x}, \vec{w}),\\
&\frac{dT(s)}{ds}=\epsilon^{2},\\
&\frac{dX(s)}{ds}=\epsilon W(s),\\
&\frac{dW(s)}{ds}=0,\\
\end{split}\right.
\end{equation}
which implies

\begin{equation}\nonumber\label{equation}\left\{
\begin{split}
&T(s)=t+\epsilon^{2}s,\\
&X(s)=\vec{x}+(\epsilon \vec{w})s,\\
&W(s)=\vec{w}.\\
\end{split}\right.
\end{equation}
Then,
\begin{equation}\label{similar-1}
\begin{split}
f_{k}(t, \vec{x}, \vec{w})=&h(\vec{x}-\frac{(\epsilon t\vec{w})}{\epsilon^{2}}, \vec{w})e^{-(1+\epsilon^{2})\frac{t}{\epsilon^{2}}} \\&+\int_{0}^{\frac{t}{\epsilon^{2}}}(\overline f_{k-1}+F)(\epsilon^{2}s, \vec{x}-\epsilon(\frac{t}{\epsilon^{2}}-s)\vec{w}, \vec{w})e^{-(1+\epsilon^{2})(\frac{t}{\epsilon^{2}}-s)}ds.
\end{split}
\end{equation}
Since $f_{0}=0$ and $h, g, F\geq 0$, we obtain $f_{1}\geq 0$. Then, we get $f_{k}\geq 0$ for any $k\geq 0$ iteratively. Set $v_{k}=f_{k}-f_{k-1}$ for $k\geq 1.$ Then, we have
\begin{equation}\nonumber
\begin{split}
v_{k+1}(t, \vec{x}, \vec{w})=\int_{0}^{\frac{t}{\epsilon^{2}}}\overline v_{k}(\epsilon^{2}s, \vec{x}-\epsilon(\frac{t}{\epsilon^{2}}-s)\vec{w}, \vec{w})e^{-(1+\epsilon^{2})(\frac{t}{\epsilon^{2}}-s)}ds.
\end{split}
\end{equation}
Since $\|\overline v_{k}\|_{L^{\infty}((0, \infty)\times\mathbb{R}^{3}\times \mathbb{S}^{2})}\leq \|v_{k}\|_{L^{\infty}((0, \infty)\times\mathbb{R}^{3}\times \mathbb{S}^{2})}$, we obtain
\begin{equation}\nonumber
\begin{split}
\|v_{k+1}\|_{L^{\infty}((0, \infty)\times\mathbb{R}^{3}\times \mathbb{S}^{2})}\leq\frac{1-e^{-(1+\epsilon^{2})\frac{t}{\epsilon^{2}}}}{1+\epsilon^{2}}\|v_{k}\|_{L^{\infty}((0, \infty)\times\mathbb{R}^{3}\times \mathbb{S}^{2})}.
\end{split}
\end{equation}
It follows that
\begin{equation}\nonumber
\|v_{k+1}\|_{L^{\infty}((0, \infty)\times\mathbb{R}^{3}\times \mathbb{S}^{2})}\leq\frac{1}{1+\epsilon^2}\|v_{k}\|_{L^{\infty}((0, \infty)\times\mathbb{R}^{3}\times \mathbb{S}^{2})},
\end{equation}
which implies $\{v_{k}\}$ is a contraction sequence. According to the fact that $v_{1}=f_{1}$, we have
\begin{equation}\nonumber
\|v_{k+1}\|_{L^{\infty}((0, \infty)\times\mathbb{R}^{3}\times \mathbb{S}^{2})}\leq(\frac{1}{1+\epsilon^2})^{k}\|f_{1}\|_{L^{\infty}((0, \infty)\times \mathbb{R}^{3}\times\mathbb{S}^{2})}
\end{equation}
for $k\geq 1.$ Therefore, $f_{k}$ converges strongly in $L^{\infty}((0, \infty)\times\Omega\times \mathbb{S}^{2})$ to a limit solution $f$ which satisfies
\begin{equation}\label{infinity estimate about limit solution}
\|f\|_{L^{\infty}((0, \infty)\times\mathbb{R}^{3}\times \mathbb{S}^{2})}\leq\sum_{k=1}^{\infty}\|v_{k}\|_{L^{\infty}((0, \infty)\times\mathbb{R}^{3}\times \mathbb{S}^{2})}\leq\frac{1+\epsilon^{2}}{\epsilon^{2}}\|f_{1}\|_{L^{\infty}((0, \infty)\times\mathbb{R}^{3}\times \mathbb{S}^{2})}.
\end{equation}
Since
\begin{equation}\nonumber
\begin{split}
f_{1}(t, \vec{x}, \vec{w})=h(\vec{x}-\frac{(\epsilon t\vec{w})}{\epsilon^{2}}, \vec{w})e^{-(1+\epsilon^{2})\frac{t}{\epsilon^{2}}} +\int_{0}^{\frac{t}{\epsilon^{2}}}F(\epsilon^{2}s, \vec{x}-\epsilon(\frac{t}{\epsilon^{2}}-s)\vec{w}, \vec{w})e^{-(1+\epsilon^{2})(\frac{t}{\epsilon^{2}}-s)}ds,
\end{split}
\end{equation}
we obtain
\begin{equation}\label{the first term about the sequence infinity estimate}
\|f_{1}\|_{L^{\infty}((0, \infty)\times\mathbb{R}^{3}\times \mathbb{S}^{2})}\leq \|h\|_{L^{\infty}(\mathbb{R}^{3}\times \mathbb{S}^{2})}+\|F\|_{L^{\infty}((0, \infty)\times \mathbb{R}^{3}\times \mathbb{S}^{2})}.
\end{equation}
It is obvious to obtain that
\begin{equation}\label{equality periodicity}
\begin{aligned}
&\|F\|_{L^{\infty}((0, \infty)\times\mathbb{R}^{3}\times \mathbb{S}^{2})}=\|F\|_{L^{\infty}((0, \infty)\times\mathbb{T}^{3}\times \mathbb{S}^{2})},
\|h\|_{L^{\infty}(\mathbb{R}^{3}\times\mathbb{S}^{2})}=\|h\|_{L^{\infty}(\mathbb{T}^{3}\times \mathbb{S}^{2})}
\\&\|f_{k}\|_{L^{\infty}((0, \infty)\times\mathbb{R}^{3}\times \mathbb{S}^{2})}=\|f_{k}\|_{L^{\infty}((0, \infty)\times\mathbb{T}^{3}\times \mathbb{S}^{2})},
\|f\|_{L^{\infty}((0, \infty)\times\mathbb{R}^{3}\times \mathbb{S}^{2})}=\|f\|_{L^{\infty}((0, \infty)\times\mathbb{T}^{3}\times \mathbb{S}^{2})},
\end{aligned}
\end{equation}
for $k\geq 0$, by the periodicity of $F, h, f_{k}$ and $f$.

According to (\ref{infinity estimate about limit solution}) and (\ref{the first term about the sequence infinity estimate}), we can easily derive the existence and estimates. Furthermore, we can get the uniqueness of equation (\ref{research equation}) according to the $L^{\infty}$ energy estimate.  \hfill$\Box$

\begin{lem}\label{main result-7}
For any $0<T<\infty$, if $D_{x}^{l}F(t, \vec{x}, \vec{w})\in C^{0}([0, T]\times \overline\Omega\times \mathbb{S}^{2}), D_{x}^{l}h(\vec{x}, \vec{w})\in C^{0}(\overline\Omega\times \mathbb{S}^{2})$ where $0\leq l\leq m$, we can obtain that the solution $f$ in Lemma \ref{main result} also satisfies that $D_{x}^{l}f\in C^{0}([0, T]\times \overline\Omega\times \mathbb{S}^{2})$ where $0\leq l\leq m$.
\end{lem}
\noindent {\bf Proof.} We can apply the operator $D_{x}^{l}$ to (\ref{similar-1}) by the fact that $D_{x}^{l}F(t, \vec{x}, \vec{w})\in C^{0}([0, T]\times \overline\Omega\times \mathbb{S}^{2}), D_{x}^{l}h(\vec{x}, \vec{w})\in C^{0}(\overline\Omega\times \mathbb{S}^{2})$ where $0\leq l\leq m$.

We can finally obtain that
\begin{equation}\nonumber
\|D_{x}^{l}v_{k+1}\|_{L^{\infty}((0, \infty)\times\mathbb{R}^{3}\times \mathbb{S}^{2})}\leq(\frac{1}{1+\epsilon^2})^{k}\|D_{x}^{l}f_{1}\|_{L^{\infty}((0, \infty)\times\mathbb{R}^{3}\times \mathbb{S}^{2})},
\end{equation}
which implies $\{D_{x}^{l}v_{k}\}$ is a contraction sequence. We can conclude that $D_{x}^{l}f_{k}$ converges strongly in $C^{0}([0, T]\times \mathbb{R}^{3}\times \mathbb{S}^{2})$ to a limit solution $D_{x}^{l}f$ which satisfies
\begin{equation}\label{infinity estimate about limit solution}
\|D_{x}^{l}f\|_{L^{\infty}((0, \infty)\times\mathbb{R}^{3}\times \mathbb{S}^{2})}\leq\sum_{k=1}^{\infty}\|D_{x}^{l}v_{k}\|_{L^{\infty}((0, \infty)\times\mathbb{R}^{3}\times \mathbb{S}^{2})}\leq\frac{1+\epsilon^{2}}{\epsilon^{2}}\|D_{x}^{l}f_{1}\|_{L^{\infty}((0, \infty)\times\mathbb{R}^{3}\times \mathbb{S}^{2})}.
\end{equation}
Since
\begin{equation}\nonumber
\begin{split}
D_{x}^{l}f_{1}(t, \vec{x}, \vec{w})=D_{x}^{l}h(\vec{x}-\frac{(\epsilon t\vec{w})}{\epsilon^{2}}, \vec{w})e^{-(1+\epsilon^{2})\frac{t}{\epsilon^{2}}} +\int_{0}^{\frac{t}{\epsilon^{2}}}D_{x}^{l}F(\epsilon^{2}s, \vec{x}-\epsilon(\frac{t}{\epsilon^{2}}-s)\vec{w}, \vec{w})e^{-(1+\epsilon^{2})(\frac{t}{\epsilon^{2}}-s)}ds,
\end{split}
\end{equation}
we obtain
\begin{equation}\label{the first term about the sequence infinity estimate}
\|D_{x}^{l}f_{1}\|_{L^{\infty}((0, \infty)\times\mathbb{R}^{3}\times \mathbb{S}^{2})}\leq \|D_{x}^{l}h\|_{L^{\infty}(\mathbb{R}^{3}\times \mathbb{S}^{2})}+\|D_{x}^{l}F\|_{L^{\infty}((0, \infty)\times\mathbb{R}^{3}\times \mathbb{S}^{2})}.
\end{equation}
Therefore, we have
\begin{equation}\label{estimate result-1-1}
\|D_{x}^{l}f\|_{L^{\infty}((0, \infty)\times\mathbb{R}^{3}\times \mathbb{S}^{2})}\leq \frac{C}{\epsilon^{2}}(\|D_{x}^{l}h\|_{L^{\infty}(\mathbb{R}^{3}\times \mathbb{S}^{2})}+\|D_{x}^{l}F\|_{L^{\infty}((0, \infty)\times\mathbb{R}^{3}\times \mathbb{S}^{2})}).
\end{equation}
which implies $D_{x}^{l}f\in C^{0}([0, T]\times\mathbb{R}^{3}\times \mathbb{S}^{2})$ where $0\leq l\leq m$. According to (\ref{equality periodicity}), we have $D_{x}^{l}f\in C^{0}([0, T]\times\overline\Omega\times \mathbb{S}^{2})$ where $0\leq l\leq m$. \hfill$\Box$

\begin{lem}\label{main result-7-bu}
For any $0<T<\infty$, if $F(t, \vec{x}, \vec{w})\in L^{2}(0, T; L^{2}(\mathbb{S}^{2}; W_{2}^{m}(\Omega))), h(\vec{x}, \vec{w})\in L^{2}(\mathbb{S}^{2}; W_{2}^{m}(\Omega))$, there exists a unique nonnegative solution $f(t, \vec{x}, \vec{w})\in L^{\infty}(0, T; L^{2}(\mathbb{S}^{2}; \\W_{2}^{m}(\Omega)))$ of the transport equation (\ref{research equation}) which satisfies
\begin{equation}\label{L2 energy estimate}
\begin{aligned}
\|f\|_{L^{\infty}(0, T; L^{2}(\mathbb{S}^{2}; W_{2}^{m}(\Omega)))}\leq C(\|h\|_{L^{2}(\mathbb{S}^{2}; W_{2}^{m}(\Omega))}+\frac{1}{\epsilon^{2}}\|F\|_{L^{2}(0, T; L^{2}(\mathbb{S}^{2}; W_{2}^{m}(\Omega)))})
\end{aligned}
\end{equation}
\end{lem}
\noindent {\bf Proof.} We can construct sequences $\{F_{k}\}_{k=1}^{\infty}\subset C^{m}([0, T]\times\overline\Omega\times\mathbb{S}^{2})$ and $\{h_{k}\}_{k=1}^{\infty}\subset C^{m}(\overline\Omega\times\mathbb{S}^{2})$ satisfying $\lim_{k\rightarrow\infty}\|F_{k}-F\|_{L^{2}(0, T; L^{2}(\mathbb{S}^{2}; W_{2}^{m}(\Omega)))}=0$ and $\lim_{k\rightarrow\infty}\|h_{k}-h\|_{L^{2}(\mathbb{S}^{2}; W_{2}^{m}(\Omega))}\\=0$. Then, we consider the following equations
\begin{equation}\label{approximating sequence equation-bu}\left\{
\begin{split}
&{\epsilon}^{2}\frac{\partial f_{k}}{\partial t}+\epsilon \vec w\cdot \nabla_{x}f_{k}+f_{k}
-\overline f_{k}+\epsilon^{2}f_{k}=F_{k}(t, \vec{x}, \vec{w}),~~\mathrm{in} \ \ (0, \infty)\times \Omega\times \mathbb{S}^{2},\\
&f_{k}(0, \vec{x}, \vec{w})=h_{k}(\vec{x}, \vec{w}),~~\mathrm{in} \ \ \Omega\times \mathbb{S}^{2}.\\
\end{split}\right.
\end{equation}
According to Lemma \ref{main result-7}, we have $D_{x}^{l}f_{k}\in C^{0}([0, T]\times\overline\Omega\times\mathbb{S}^{2})$ where $0\leq l\leq m$. Applying the operator on both sides of $(\ref{approximating sequence equation-bu})_{1}$, multiplying $D_{x}^{l}f_{k}$ on both sides of $(\ref{approximating sequence equation-bu})_{1}$ and integrating on $(0, t)\times\Omega\times\mathbb{S}^{2}$ for $t\in [0, T]$, then, we get
\begin{equation}\nonumber
\begin{aligned}
&\int_{\Omega}\int_{\mathbb{S}^{2}}|D_{x}^{l}f_{k}|^{2}(t)\mathrm{d}\vec x\mathrm{d}\vec w+\int_{0}^{t}\int_{\Omega}\int_{\mathbb{S}^{2}}\frac{|D_{x}^{l}(f_{k}-\overline f_{k})|^{2}}{\epsilon^{2}}\mathrm{d} t\mathrm{d} \vec x\mathrm{d}\vec w+\int_{0}^{t}\int_{\Omega}\int_{\mathbb{S}^{2}}|D_{x}^{l}f_{k}|^{2}\\=&\int_{\Omega}\int_{\mathbb{S}^{2}}|D_{x}^{l}h_{k}|^{2}\mathrm{d}\vec x\mathrm{d} \vec w
+\frac{1}{\epsilon^{2}}\int_{0}^{t}\int_{\Omega}\int_{\mathbb{S}^{2}}D_{x}^{l}F_{k}D_{x}^{l}f_{k}\mathrm{d}t\mathrm{d}\vec x\mathrm{d}\vec w
\\\leq& \int_{\Omega}\int_{\mathbb{S}^{2}}|D_{x}^{l}h_{k}|^{2}\mathrm{d}\vec x\mathrm{d} \vec w
+\frac{2}{\epsilon^{4}}\int_{0}^{t}\int_{\Omega}\int_{\mathbb{S}^{2}}|D_{x}^{l}F_{k}|^{2}\mathrm{d}t\mathrm{d}\vec x\mathrm{d}\vec w+\frac{1}{2}\int_{0}^{t}\int_{\Omega}\int_{\mathbb{S}^{2}}|D_{x}^{l}f_{k}|^{2}\mathrm{d}t\mathrm{d}\vec x\mathrm{d}\vec w.
\end{aligned}
\end{equation}
So,
\begin{equation}\nonumber
\begin{aligned}
\sup_{t\in[0, T]}\int_{\Omega}\int_{\mathbb{S}^{2}}|D_{x}^{l}f_{k}|^{2}(t)\mathrm{d}\vec x\mathrm{d}\vec w\leq \int_{\Omega}\int_{\mathbb{S}^{2}}|D_{x}^{l}h_{k}|^{2}\mathrm{d}\vec x\mathrm{d} \vec w
+\frac{2}{\epsilon^{4}}\int_{0}^{T}\int_{\Omega}\int_{\mathbb{S}^{2}}|D_{x}^{l}F_{k}|^{2}\mathrm{d}t\mathrm{d}\vec x\mathrm{d}\vec w.
\end{aligned}
\end{equation}
Letting $k\rightarrow\infty$, we have
\begin{equation}\nonumber
\begin{aligned}
\sup_{t\in[0, T]}\int_{\Omega}\int_{\mathbb{S}^{2}}|D_{x}^{l}f|^{2}(t)\mathrm{d}\vec x\mathrm{d}\vec w\leq \int_{\Omega}\int_{\mathbb{S}^{2}}|D_{x}^{l}h|^{2}\mathrm{d}\vec x\mathrm{d} \vec w
+\frac{2}{\epsilon^{4}}\int_{0}^{T}\int_{\Omega}\int_{\mathbb{S}^{2}}|D_{x}^{l}F|^{2}\mathrm{d}t\mathrm{d}\vec x\mathrm{d}\vec w,
\end{aligned}
\end{equation}
which implies (\ref{L2 energy estimate}).

Note that $D_{x}^{l}f_{k}\rightarrow D_{x}^{l}f$ in $L^{\infty}(0, T; L^{2}(\mathbb{S}^{2}\times\Omega))$ in the weak-$\ast$ sense up to a subsequence. The uniqueness of equation (\ref{research equation}) can be ensured by (\ref{L2 energy estimate}).\hfill$\Box$

\subsection{The proof of the first part of Theorem \ref{definition weak solution-0}} We can write (\ref{research equations}) in the following form
\begin{equation}\label{research equations-3}\left\{
\begin{split}
&\partial_{t}\rho^{\epsilon}+\vec{u}^{\epsilon}\cdot\nabla_{x}\rho^{\epsilon}+\rho^{\epsilon}\mathrm{div}_{x}\vec{u}^{\epsilon}=0,~~\mathrm{in} \ \ (0,T)\times \Omega,\\
&\partial_{t}\vec{u}^{\epsilon}+\vec{u}^{\epsilon}\cdot\nabla_{x}\vec{u}^{\epsilon}+\nabla_{x}\theta^{\epsilon}
+\frac{\theta^{\epsilon}}{\rho^{\epsilon}}\nabla_{x}\rho^{\epsilon}=\frac{1}{\rho^{\epsilon}}\Big\langle\Big(\frac{1}{\epsilon}+\epsilon\Big)\vec{w}(f^{\epsilon}-\overline f^{\epsilon})\Big\rangle,~~\mathrm{in} \ \ (0,T)\times \Omega,\\
&\partial_{t}\theta^{\epsilon}+\vec{u}^{\epsilon}\cdot\nabla_{x}\theta^{\epsilon}+\theta^{\epsilon}\mathrm{div}_{x}\vec{u}^{\epsilon}
=\frac{1}{\rho^{\epsilon}}(\overline f^{\epsilon}-B(\theta^{\epsilon}))\\&\hspace{5.5cm}-\frac{1}{\rho^{\epsilon}}\Big(\frac{1}{\epsilon}+\epsilon\Big)\langle\vec{w}(f^{\epsilon}-\overline f^{\epsilon})\rangle\cdot\vec{u}^{\epsilon},~~\mathrm{in} \ \ (0,T)\times \Omega, \\
&\partial_{t}f^{\epsilon}+\frac{1}{\epsilon}\vec w\cdot \nabla_{x}f^{\epsilon}+\frac{1}{\epsilon^{2}}(f^{\epsilon}-\overline f^{\epsilon})+f^{\epsilon}=B(\theta^{\epsilon}),~~\mathrm{in} \ \ (0,T)\times \Omega\times \mathbb{S}^{2},\\
&\rho^{\epsilon}(0, \vec{x})=\rho^{0}(\vec x),\;\vec{u}^{\epsilon}(0, \vec{x})=\vec{u}^{0}(\vec x), \;\theta^{\epsilon}(0, \vec{x})=\theta^{0}(\vec{x}),~~\mathrm{in}  \ \ \Omega,\\
&f^{\epsilon}(0, \vec{x}, \vec{w})=h(\vec{x}, \vec{w}),~~\mathrm{in} \ \ \Omega\times \mathbb{S}^{2}.\\
\end{split}\right.
\end{equation}
Denote the vector and matrix
\begin{equation}\nonumber
V^{\epsilon}=(\rho^{\epsilon}, u_{1}^{\epsilon}, u_{2}^{\epsilon}, u_{3}^{\epsilon}, \theta^{\epsilon})^{t}, \tilde{A}_{j}(V^{\epsilon})=\{\tilde{a}_{mn}\}_{5\times5}
\end{equation}
where $\tilde{a}_{ii}=u_{j}^{\epsilon}$, $\tilde{a}_{1(j+1)}=\rho^{\epsilon}$, $\tilde{a}_{(j+1)1}=\frac{\theta^{\epsilon}}{\rho^{\epsilon}}$, $\tilde{a}_{(j+1)5}=1$, $\tilde{a}_{5(j+1)}=\theta^{\epsilon}$ for $j=1, 2, 3,$ and the rest elements of $\tilde{a}_{ij}$ are set to 0. Define the vector $H(V^{\epsilon}, f^{\epsilon})$ by
\begin{equation}\nonumber
H(V^{\epsilon}, f^{\epsilon})=(g_{0}^{\epsilon}, g_{1}^{\epsilon}, g_{2}^{\epsilon}, g_{3}^{\epsilon}, g_{4}^{\epsilon})^{t}
\end{equation}
with $g_{0}^{\epsilon}=0$, and for $j=1, 2, 3,$
\begin{equation}\nonumber
g_{j}^{\epsilon}=\frac{1}{\rho^{\epsilon}}\Big\langle\Big(\frac{1}{\epsilon}+\epsilon\Big)w_{j}(f^{\epsilon}-\overline f^{\epsilon})\Big\rangle,
\end{equation}
\begin{equation}\nonumber
g_{4}^{\epsilon}=\frac{1}{\rho^{\epsilon}}(\overline f^{\epsilon}-B(\theta^{\epsilon}))-\frac{1}{\rho^{\epsilon}}\Big(\frac{1}{\epsilon}+\epsilon\Big)\langle\vec{w}(f^{\epsilon}-\overline f^{\epsilon})\rangle\cdot\vec{u}^{\epsilon}.
\end{equation}
Therefore, we can rewrite $(\ref{research equations-3})_{1\sim 3}$ as
\begin{equation}\label{research equations-3-1}
\frac{\partial V^{\epsilon}}{\partial t}+\sum_{j=1}^{3}\tilde{A}_{j}(V^{\epsilon})\frac{\partial V^{\epsilon}}{\partial x_{j}}=H(V^{\epsilon}, f^{\epsilon}).
\end{equation}
We shall study the Cauchy problem for $(\ref{research equations-3})_{4}$ and (\ref{research equations-3-1}) together with the initial data
\begin{equation}\label{research equations-3-3}
V^{\epsilon}(0, \vec x)=V^{0}(\vec x),~~~~~~f^{\epsilon}(0, \vec x, \vec w)=h(\vec x, \vec w).
\end{equation}
First, we symmetrize $(\ref{research equations-3-1})$ by multiplying it with the symmetrizing matrix defined by
\begin{equation}\label{xishujuzhen-1}
A_{0}(V^{\epsilon})
=
\left[
\begin{array}{cccc}
(\rho^{\epsilon})^{-1}&0&0\\
0&\frac{\rho^{\epsilon}}{\theta^{\epsilon}}\mathbb{I}_{3\times3}&0\\
0&0&\frac{\rho^{\epsilon}}{(\theta^{\epsilon})^{2}}
\end{array}
\right].
\end{equation}
Therefore, we shall prove the first part of Theorem \ref{definition weak solution-0} for the quasilinear symmetric system:
\begin{equation}\label{research equations-3-2}
A_{0}(V^{\epsilon})\frac{\partial V^{\epsilon}}{\partial t}+\sum_{j=1}^{3}A_{j}(V^{\epsilon})\frac{\partial V^{\epsilon}}{\partial x_{j}}=F(V^{\epsilon}, f^{\epsilon})
\end{equation}
coupled to $(\ref{research equations-3})_{4}$, where $A_{j}(V^{\epsilon}):=A_{0}(V^{\epsilon})\tilde{A}_{j}(V^{\epsilon})$ is symmetric, and
\begin{equation}\nonumber
F(V^{\epsilon}, f^{\epsilon})=A_{0}(V^{\epsilon})H(V^{\epsilon}, f^{\epsilon})=(l_{0}^{\epsilon}, l_{1}^{\epsilon}, l_{2}^{\epsilon}, l_{3}^{\epsilon}, l_{4}^{\epsilon})^{t}
\end{equation}
with $l_{0}=0$, and for $j=1, 2, 3$,
\begin{equation}\nonumber
l_{j}^{\epsilon}=\frac{1}{\theta^{\epsilon}}\Big\langle\Big(\frac{1}{\epsilon}+\epsilon\Big)w_{j}(f^{\epsilon}-\overline f^{\epsilon})\Big\rangle,
\end{equation}
\begin{equation}\nonumber
l_{4}^{\epsilon}=\frac{1}{(\theta^{\epsilon})^{2}}(\overline f^{\epsilon}-B(\theta^{\epsilon}))-\frac{1}{(\theta^{\epsilon})^{2}}\Big(\frac{1}{\epsilon}+\epsilon\Big)\langle\vec{w}(f^{\epsilon}-\overline f^{\epsilon})\rangle\cdot\vec{u}^{\epsilon}.
\end{equation}
In the sequel, we construct a solution to (\ref{research equations-3-2}), $(\ref{research equations-3})_{4}$ and $(\ref{research equations-3-3})$. For $k\geq 0$ and $k\in\mathbb{Z}$, we define $V_{k+1}^{\epsilon}(t, \vec x)$ and $f_{k+1}^{\epsilon}$ inductively as the solution of the following linearized system:
\begin{equation}\label{research equations-4}\left\{
\begin{split}
&A_{0}(V_{k})\frac{\partial V_{k+1}}{\partial t}+\sum_{j=1}^{3}A_{j}(V_{k})\frac{\partial V_{k+1}}{\partial x_{j}}=F(V_{k}, f_{k})~~\mathrm{in} \ \ (0,T)\times \Omega,\\
&\partial_{t}f_{k+1}+\frac{1}{\epsilon}\vec w\cdot \nabla_{x}f_{k+1}+\frac{1}{\epsilon^{2}}(f_{k+1}-\overline f_{k+1})+f_{k+1}=B(\theta_{k})~~\mathrm{in} \ \ (0,T)\times \Omega\times \mathbb{S}^{2},\\
&V_{k+1}(0, \vec x)=V^{0}(\vec x)~~\mathrm{in}  \ \ \Omega,\\
&f_{k+1}(0, \vec{x}, \vec{w})=h(\vec{x}, \vec{w})~~\mathrm{in} \ \ \Omega\times \mathbb{S}^{2},\\
\end{split}\right.
\end{equation}
where we drop the superscripts of $V_{k}^{\epsilon}, V_{k+1}^{\epsilon}$, $f_{k}^{\epsilon}$ and $f_{k+1}^{\epsilon}$ for simplicity.

Define the map
\begin{equation}\nonumber
S: (f_{k}, V_{k})\in G\rightarrow (f_{k+1}, V_{k+1}):=S(f_{k}, V_{k})
\end{equation}
where $(f_{k+1}, V_{k+1})$ satisfies (\ref{research equations-4}).

First, we note that the Theorem I in \cite{cauchy-problem} can be shown in a similar manner in $\mathbb{T}^{3}$.
Choosing $V_{0}=V^{0}$ and $f_{0}=0$ for some $T\in(0, \infty),$ we get $V_{1}\in C^{0}([0, T]; W_{2}^{6}(\Omega))$ by the Theorem I in \cite{cauchy-problem}. Then we have $B(\theta_{1})\in C^{0}([0, T]; W_{2}^{6}(\Omega))$. Thus, we derive taht $f_{1}\in L^{\infty}((0, T); L^{2}(\mathbb{S}^{2}; W_{2}^{6}(\Omega)))\cap C^{0}([0, T]; L^{2}(S^{2}; W_{2}^{5}(\Omega)))$ according to Lemma \ref{main result-7-bu} combing with Aubin-Lions-Simon lemma. So we have $F(V_{1}, f_{1})\in L^{\infty}((0, T); W_{2}^{6}(\Omega))\cap C^{0}([0, T]; W_{2}^{5}(\Omega))$. It follows that $V_{2}\in C^{0}([0, T]; W_{2}^{6}(\Omega))$. Similarly, we can derive that $f_{2}\in L^{\infty}((0, T); L^{2}(\mathbb{S}^{2}; W_{2}^{6}(\Omega)))\cap C^{0}([0, T]; L^{2}(S^{2}; W_{2}^{5}(\Omega)))$. Finally, we obtain that $V_{k+1}, B(\theta_{k+1})\in C^{0}([0, T]; W_{2}^{6}(\Omega))$ and $f_{k+1}\in L^{\infty}((0, T); L^{2}(\mathbb{S}^{2}; W_{2}^{6}(\Omega)))\cap C^{0}([0, T]; L^{2}(S^{2};\\ W_{2}^{5}(\Omega))),$ for any $k\geq 0$ iteratively although these estimates and the existence time $T$  depend on $\epsilon$ and $k$. Next, we will show uniform estimates of solutions to the system (\ref{research equations-4}).

Set
\begin{equation}\nonumber
r_{k}=\|f_{k}\|_{L^{\infty}(0, T; L^{2}(\mathbb{S}^{2}; W_{2}^{6}(\Omega)))}+\|V_{k}\|_{L^{\infty}(0, T; W_{2}^{6}(\Omega))}.
\end{equation}
It is easy to obtain that there exists small $T=T(k, \epsilon)$ such that
\begin{equation}\nonumber
\theta_{k}, \rho_{k}\geq M_{k}>0.
\end{equation}
We will prove that there exist $T$ independent of $k$ and $\epsilon$ such that
\begin{equation}\nonumber
\theta_{k}, \rho_{k}\geq M_{k}\geq M>0
\end{equation}
and
\begin{equation}\nonumber
r_{k}\leq r
\end{equation}
where $M$ and $r$ are independent of $k$ and $\epsilon$ for any $k\geq 0$. Note that $r_{k}, M_{k}, f_{k}, V_{k}$ will disappear if $k<0$.

Applying the operator $D_{x}^{\gamma}$ to the equation $(\ref{research equations-4})_{2}$, multiplying the resulting equation with $2(D_{x}^{\gamma}f_{k+1})$ in $L^{2}(\Omega\times \mathbb{S}^{2})$, we have
\begin{equation}\label{Eq 1}
\begin{split}
&\frac{\mathrm{d}}{\mathrm{d}t}\int_{\Omega}\int_{\mathbb{S}^{2}}|D_{x}^{\gamma}f_{k+1}|^{2}\mathrm{d}\vec{w}\mathrm{d}\vec{x}
+\frac{2}{\epsilon^{2}}\int_{\Omega}\int_{\mathbb{S}^{2}}|D_{x}^{\gamma}f_{k+1}-D_{x}^{\gamma}\overline f_{k+1}|^{2}\mathrm{d}\vec{x}\mathrm{d}\vec{w}
+2\int_{\Omega}\int_{\mathbb{S}^{2}}(D_{x}^{\gamma}f_{k+1})^{2}\mathrm{d}\vec{w}\mathrm{d}\vec{x}
\\=&2\int_{\Omega}\int_{\mathbb{S}^{2}}D_{x}^{\gamma}f_{k+1}D_{x}^{\gamma}(B(\theta_{k}))\mathrm{d}\vec{x}\mathrm{d}\vec{w}
\end{split}
\end{equation}
where $0\leq\gamma\leq 6.$
Integrating (\ref{Eq 1}) on $(0, t)$ and then taking the sup over $t\in(0, T)$, we have
\begin{equation}\label{gronwall's inequality 1-1}
\begin{split}
&\|f_{k+1}\|_{L^{\infty}(0, T; L^{2}(\mathbb{S}^{2}; W_{2}^{6}(\Omega)))}^{2}+\|(f_{k+1}-\overline f_{k+1})/\epsilon\|_{L^{2}(0, T; L^{2}(\mathbb{S}^{2}; W_{2}^{6}(\Omega)))}^{2}
\\\leq&C\|\theta_{k}\|_{L^{\infty}(0, T; W_{2}^{6}(\Omega))}^{8}T+\|h\|_{L^{2}(\mathbb{S}^{2}; W_{2}^{6}(\Omega))}^{2}
\\\leq&Cr_{k}^{8}T+\|h\|_{L^{2}(\mathbb{S}^{2}; W_{2}^{6}(\Omega))}^{2}.
\end{split}
\end{equation}
Set $U_{k+1}=V_{k+1}-V^{0}$. Then,
\begin{equation}\label{research equations-4-1-0}
\begin{split}
A_{0}(V_{k})\frac{\partial U_{k+1}}{\partial t}+\sum_{j=1}^{3}A_{j}(V_{k})\frac{\partial U_{k+1}}{\partial x_{j}}=F(V_{k}, f_{k})-\sum_{j=1}^{3}A_{j}(V_{k})\frac{\partial V^{0}}{\partial x_{j}},
\end{split}
\end{equation}
where $U_{k+1}(0, \vec x)=0.$ Moreover, applying the operator $A_{0}(V_{k})D_{x}^{\gamma}A_{0}^{-1}(V_{k})$ to the equation $(\ref{research equations-4-1-0})$, we obtain
\begin{equation}\label{research equations-4-1}
\begin{split}
A_{0}(V_{k})\frac{\partial D_{x}^{\gamma}U_{k+1}}{\partial t}+\sum_{j=1}^{3}A_{j}(V_{k})\frac{\partial D_{x}^{\gamma}U_{k+1}}{\partial x_{j}}=&A_{0}(V_{k})D_{x}^{\gamma}A_{0}^{-1}(V_{k})F(V_{k}, f_{k})+F_{\gamma}\\&-\sum_{j=1}^{3}A_{0}(V_{k})D_{x}^{\gamma}A_{0}^{-1}(V_{k})A_{j}(V_{k})\frac{\partial V^{0}}{\partial x_{j}},
\end{split}
\end{equation}
where $F_{\gamma}=-\sum_{j=1}^{3}\sum_{0\leq\beta\leq \gamma-1}A_{0}(V_{k})D_{x}^{\gamma-\beta}(A_{0}^{-1}(V_{k})A_{j}(V_{k}))\frac{\partial D_{x}^{\beta}U_{k+1}}{\partial x_{j}}$.

Define
\begin{equation}\nonumber
E(t):=\frac{1}{2}\int_{\Omega}(A_{0}(V_{k})D_{x}^{\gamma}U_{k+1})\cdot D_{x}^{\gamma}U_{k+1}\mathrm{d}\vec x.
\end{equation}
So,
\begin{equation}\label{energy equality-2}
\begin{aligned}
\frac{\mathrm{d}E(t)}{\mathrm{d}t}=&\int_{\Omega}\Big(\frac{1}{2}\partial_{t}A_{0}(V_{k})
D_{x}^{\gamma}U_{k+1}\cdot D_{x}^{\gamma}U_{k+1}+\frac{1}{2}\sum_{j=1}^{3}\frac{\partial A_{j}}{\partial x_{j}}D_{x}^{\gamma}U_{k+1}\cdot D_{x}^{\gamma}U_{k+1}\\&-\sum_{j=1}^{3}A_{0}(V_{k})D_{x}^{\gamma}A_{0}^{-1}(V_{k})A_{j}(V_{k})\frac{\partial V^{0}}{\partial x_{j}}\cdot D_{x}^{\gamma}U_{k+1}\\&+A_{0}(V_{k})D_{x}^{\gamma}(A_{0}^{-1}(V_{k})F(V_{k}, f_{k}))\cdot D_{x}^{\gamma}U_{k+1}+F_{\gamma}\cdot D_{x}^{\gamma}U_{k+1}\Big)\mathrm{d}\vec x.
\end{aligned}
\end{equation}
According to (\ref{xishujuzhen-1}), we have
\begin{equation}\nonumber
\partial_{t}A_{0}(V_{k})
=
\left[
\begin{array}{cccc}
-\frac{\partial_{t}\rho_{k}}{\rho_{k}^{2}}&0&0\\
0&\Big(\frac{\partial_{t}\rho_{k}}{\theta_{k}}-\frac{\rho_{k}\partial_{t}\theta_{k}}{\theta_{k}^{2}}\Big)\mathbb{I}_{3\times3}&0\\
0&0&\frac{\partial_{t}\rho_{k}}{(\theta_{k})^{2}}-\frac{2\rho_{k}\partial_{t}\theta_{k}}{\theta_{k}^{3}}
\end{array}
\right].
\end{equation}
For $k=0$, $\partial_{t}A_{0}(V_{0})=0$.
Since $\partial_{t}\rho_{k}=-\vec{u}_{k-1}\cdot\nabla_{x}\rho_{k}-\rho_{k-1}\mathrm{div}_{x}\vec{u}_{k}$ and $\partial_{t}\theta_{k}=-\theta_{k-1}\mathrm{div}_{x}\vec{u}_{k}-\vec{u}_{k-1}\cdot\nabla_{x}\theta_{k}
+\frac{1}{\rho_{k-1}}(\overline f_{k-1}-B(\theta_{k-1}))-\frac{1}{\rho_{k-1}}(\frac{1}{\epsilon}+\epsilon)\langle\vec{w}(f_{k-1}-\overline f_{k-1})\rangle\cdot\vec{u}_{k-1}$ for $k\geq 1$, then, we obtain
\begin{equation}\nonumber
\begin{aligned}
&\|\partial_{t}A_{0}(V_{k})\|_{L^{2}(0, T; L^{\infty}(\Omega))}\\\leq& \Big\|-\frac{\partial_{t}\rho_{k}}{\rho_{k}^{2}}\Big\|_{L^{2}(0, T; L^{\infty}(\Omega))}+\Big\|\frac{\partial_{t}\rho_{k}}{\theta_{k}}-\frac{\rho_{k}\partial_{t}\theta_{k}}{\theta_{k}^{2}}\Big\|_{L^{2}(0, T; L^{\infty}(\Omega))}+\Big\|\frac{\partial_{t}\rho_{k}}{(\theta_{k})^{2}}-\frac{2\rho_{k}\partial_{t}\theta_{k}}{\theta_{k}^{3}}\Big\|_{L^{2}(0, T; L^{\infty}(\Omega))}
\\\leq& C\Big(\frac{1}{M_{k}^{3}}+\frac{1}{M_{k}^{2}}+\frac{1}{M_{k}}+1\Big)(\|\partial_{t}\rho_{k}\|_{L^{2}(0, T; L^{\infty}(\Omega))}+\|\rho_{k}\partial_{t}\theta_{k}\|_{L^{2}(0, T; L^{\infty}(\Omega))})
\\\leq& C\Big(\frac{1}{M_{k}^{3}}+1\Big)(1+r_{k})(\|-\vec{u}_{k-1}\cdot\nabla_{x}\rho_{k}-\rho_{k-1}\mathrm{div}_{x}\vec{u}_{k}\|_{L^{\infty}((0, T)\times\Omega)}\sqrt{T}\\&+\|-\theta_{k-1}\mathrm{div}_{x}\vec{u}_{k}-\vec{u}_{k-1}\cdot\nabla_{x}\theta_{k}
+\frac{1}{\rho_{k-1}}(\overline f_{k-1}-B(\theta_{k-1}))\|_{L^{\infty}((0, T)\times\Omega)}\sqrt{T}\\&+\|-\frac{1}{\rho_{k-1}}\Big(\frac{1}{\epsilon}+\epsilon\Big)\langle\vec{w}(f_{k-1}-\overline f_{k-1})\rangle\cdot\vec{u}_{k-1}\|_{L^{2}(0, T; L^{\infty}(\Omega))})
\\\leq& C\Big(\frac{1}{M_{k}^{3}}+1\Big)(1+r_{k})(r_{k-1}r_{k}+r_{k-1}r_{k})\sqrt{T}+\Big(r_{k-1}r_{k}+r_{k-1}r_{k}
+\frac{1}{M_{k-1}}(r_{k-1}+r_{k-1}^{4})\Big)\sqrt{T}\\&+\frac{1}{M_{k-1}}\sqrt{r_{k-2}^{8}T+\|h\|_{L^{2}(\mathbb{S}^{2}; W_{2}^{6}(\Omega))}^{2}}
\\\leq& C(h)\Big(1+\frac{1}{M_{k}^{3}}+\frac{1}{M_{k-1}}\Big)(1+r_{k-2}^{4}+r_{k-1}^{4}+r_{k}^{4})\sqrt{1+T}
\end{aligned}
\end{equation}
according to the Sobolev embedding theorem, Young's inequality and (\ref{gronwall's inequality 1-1}).
Here we note that
\begin{equation}\nonumber
\begin{aligned}
\Big\|\frac{1}{\epsilon}\langle\vec{w}(f_{k-1}-\overline f_{k-1})\rangle\Big\|_{L^{2}(0, T; L^{\infty}(\Omega))}&\leq
C\Big\|\frac{1}{\epsilon}\langle\vec{w}(f_{k-1}-\overline f_{k-1})\rangle\Big\|_{L^{2}(0, T; W_{2}^{6}(\Omega))}\\&\leq C\Big\|\frac{1}{\epsilon}(f_{k-1}-\overline f_{k-1})\Big\|_{L^{2}(0, T; L^{2}(\mathbb{S}^{2}; W_{2}^{6}(\Omega)))}.
\end{aligned}
\end{equation}
We can also get
\begin{equation}\label{zaibu-1}
\begin{aligned}
&\|A_{0}(V_{k})\|_{L^{\infty}((0, T)\times\Omega)}\\\leq& \Big\|\frac{1}{\rho_{k}}\Big\|_{L^{\infty}((0, T)\times\Omega)}+\Big\| \frac{\rho_{k}}{\theta_{k}}\Big\|_{L^{\infty}((0, T)\times\Omega)}+\Big\|\frac{\rho_{k}}{\theta_{k}^{2}}\Big\|_{L^{\infty}((0, T)\times\Omega)}
\\\leq& C(1+\frac{1}{M_{k}^{2}})(1+r_{k}),
\end{aligned}
\end{equation}
\begin{equation}\label{zaibu-2}
\begin{aligned}
&\|D_{x}^{\gamma}A_{0}(V_{k})\|_{L^{\infty}((0, T); L^{2}(\Omega))}\\\leq&\Big\|D_{x}^{\gamma}\Big(\frac{1}{\rho_{k}}\Big)\Big\|_{L^{\infty}((0, T); L^{2}(\Omega))}+\Big\|D_{x}^{\gamma}\Big(\frac{\rho_{k}}{\theta_{k}}\Big)\Big\|_{L^{\infty}((0, T); L^{2}(\Omega))}+\Big\|D_{x}^{\gamma}\Big(\frac{\rho_{k}}{\theta_{k}^{2}}\Big)\Big\|_{L^{\infty}((0, T); L^{2}(\Omega))}
\\\leq& C\Big(1+\frac{1}{M_{k}^{8}}\Big)(1+r_{k}^{7})
\end{aligned}
\end{equation}
and
\begin{equation}\label{zaibu-3}
\begin{aligned}
&\|D_{x}^{\gamma}(A_{0}^{-1}(V_{k}))\|_{L^{\infty}((0, T); L^{2}(\Omega))}\\\leq& C(\|D_{x}^{\gamma}\rho_{k}\|_{L^{\infty}(0, T; L^{2}(\Omega))}+\Big\|D_{x}^{\gamma}\Big(\frac{\theta_{k}}{\rho_{k}}\Big)\Big\|_{L^{\infty}(0, T; L^{2}(\Omega))}+\Big\|D_{x}^{\gamma}\Big(\frac{\theta_{k}^{2}}{\rho_{k}}\Big)\Big\|_{L^{\infty}(0, T; L^{2}(\Omega))})\\\leq& C\Big(1+\frac{1}{M_{k}^{7}}\Big)(1+r_{k}^{8}).
\end{aligned}
\end{equation}
According to (\ref{xishujuzhen-1}), we have
\begin{equation}\nonumber
\lambda_{\min}=\min\Big\{\frac{1}{\rho_{k}}, \frac{\rho_{k}}{\theta_{k}}, \frac{\rho_{k}}{\theta_{k}^{2}}\Big\}\geq \frac{\min\{1, M_{k}\}}{C(1+r_{k}^{2})},
\end{equation}
where $\lambda_{\min}$ is the minimal eigenvalue of $A_{0}(V_{k})$. So,
\begin{equation}\label{zaibu-4}
\frac{1}{\lambda_{\min}}\leq  C\Big(1+\frac{1}{M_{k}}\Big)(1+r_{k}^{2}).
\end{equation}

For $A_{j}(V_{k})$, we have
\begin{equation}\label{zaibu-5}
\begin{aligned}
\Big\|\frac{\partial A_{j}(V_{k})}{\partial x_{j}}\Big\|_{L^{\infty}((0, T)\times\Omega)}\leq& C\Big(\sum_{j=1}^{3}\Big\|\frac{\partial}{\partial x_{j}}\Big( \frac{u_{k}^{j}}{\rho_{k}}\Big)\Big\|_{L^{\infty}((0, T)\times\Omega)}+\sum_{j=1}^{3}\Big\|\frac{\partial}{\partial x_{j}}\Big(\frac{\rho_{k}u_{k}^{j}}{\theta_{k}}\Big)\Big\|_{L^{\infty}((0, T)\times\Omega)}\\&+\sum_{j=1}^{3}\Big\|\frac{\partial}{\partial x_{j}}\Big(\frac{\rho_{k}u_{k}^{j}}{\theta_{k}^{2}}\Big)\Big\|_{L^{\infty}((0, T)\times\Omega)}+\sum_{j=1}^{3}\Big\|\frac{\partial}{\partial x_{j}}\Big(\frac{\rho_{k}}{\theta_{k}}\Big)\Big\|_{L^{\infty}((0, T)\times\Omega)}\Big)\\\leq& C\Big(1+\frac{1}{M_{k}^{3}}\Big)(1+r_{k}^{3})
\end{aligned}
\end{equation}
and
\begin{equation}\label{zaibu-6}
\begin{aligned}
\|D_{x}^{\gamma}A_{j}(V_{k})\|_{L^{\infty}(0, T; L^{2}(\Omega))}\leq& C\Big(\sum_{j=1}^{3}\Big\|D_{x}^{\gamma}\Big( \frac{u_{k}^{j}}{\rho_{k}}\Big)\Big\|_{L^{\infty}(0, T; L^{2}(\Omega))}+\sum_{j=1}^{3}\Big\|D_{x}^{\gamma}\Big(\frac{\rho_{k}u_{k}^{j}}{\theta_{k}}\Big)\Big\|_{L^{\infty}(0, T; L^{2}(\Omega))}\\&+\sum_{j=1}^{3}\Big\|D_{x}^{\gamma}\Big(\frac{\rho_{k}u_{k}^{j}}{\theta_{k}^{2}}\Big)\Big\|_{L^{\infty}(0, T; L^{2}(\Omega))}+\sum_{j=1}^{3}\Big\|D_{x}^{\gamma}\Big(\frac{\rho_{k}}{\theta_{k}}\Big)\Big\|_{L^{\infty}(0, T; L^{2}(\Omega))}\Big)\\\leq& C\Big(1+\frac{1}{M_{k}^{8}}\Big)(1+r_{k}^{8}).
\end{aligned}
\end{equation}
Combining with (\ref{zaibu-2}), (\ref{zaibu-4}) and (\ref{zaibu-6}),  we have
\begin{equation}\nonumber
\begin{aligned}
&\int_{0}^{T}\int_{\Omega}-\sum_{j=1}^{3}(A_{0}(V_{k})D_{x}^{\gamma}A_{0}^{-1})(V_{k})\Big(A_{j}(V_{k})\frac{\partial V^{0}}{\partial x_{j}}\Big)\cdot D_{x}^{\gamma}U_{k+1}\mathrm{d}\vec x \mathrm{d}t\\\leq& C\Big(\Big\|-\sum_{j=1}^{3}(A_{0}(V_{k})D_{x}^{\gamma}A_{0}^{-1})(V_{k})\Big(A_{j}(V_{k})\frac{\partial V^{0}}{\partial x_{j}}\Big)\Big\|_{L^{\infty}(0, T; L^{2}(\Omega))}^{2}+\|U_{k+1}\|_{L^{\infty}(0, T; W_{2}^{6}(\Omega))}^{2}\Big)T
\\\leq& C\Big((1+(\lambda_{min}^{-1}\|A_{0}\|_{L^{\infty}(0, T; W_{2}^{6}(\Omega))})^{6})\Big\|(A_{j}(V_{k})\frac{\partial V^{0}}{\partial x_{j}}\Big\|_{L^{\infty}(0, T; W_{2}^{5}(\Omega))}^{2}+\|U_{k+1}\|_{L^{\infty}(0, T; W_{2}^{6}(\Omega))}^{2}\Big)T
\\\leq& C\Big(1+\frac{1}{M_{k}^{62}}\Big)(1+r_{k}^{62})(\|V^{0}\|_{W_{2}^{6}(\Omega)}^{2}+\|U_{k+1}\|_{L^{\infty}(0, T; W_{2}^{6}(\Omega))}^{2})T.
\end{aligned}
\end{equation}

For $D_{x}^{\gamma}F(V_{k}, f_{k})$, we have
\begin{equation}\nonumber
\begin{aligned}
\int_{0}^{T}\int_{\Omega}|D_{x}^{\gamma}F(V_{k}, f_{k})|^{2}\mathrm{d}\vec x \mathrm{d}t\leq& C\int_{0}^{T}\Big\|\frac{1}{\theta_{k}}\frac{1}{\epsilon}\langle w_{j}(f_{k}-\overline f_{k})\rangle\Big\|_{W_{2}^{6}(\Omega)}^{2}+\Big\|\frac{1}{\theta_{k}^{2}}(\overline f_{k}-B(\theta_{k}))\Big\|_{W_{2}^{6}(\Omega)}^{2}\mathrm{d}t\\&+C\int_{0}^{T}\Big\|\frac{1}{\theta_{k}^{2}}\frac{1}{\epsilon}\langle w_{j}(f_{k}-\overline f_{k})\rangle\vec{u}_{k}\Big\|_{W_{2}^{6}(\Omega)}^{2}\mathrm{d}t
\\\leq& C\Big(\Big\|\frac{1}{\theta_{k}}\Big\|_{L^{\infty}(0, T;W_{2}^{6}(\Omega))}^{2}\Big\|\frac{1}{\epsilon}(f_{k}-\overline f_{k})\Big\|_{L^{2}(0, T; L^{2}(\mathbb{S}^{2}; W_{2}^{6}(\Omega)))}^{2}
\\&+\Big\|\frac{1}{\theta_{k}^{2}}\Big\|_{L^{\infty}(0, T;W_{2}^{6}(\Omega))}^{2}\|\overline f_{k}-B(\theta_{k})\|_{L^{\infty}(0, T;W_{2}^{6}(\Omega))}^{2}T\\&+\Big\|\frac{1}{\theta_{k}^{2}}\Big\|_{L^{\infty}(0, T;W_{2}^{6}(\Omega))}^{2}\Big\|\frac{1}{\epsilon}(f_{k}-\overline f_{k})\Big\|_{L^{2}(0, T; L^{2}(\mathbb{S}^{2}; W_{2}^{6}(\Omega)))}^{2}\\&\times\|\vec{u}_{k}\|_{L^{\infty}(0, T;W_{2}^{6}(\Omega))}^{2}\Big)
\\\leq & C\Big(\Big(1+\frac{1}{M_{k}^{16}}\Big)(1+r_{k}^{20})(1+r_{k-1}^{8})(1+T)(1+\|h\|_{L^{2}(\mathbb{S}^{2}; W_{2}^{6}(\Omega))}^{2})\Big)
\end{aligned}
\end{equation}
combining with inequality (\ref{gronwall's inequality 1-1}).

For $F_{\gamma}\cdot D_{x}^{\gamma}U_{k+1}$, we have
\begin{equation}\nonumber
\begin{aligned}
\int_{0}^{T}\int_{\Omega}F_{\gamma}\cdot D_{x}^{\gamma}U_{k+1}\mathrm{d}\vec x \mathrm{d}t\leq& C\sum_{j=1}^{3}\|A_{0}(V_{k})\|_{L^{\infty}((0, T)\times\Omega)}\|A_{0}^{-1}(V_{k})\|_{L^{\infty}(0, T; W_{2}^{6}(\Omega))}\\&\times\|A_{j}(V_{k})\|_{L^{\infty}(0, T; W_{2}^{6}(\Omega))}\|U_{k+1}\|_{L^{\infty}(0, T; W_{2}^{6}(\Omega))}^{2}T.
\end{aligned}
\end{equation}
Then, we get
\begin{equation}\nonumber
\begin{aligned}
\int_{0}^{T}\int_{\Omega}F_{\gamma}\cdot D_{x}^{\gamma}U_{k+1}\mathrm{d}\vec x \mathrm{d}t\leq& C\Big(1+\frac{1}{M_{k}^{2}}\Big)(1+r_{k})\Big(1+\frac{1}{M_{k}^{7}}\Big)(1+r_{k}^{8})\Big(1+\frac{1}{M_{k}^{8}}\Big)(1+r_{k}^{8})\\&\times\|U_{k+1}\|_{L^{\infty}(0, T; W_{2}^{6}(\Omega))}^{2}T\\\leq& C\Big(1+\frac{1}{M_{k}^{17}}\Big)(1+r_{k}^{17})\|U_{k+1}\|_{L^{\infty}(0, T; W_{2}^{6}(\Omega))}^{2}T
\end{aligned}
\end{equation}
combining with (\ref{zaibu-1}), (\ref{zaibu-3}) and (\ref{zaibu-6}).

Integrating (\ref{energy equality-2}) over $(0, t)$ and taking the sup over $t\in (0, T)$, we have
\begin{equation}
\begin{aligned}
\|U_{k+1}\|_{L^{\infty}(0, T; W_{2}^{6}(\Omega))}^{2}\leq& C\lambda_{min}^{-1}\Big(\|\partial_{t}A_{0}(V_{k})\|_{L^{2}(0, T; L^{\infty}(\Omega)))}\|U_{k+1}\|_{L^{\infty}(0, T; W_{2}^{6}(\Omega))}^{2}\sqrt{T}\\&+\sum_{j=1}^{3}\|D_{x}A_{j}(V_{k})\|_{L^{\infty}((0, T)\times\Omega))}\|U_{k+1}\|_{L^{\infty}(0, T; W_{2}^{6}(\Omega))}^{2}T
\\&+\|A_{0}(V_{k})\|_{L^{\infty}((0, T)\times\Omega)}\|A_{0}^{-1}(V_{k})\|_{L^{\infty}(0, T; W_{2}^{6}(\Omega))}\|F(V_{k}, f_{k})\|_{L^{2}(0, T; W_{2}^{6}(\Omega))}\\&\times\|U_{k+1}\|_{L^{\infty}(0, T; W_{2}^{6}(\Omega))}\sqrt{T}+C\sum_{j=1}^{3}\|A_{0}(V_{k})\|_{L^{\infty}((0, T)\times\Omega)}\\&\times\|A_{0}^{-1}(V_{k})\|_{L^{\infty}(0, T; W_{2}^{6}(\Omega))}\|A_{j}(V_{k})\|_{L^{\infty}(0, T; W_{2}^{6}(\Omega))}\|U_{k+1}\|_{L^{\infty}(0, T; W_{2}^{6}(\Omega))}^{2}T\\&+C\Big(\Big\|-\sum_{j=1}^{3}(A_{0}(V_{k})D_{x}^{\gamma}A_{0}^{-1})(V_{k})\Big(A_{j}(V_{k})\frac{\partial V^{0}}{\partial x_{j}}\Big)\Big\|_{L^{\infty}(0, T; L^{2}(\Omega))}^{2}\\&+\|U_{k+1}\|_{L^{\infty}(0, T; W_{2}^{6}(\Omega))}^{2}\Big)T\Big).\label{ineq 1}
\end{aligned}
\end{equation}
For $0<T<1$, we obtain that
\begin{equation}\nonumber
\begin{aligned}
\|U_{k+1}\|_{L^{\infty}(0, T; W_{2}^{6}(\Omega))}^{2}\leq& C\Big(1+\frac{1}{M_{k}}\Big)(1+r_{k}^{2})
\Big\{\Big(1+\frac{1}{M_{k}^{3}}+\frac{1}{M_{k-1}}\Big)(1+r_{k-2}^{4}+r_{k-1}^{4}+r_{k}^{4})\sqrt{1+T}\\&\times\|U_{k+1}\|_{L^{\infty}(0, T; W_{2}^{6}(\Omega))}^{2}\sqrt{T}+\Big(1+\frac{1}{M_{k}^{3}}\Big)(1+r_{k}^{3})\|U_{k+1}\|_{L^{\infty}(0, T; W_{2}^{6}(\Omega))}^{2}T\\&+\Big(\sqrt{\Big(1+\frac{1}{M_{k}^{16}}\Big)(1+r_{k}^{20})(1+r_{k-1}^{8})(1+T)}(1+C\|h\|_{L^{2}(\mathbb{S}^{2}; W_{2}^{6}(\Omega))})\Big)\\&\times\|U_{k+1}\|_{L^{\infty}(0, T; W_{2}^{6}(\Omega))}\Big(\Big(1+\frac{1}{M_{k}^{2}}\Big)(1+r_{k})\Big(1+\frac{1}{M_{k}^{7}}\Big)(1+r_{k}^{8})\Big)\sqrt{T}
\\&+\Big(1+\frac{1}{M_{k}^{17}}\Big)(1+r_{k}^{17})\|U_{k+1}\|_{L^{\infty}(0, T; W_{2}^{6}(\Omega))}^{2}T+\Big(1+\frac{1}{M_{k}^{62}}\Big)(1+r_{k}^{62})\\&\times(\|V^{0}\|_{W_{2}^{6}(\Omega)}^{2}+\|U_{k+1}\|_{L^{\infty}(0, T; W_{2}^{6}(\Omega))}^{2})T\Big\}
\\\leq& C(a, h, V^{0})\Big(1+\frac{1}{M_{k-1}}\Big)\Big(1+\frac{1}{M_{k}^{63}}\Big)(1+r_{k-2}^{4})(1+r_{k-1}^{4})(1+r_{k}^{64})\\&\times(\|U_{k+1}\|_{L^{\infty}(0, T; W_{2}^{6}(\Omega))}^{2}+1)\sqrt{T}.
\end{aligned}
\end{equation}
Using (\ref{gronwall's inequality 1-1}), we have
\begin{equation}\label{jiaeq-1}
\begin{aligned}
&\|U_{k+1}\|_{L^{\infty}(0, T; W_{2}^{6}(\Omega))}^{2}+\|f_{k+1}\|_{L^{\infty}(0, T; L^{2}(\mathbb{S}^{2}; W_{2}^{6}(\Omega)))}^{2}\\\leq& C(a, h, V^{0}, M_{k-1}, M_{k}, r_{k-2}, r_{k-1}, r_{k})(\|U_{k+1}\|_{L^{\infty}(0, T; W_{2}^{6}(\Omega))}^{2}+\|f_{k+1}\|_{L^{\infty}(0, T; L^{2}(\mathbb{S}^{2}; W_{2}^{6}(\Omega)))}^{2}+1)\sqrt{T}.
\end{aligned}
\end{equation}
For small $T=T(k)$, we get that
\begin{equation}\nonumber
\begin{aligned}
&\|U_{k+1}\|_{L^{\infty}(0, T; W_{2}^{6}(\Omega))}^{2}+\|f_{k+1}\|_{L^{\infty}(0, T; L^{2}(\mathbb{S}^{2}; W_{2}^{6}(\Omega)))}^{2}\\\leq& \frac{C(a, h, V^{0}, M_{k-1}, M_{k}, r_{k-2}, r_{k-1}, r_{k})\sqrt{T}}{1-C(a, h, V^{0}, M_{k-1}, M_{k}, r_{k-2}, r_{k-1}, r_{k})\sqrt{T}}\\:=&\frac{1}{2}(C'(a, h, V^{0}, M_{k-1}, M_{k}, r_{k-2}, r_{k-1}, r_{k}, T)T^{1/4})^{2}
\end{aligned}
\end{equation}
where $C^{'}(a, h, V^{0}, M_{k-1}, M_{k}, r_{k-2}, r_{k-1}, r_{k}, T)$ decreases as $T$ decreases or $M_{k-1}$ and $M_{k}$ increase. So it holds that
\begin{equation}\nonumber
r_{k+1}\leq C^{'}(a, h, V^{0}, M_{k-1}, M_{k}, r_{k-2}, r_{k-1}, r_{k}, T)T^{1/4}+\|V^{0}\|_{W_{2}^{6}(\Omega)}.
\end{equation}

Since $\partial_{t}\rho_{k+1}=-\vec{u}_{k}\cdot\nabla_{x}\rho_{k+1}-\rho_{k}\mathrm{div}_{x}\vec{u}_{k+1}$ and $\partial_{t}\theta_{k+1}=-\theta_{k}\mathrm{div}_{x}\vec{u}_{k+1}-\vec{u}_{k}\cdot\nabla_{x}\theta_{k+1}
+\frac{1}{\rho_{k}}(\overline f_{k}-B(\theta_{k}))-\frac{1}{\rho_{k}}(\frac{1}{\epsilon}+\epsilon)\langle\vec{w}(f_{k}-\overline f_{k})\rangle\cdot\vec{u}_{k}$ for $k\geq 1$, by egrating the above two equalities over $(0,t)$ and taking the inf on $(0, T)\times\Omega$,  we obtain that
\begin{equation}\label{lower bound-1}
\begin{aligned}
\inf_{t, \vec x}\rho_{k+1}(t, \vec x)&\geq \rho^{0}(\vec x)-CT(1+r_{k}r_{k+1}+r_{k}r_{k+1})\\&\geq a-CT(1+r_{k}r_{k+1})
\\&\geq a-C''(a, h, V^{0}, M_{k-1}, M_{k}, r_{k-1}, r_{k}, T)T
\end{aligned}
\end{equation}
and
\begin{equation}\label{lower bound-2}
\begin{aligned}
\inf_{t, \vec x}\theta_{k+1}(t, \vec x)\geq& \theta^{0}(\vec x)-CT(1+\frac{1}{M_{k}})(1+r_{k}r_{k+1}+r_{k}r_{k+1}+r_{k}+r_{k}^{4})\\&-C(1+\frac{1}{M_{k}})(1+r_{k-1}^{4})(1+r_{k})\sqrt{(1+T)}\sqrt{T}\\\geq& a-CT(1+\frac{1}{M_{k}})(1+r_{k}r_{k+1}+r_{k}^{3}r_{k+1}+r_{k}+r_{k}^{4})\\&-C(1+\frac{1}{M_{k}})(1+r_{k-1}^{4})(1+r_{k})\sqrt{(1+T)}\sqrt{T}\\\geq& a-C'''(a, h, V^{0}, M_{k-1}, M_{k}, r_{k-2}, r_{k-1}, r_{k}, T)\sqrt{T}
\end{aligned}
\end{equation}
where $C''$ and $C'''$ decrease as $T$ decreases or $M_{k-1}$ and $M_{k}$ increases.

Note that $\partial_{t}A_{0}(V_{k})=0$ for $k=0$ implies the terms $1/M_{-1}$ and $r_{-1}$ won't appear in (\ref{lower bound-2}).
Then, we have
\begin{equation}\label{jiaeq-1}
\begin{aligned}
\|V_{1}\|_{L^{\infty}(0, T; W_{2}^{6}(\Omega))}+\|f_{1}\|_{L^{\infty}(0, T; L^{2}(\mathbb{S}^{2}; W_{2}^{6}(\Omega)))}\leq& C'(a, b, h, V^{0}, T)T^{1/4}+\|V^{0}\|_{W_{2}^{6}(\Omega)}.
\end{aligned}
\end{equation}
Note that $M_{0}=a$ and $r_{0}=\|V^{0}\|_{W_{2}^{6}(\Omega)}$. Then, we can choose $T=T_{1}$ where $T_{1}$ sufficiently small, such that

\begin{equation}\nonumber
r_{1}< C'(a, h, V^{0}, T_{1})T_{1}^{1/4}+\|V^{0}\|_{W_{2}^{6}(\Omega)}+1:=r
\end{equation}
According to (\ref{lower bound-1}) and (\ref{lower bound-2}), we can choose $T=T_{2}$ where $T_{2}$ is small enough such that
\begin{equation}\nonumber
\inf_{t, \vec x}\rho_{1}(t, \vec x), \inf_{t, \vec x}\theta_{1}(t, \vec x)\geq a-\sup_{0\leq r_{0}\leq r; M_{0}\geq a/2}\hat C(a, h, V^{0}, r_{0}, T)\sqrt T\geq a/2
\end{equation}
where $\hat{C}(\cdot)=\max\{C''(\cdot), C'''(\cdot)\}$.

We can also choose $T=T_{3}$ where $T_{3}$ is sufficiently small such that
\begin{equation}\label{jiaeq-2}
\begin{aligned}
\|V_{2}\|_{L^{\infty}(0, T; W_{2}^{6}(\Omega))}+\|f_{2}\|_{L^{\infty}(0, T; L^{2}(\mathbb{S}^{2}; W_{2}^{6}(\Omega)))}\leq& \sup_{\begin{matrix}
  0\leq r_{0}, r_{1}\leq r, \\
  M_{0}, M_{1}\geq a/2
\end{matrix}}C'(a, b, h, V^{0}, r_{0}, r_{1}, T)T^{1/4}\\&+\|V^{0}\|_{W_{2}^{6}(\Omega)}\\\leq& r.
\end{aligned}
\end{equation}
So we conclude that $r_{2}\leq r$.

Now we choose $T=T_{4}$ where $T_{4}$ is small enough such that
\begin{equation}\nonumber
\inf_{t, \vec x}\rho_{2}(t, \vec x), \inf_{t, \vec x}\theta_{2}(t, \vec x)\geq a-\sup_{\begin{matrix}
  0\leq r_{0}, r_{1}\leq r, \\
  M_{0}, M_{1}\geq a/2
\end{matrix}}\hat C(a, b, h, V^{0}, r_{1}, T)\sqrt T\geq a/2.
\end{equation}
So it implies that $M_{2}\geq M=a/2.$

Similarly, we can choose $T=T_{5}$ where $T_{5}$ is small enough such that
\begin{equation}\label{jiaeq-2-1}
\begin{aligned}
\|V_{3}\|_{L^{\infty}(0, T; W_{2}^{6}(\Omega))}+\|f_{3}\|_{L^{\infty}(0, T; L^{2}(\mathbb{S}^{2}; W_{2}^{6}(\Omega)))}\leq& \sup_{\begin{matrix}
  0\leq r_{0}, r_{1}, r_{2}\leq r, \\
   M_{1}, M_{2}\geq a/2
\end{matrix}}C'(a, b, h, V^{0}, r_{0}, r_{1},r_{2}, T)T^{1/4}\\&+\|V^{0}\|_{W_{2}^{6}(\Omega)}\\\leq& r
\end{aligned}
\end{equation}
and
\begin{equation}
\inf_{t, \vec x}\rho_{3}(t, \vec x), \inf_{t, \vec x}\theta_{3}(t, \vec x)\geq a-\sup_{\begin{matrix}
  0\leq r_{0}, r_{1}, r_{2}\leq r, \\
  M_{1}, M_{2}\geq a/2
\end{matrix}}\hat C(a, b, h, V^{0}, r_{0}, r_{1}, r_{2}, T)\sqrt T\geq a/2.
\end{equation}
So we have that $r_{3}\leq r$ and $M_{3}\geq M=a/2.$

Choose $T=T'=\min\{1, T_{1}, T_{2}, T_{3}, T_{4}, T_{5}\}$ and assume that $0\leq r_{k}\leq r$ and $M_{k}\geq M=a/2$. Then, we have
\begin{equation}\label{jiaeq-2}
\begin{aligned}
r_{k+1}\leq \sup_{\begin{matrix}
  0\leq r_{k-2}, r_{k-1}, r_{k}\leq r, \\
  M_{k-1}, M_{k}\geq a/2
\end{matrix}}C'(a, b, h, V^{0}, M_{k-1}, M_{k}, r_{k-2}, r_{k-1}, r_{k}, T)T^{1/4}+\|V^{0}\|_{W_{2}^{6}(\Omega)}\leq r.
\end{aligned}
\end{equation}
and
\begin{equation}\nonumber
\begin{aligned}
\inf_{t, \vec x}\rho_{k+1}(t, \vec x), \inf_{t, \vec x}\theta_{2}(t, \vec x)\geq& a-\sup_{\begin{matrix}
  0\leq r_{k-2}, r_{k-1}, r_{k}\leq r, \\
  M_{k-1}, M_{k}\geq a/2
\end{matrix}}\hat C(a, b, h, V^{0}, M_{k-1}, M_{k}, r_{k-2}, r_{k-1}, r_{k}, T)\sqrt T\\\geq& a/2.
\end{aligned}
\end{equation}
Therefore, we have proved $(f_{k+1}, V_{k+1})\in G$ and $\rho_{k+1}, \theta_{k+1}\geq M=a/2$ for any $k\geq 0$.

Define $P_{k+1}=f_{k+1}-f_{k}, Q_{k+1}=V_{k+1}-V_{k}$. Then we get the system about $(P_{k+1}, Q_{k+1})$ as follows
\begin{equation}\label{equations about Pk and Qk}\left\{
\begin{split}
&A_{0}(V_{k})\partial_{t}Q_{k+1}+\sum_{j=1}^{3}A_{j}(V_{k})\frac{\partial Q_{k+1}}{\partial x_{j}}
=F(V_{k}, f_{k})-F(V_{k-1}, f_{k-1})+L\ \ \mathrm{in} \ \ (0, T]\times \Omega, \\
&\partial_{t}P_{k+1}+\frac{1}{\epsilon}\vec w\cdot\nabla_{x}P_{k+1}+\frac{1}{\epsilon^{2}}(P_{k+1}-\overline P_{k+1})+P_{k+1}=B(\theta_{k})-B(\theta_{k-1})\ \ \mathrm{in} \ \ (0, T]\times \Omega\times \mathbb{S}^{2},\\
&P_{k+1}(0, \vec{x}, \vec{w})=0\ \ \mathrm{in} \ \ \Omega\times \mathbb{S}^{2},\\
&Q_{k+1}(0, \vec{x})=0\ \ \mathrm{in}  \ \ \Omega,\\
\end{split}\right.
\end{equation}
where
\begin{equation}\nonumber
L=-(A_{0}(V_{k})-A_{0}(V_{k-1}))\frac{\partial V_{k}}{\partial t}-\sum_{j=1}^{3}(A_{j}(V_{k})-A_{j}(V_{k-1}))\frac{\partial V_{k}}{\partial x_{j}}.
\end{equation}

Multiplying $P_{k+1}$ on both sides of $(\ref{equations about Pk and Qk})_{1}$, we have
\begin{equation}\label{Eq 2}
\begin{split}
&\frac{\mathrm{d}}{\mathrm{d}t}\int_{\Omega}\int_{\mathbb{S}^{2}}|P_{k+1}|^{2}\mathrm{d}\vec{w}\mathrm{d}\vec{x}
+\frac{2}{\epsilon^{2}}\int_{\Omega}\int_{\mathbb{S}^{2}}|P_{k+1}-\overline P_{k+1}|^{2}\mathrm{d}\vec{w}\mathrm{d}\vec{x}+\int_{\Omega}\int_{\mathbb{S}^{2}}|P_{k+1}|^{2}\mathrm{d}\vec{w}\mathrm{d}\vec{x}
\\=&2\int_{\Omega}\int_{\mathbb{S}^{2}}(B(\theta_{k})-B(\theta_{k-1}))P_{k+1}\mathrm{d}\vec{w}\mathrm{d}\vec{x}.
\end{split}
\end{equation}
Integrating (\ref{Eq 2}) over $(0, t)$ and then taking the sup over $t\in(0, T)$, we have
\begin{equation}\label{gronwall's inequality 6}
\begin{split}
&\|P_{k+1}\|_{L^{\infty}(0, T; L^{2}(\mathbb{S}^{2}; L^{2}(\Omega))}^{2}+\|(P_{k+1}-\overline P_{k+1})/\epsilon\|_{L^{2}(0, T; L^{2}(\mathbb{S}^{2}; L^{2}(\Omega)))}^{2}
\\\leq&C(r, \Omega)T\|Q_{k}\|_{L^{\infty}(0, T; L^{2}(\Omega))}^{2}
\end{split}
\end{equation}
in terms of the Young's inequality. Then we get that
\begin{equation}\label{eq 3}
\begin{split}
\|P_{k+1}\|_{L^{\infty}(0, T; L^{2}(\mathbb{S}^{2}; L^{2}(\Omega))}+\|(P_{k+1}-\overline P_{k+1})/\epsilon\|_{L^{2}(0, T; L^{2}(\mathbb{S}^{2}; L^{2}(\Omega)))}
\leq&C(r, \Omega)\sqrt{T}\|Q_{k}\|_{L^{\infty}(0, T; L^{2}(\Omega))}.
\end{split}
\end{equation}
From Sobolev's embedding theorem and Taylor's expansion theorem, we deduce that
\begin{equation}\nonumber
\|L\|_{L^{2}(0, T; L^{2}(\Omega))}\leq C(r)\|Q_{k}\|_{L^{\infty}(0, T; L^{2}(\Omega))}\sqrt{T}.
\end{equation}
We find from the mean value theorem and the inequalities (\ref{gronwall's inequality 1-1}) and (\ref{eq 3}) that
\begin{equation}
\begin{aligned}
&\|F(V_{k}, f_{k})-F(V_{k-1}, f_{k-1})\|_{L^{2}(0, T; L^{2}(\Omega))}
\\\leq& C(a, r)\Big\{\|(f_{k}-\overline f_{k})/\epsilon\|_{L^{2}(0, T; L^{2}(\mathbb{S}^{2}; L^{2}(\Omega)))}\|Q_{k}\|_{L^{2}(0, T; L^{2}(\Omega))}\\&+\|(P_{k}-\overline P_{k})/\epsilon\|_{L^{2}(0, T; L^{2}(\mathbb{S}^{2}; L^{2}(\Omega)))}\|V_{k-1}\|_{L^{2}(0, T; L^{2}(\Omega))}\Big\}
\\\leq& C(a, r)\sqrt{T}(\|Q_{k}\|_{L^{\infty}(0, T; L^{2}(\Omega))}+\|Q_{k-1}\|_{L^{\infty}(0, T; L^{2}(\Omega))}).
\end{aligned}
\end{equation}
Define
\begin{equation}\nonumber
E(t):=\frac{1}{2}\int_{\Omega}(A_{0}(V_{k})Q_{k+1})\cdot Q_{k+1}\mathrm{d}\vec x.
\end{equation}
Similarly, we can obtain
\begin{equation}\nonumber
\begin{aligned}
\|Q_{k+1}\|_{L^{\infty}(0, T; L^{2}(\Omega))}\leq& C(a, r)(\|Q_{k+1}\|_{L^{\infty}(0, T; L^{2}(\Omega))}\sqrt{T}+\|Q_{k}\|_{L^{\infty}(0, T; L^{2}(\Omega))}\sqrt{T}\\&+\|Q_{k-1}\|_{L^{\infty}(0, T; L^{2}(\Omega))}\sqrt{T}).
\end{aligned}
\end{equation}
after a careful calculation.
Then,
\begin{equation}\nonumber
\begin{aligned}
&\|P_{k+1}\|_{L^{\infty}(0, T; L^{2}(\mathbb{S}^{2}; L^{2}(\Omega))}+\|Q_{k+1}\|_{L^{\infty}(0, T; L^{2}(\Omega))}\\\leq& C(a, r)(\|Q_{k+1}\|_{L^{\infty}(0, T; L^{2}(\Omega))}\sqrt{T}+\|Q_{k}\|_{L^{\infty}(0, T; L^{2}(\Omega))}\sqrt{T}+\|Q_{k-1}\|_{L^{\infty}(0, T; L^{2}(\Omega))}\sqrt{T}
\\&+\|P_{k}\|_{L^{\infty}(0, T; L^{2}(\mathbb{S}^{2}; L^{2}(\Omega))}\sqrt{T}).
\end{aligned}
\end{equation}
Choose small $T=T_{6}$ such that $C(a, r)\sqrt{T}<1$. Then we have
\begin{equation}\nonumber
\begin{aligned}
&\|P_{k+1}\|_{L^{\infty}(0, T; L^{2}(\mathbb{S}^{2}; L^{2}(\Omega))}+\|Q_{k+1}\|_{L^{\infty}(0, T; L^{2}(\Omega))}\\\leq& C(a, r)(\|Q_{k}\|_{L^{\infty}(0, T; L^{2}(\Omega))}\sqrt{T}+\|Q_{k-1}\|_{L^{\infty}(0, T; L^{2}(\Omega))}\sqrt{T}+\|P_{k}\|_{L^{\infty}(0, T; L^{2}(\mathbb{S}^{2}; L^{2}(\Omega))}\sqrt{T}).
\end{aligned}
\end{equation}

Now, we choose small $T=T_{7}$ such that $d=C(a, r)\sqrt{T}<1/2$. Then we obtain that
\begin{equation}\nonumber
\begin{split}
&\|P_{k+1}\|_{L^{\infty}(0, T; L^{2}(\mathbb{S}^{2}; W_{2}^{6}(\Omega))}+\|Q_{k+1}\|_{L^{\infty}(0, T; W_{2}^{6}(\Omega))}
\\\leq&d(\|P_{k}\|_{L^{\infty}(0, T; L^{2}(\mathbb{S}^{2}; W_{2}^{6}(\Omega))}+\|Q_{k}\|_{L^{\infty}(0, T; W_{2}^{6}(\Omega))}+\|Q_{k-1}\|_{L^{\infty}(0, T; W_{2}^{6}(\Omega))}).
\end{split}
\end{equation}
So it holds that
\begin{equation}\nonumber
\begin{split}
&\sum_{k=0}^{\infty}\|P_{k}\|_{L^{\infty}(0, T; L^{2}(\mathbb{S}^{2}; W_{2}^{6}(\Omega))}+\|Q_{k}\|_{L^{\infty}(0, T; W_{2}^{6}(\Omega))}
<\infty.
\end{split}
\end{equation}
Therefore, we have $(f_{k}, V_{k})$ converges strongly to $(f^{\epsilon}, V^{\epsilon})$ in $L^{\infty}(0, T; L^{2}(\mathbb{S}^{2}; L^{2}(\Omega)))\times L^{\infty}(0, T; L^{2}(\Omega))$. So, we can choose
\begin{equation}\nonumber
T=\min\{1, T_{1}, T_{2}, T_{3}, T_{4}, T_{5}, T_{6}, T_{7}\}
\end{equation}
which is independent of $k$ and $\epsilon$ to satisfy all the above uniform estimates and contraction principle. Therefore, we can infer that $(f^{\epsilon}, V^{\epsilon})$ is the solution of (\ref{research equations}) where $(f^{\epsilon}, V^{\epsilon})\in G$.

By the fact that $(f_{k}, V_{k})\in G$, we can get the strong convergence of $(f_{k}, V_{k})$ in $C^{0}([0, T]; L^{2}(\mathbb{S}^{2}; W_{2}^{6-\delta}(\Omega)))\times C^{0}([0, T]; W^{6-\delta}(\Omega))$ as $k\rightarrow\infty$ for any $\delta>0$ by Sobolev's theorem.

Finally, the uniqueness of (\ref{research equations}) can be obtained by the energy method.  We have finished the proof of the existence and uniqueness of the solution to the  system (\ref{research equations}).

We write $f^{\epsilon}=\langle f^{\epsilon}\rangle+\epsilon\hat g^{\epsilon},$ where $\hat g^{\epsilon}=\frac{1}{\epsilon}(f^{\epsilon}-\langle f^{\epsilon}\rangle).$ We define the macroscopic quantities $\vec{J}^{\epsilon}=\frac{1}{\epsilon}\langle\vec{w}f^{\epsilon}\rangle.$ Thus,
\begin{equation}\nonumber
\vec{J}^{\epsilon}=\langle \vec{w}\hat g^{\epsilon}\rangle.
\end{equation}

\begin{lem}\label{main result-6}
The following properties are satisfied:

\noindent$\mathrm{(i)}$ $\{f^{\epsilon}\}_{\epsilon>0}$ is bounded in $L^{\infty}(0, T; L^{2}(\mathbb{S}^{2}; W_{2}^{6}(\Omega)))$,

\noindent$\mathrm{(ii)}$ $\{\overline f^{\epsilon}\}_{\epsilon>0}$ is bounded in $L^{\infty}(0, T; W_{2}^{6}(\Omega))$,

\noindent$\mathrm{(iii)}$ $\{\hat g^{\epsilon}\}_{\epsilon>0}$ is bounded in $L^{2}(0, T; L^{2}(\mathbb{S}^{2}; W_{2}^{6}(\Omega)))$,

\noindent$\mathrm{(iv)}$ $\{\vec{J}^{\epsilon}\}_{\epsilon>0}$ is bounded in $L^{2}(0, T; W_{2}^{6}(\Omega))$,

\noindent$\mathrm{(v)}$ $\{V^{\epsilon}\}_{\epsilon>0}$ is bounded $L^{\infty}(0, T; W_{2}^{6}(\Omega))$.
\end{lem}
\noindent\textbf{Proof.}  According to the proof process in the first part of Theorem \ref{definition weak solution-0}, we can get $\mathrm{(i)}$, $\mathrm{(ii)}$, $\mathrm{(iii)}$, $\mathrm{(iv)}$ and $\mathrm{(v)}$.

\section{Moment equations}\label{si}
In this section, we discuss the equations satisfied by the moments of $f^{\epsilon}$.
We integrate $(\ref{research equations})_{1}$ over $\mathbb{S}^{2}$ to obtain
\begin{equation}\label{moment equation}
\partial_{t}\overline f^{\epsilon}+\mathrm{div}_{x}\vec{J}^{\epsilon}=B(\theta^{\epsilon})-\overline f^{\epsilon},
\end{equation}
where $div_{x}\vec{J}^{\epsilon}=\frac{1}{\epsilon}\frac{1}{4\pi}\int_{\mathbb{S}^{2}}\vec{w}\cdot\nabla_{x}f^{\epsilon}\mathrm{d}\vec w=\frac{1}{4\pi}\int_{\mathbb{S}^{2}}\vec{w}\cdot\nabla_{x}\hat {g}^{\epsilon}\mathrm{d}\vec w\in L^{2}(0, T; W_{2}^{5}(\Omega))$.

Multiplying $(\ref{research equations})_{1}$ by $\vec{w}$, integrating on $\mathbb{S}^{2}$, we can derive
\begin{equation}\label{moment equation-2}
\nabla_{x}\overline f^{\epsilon}=(\langle\vec{w}\otimes\vec{w}\rangle)^{-1}(-\epsilon^{2}\partial_{t}\vec{J}^{\epsilon}
-\epsilon\mathrm{div}_{x}(\langle\vec{w}\otimes\vec{w}\hat{g}^{\epsilon}\rangle)+\frac{1}{\epsilon}\langle\vec{w}(\overline f^{\epsilon}-f^{\epsilon})\rangle+\epsilon \langle\vec w(B(\theta^{\epsilon})-f^{\epsilon})\rangle
\end{equation}
after a careful calculation.
\begin{lem}\label{main result-3}
The sequence $\{\overline f^{\epsilon}\}_{\epsilon>0}$ lies in a (strong) compact set of $C^{0}([0, T]; W_{2}^{5}(\Omega))$. It is also compact in $C^{0}_{weak}([0, T]; W_{2}^{6}(\Omega)).$

\noindent {\bf Proof.} Since $\overline f^{\epsilon}$ is bounded in $L^{\infty}(0, T; W_{2}^{6}(\Omega))$ and $\partial_{t}\overline f^{\epsilon}$ is bounded in $L^{2}(0, T; W_{2}^{5}(\Omega))$, we can get that $\{\overline f^{\epsilon}\}_{\epsilon>0}$ lies a compact set of $C^{0}([0, T]; W_{2}^{5}(\Omega))$ by the Aubin-Lions-Simon Lemma. According to the above conditions about $\overline f^{\epsilon}$ and $\partial_{t}\overline f^{\epsilon}$, we conclude that $\{\overline f^{\epsilon}\}_{\epsilon>0}$ lies a compact set of $C^{0}_{weak}([0, T]; W_{2}^{6}(\Omega)).$ \hfill$\Box$
\end{lem}

\begin{lem}\label{main result-3-1}
The sequence $\{V^{\epsilon}\}_{\epsilon>0}$ lie in a (strong) compact set of $C^{0}([0, T]; W_{2}^{5}(\Omega))$ for any $\delta>0.$ It is also compact in $C^{0}_{weak}([0, T]; W_{2}^{6}(\Omega)).$

\noindent {\bf Proof.} Since $V^{\epsilon}$ is bounded in $L^{\infty}(0, T; W_{2}^{6}(\Omega))$ and $\partial_{t}V^{\epsilon}$ is bounded in $L^{2}(0, T; W_{2}^{5}(\Omega))$, we can get that $\{V^{\epsilon}\}_{\epsilon>0}$ lies a compact set of $C^{0}([0, T]; W_{2}^{5}(\Omega))$ by the Aubin-Lions-Simon Lemma. According to the above conditions about $V^{\epsilon}$ and $\partial_{t}V^{\epsilon}$, we can obtain that $\{V^{\epsilon}\}_{\epsilon>0}$ lies a compact set of $C^{0}_{weak}([0, T]; W_{2}^{6}(\Omega)).$\hfill$\Box$
\end{lem}
\noindent\textbf{Passage to the limit.} We are now ready to pass to the limit in (1.1). By Lemmas \ref{main result}, \ref{main result-3} and \ref{main result-3-1}, up to a subsequence, we have

\noindent$\mathrm{(i)}$$f^{\epsilon}\rightharpoonup f$ in $L^{\infty}(0, T; L^{2}(\mathbb{S}^{2}; W_{2}^{6}(\Omega)))$ in the $weak-\ast$ sense,

\noindent$\mathrm{(ii)}$$\overline f^{\epsilon}\rightarrow \overline f$ strongly in $C^{0}([0, T]; W_{2}^{5}(\Omega))$ and $C^{0}_{weak}([0, T]; W_{2}^{6}(\Omega))$,

\noindent$\mathrm{(iii)}$$\vec{J}^{\epsilon}\rightharpoonup \vec{J}$ weakly in $L^{2}(0, T; W_{2}^{6}(\Omega))$,

\noindent$\mathrm{(iv)}$$V^{\epsilon}\rightarrow V$ weakly in $L^{\infty}(0, T; W_{2}^{6}(\Omega))$ in the $weak-\ast$ sense,

\noindent$\mathrm{(v)}$$V^{\epsilon}\rightarrow V$ strongly in $C^{0}([0, T]; W_{2}^{5}(\Omega))$ and $C^{0}_{weak}([0, T]; W_{2}^{6}(\Omega))$.

The limit $f$ is the macroscopic quantity $\overline f$ since we have
$f^{\epsilon}=\overline f^{\epsilon}+\epsilon \hat{g}^{\epsilon}\rightharpoonup f=\overline f+0$ in $L^{2}(0, T; L^{2}(\mathbb{S}^{2}; W_{2}^{6}(\Omega)))$ in the weak sense.
\begin{rem}
In fact, we have $f^{\epsilon}\rightarrow \overline f$ strongly in $L^{2}(0, T; L^{2}(\mathbb{S}^{2}; W_{2}^{5}(\Omega)))$. We can obtain the uniform estimate $\|f^{\epsilon}-\overline f^{\epsilon}\|_{L^{2}(0, T; L^{2}(\mathbb{S}^{2}; W_{2}^{5}(\Omega)))}=O(\epsilon)$ according to Lemma \ref{main result-6}. Moreover, $\|\overline f^{\epsilon}-\overline f\|_{L^{2}(0, T; L^{2}(\mathbb{S}^{2}; W_{2}^{5}(\Omega)))}=C\|\overline f^{\epsilon}-\overline f\|_{L^{2}(0, T; W_{2}^{5}(\Omega))}$. Therefore, we have $\|f^{\epsilon}-\overline f\|_{L^{2}(0, T; L^{2}(\mathbb{S}^{2};W_{2}^{5}(\Omega)))}\leq C\|\overline f^{\epsilon}-\overline f\|_{L^{2}(0, T; W_{2}^{5}(\Omega))}+O(\epsilon)$ which implies the desired strong convergence.
\end{rem}

\begin{lem}\label{main result-3-bu-bu}
The limit $\overline f$ of $\overline f^{\epsilon}$ satisfies $\overline f\in L^{2}(0, T; W_{2}^{7}(\Omega))$.

\noindent {\bf Proof.} For any periodic function $\varphi\in C^{\infty}((0, T)\times\Omega)$ about the space variable $\vec x$ which satisfies $\varphi(0)=\varphi(T)=0$, applying the operator $D_{x}^{\alpha}$ on both sides of (\ref{moment equation-2}) for $0\leq \alpha\leq 6$, we have
\begin{equation}\nonumber
\begin{aligned}
\Big|\int_{0}^{T}\int_{\Omega}D_{x}^{\alpha}\overline f\nabla_{x}\varphi \mathrm{d}\vec x\mathrm{d}t\Big|=&\lim_{\epsilon\rightarrow 0}\Big|\int_{0}^{T}\int_{\Omega}D_{x}^{\alpha}\overline f^{\epsilon}\nabla_{x}\varphi \mathrm{d}\vec x\mathrm{d}t\Big|\\=&3\lim_{\epsilon\rightarrow 0}\Big|\int_{0}^{T}\int_{\Omega}(\epsilon^{2}D_{x}^{\alpha}\vec{J}^{\epsilon}\partial_{t}\varphi+\epsilon\langle\vec w\otimes\vec wD_{x}^{\alpha}\hat{g^{\epsilon}}\rangle\nabla_{x}\varphi)\\&+\epsilon\langle\vec wD_{x}^{\alpha}(B(\theta^{\epsilon})-\overline f^{\epsilon})\rangle\varphi+\frac{1}{\epsilon}\langle\vec wD_{x}^{\alpha}(\overline f^{\epsilon}-f^{\epsilon})\rangle\varphi\mathrm{d}\vec x\mathrm{d}t\Big|
\\\leq&\sup_{\epsilon>0}\Big\|\frac{1}{\epsilon}\langle\vec wD_{x}^{\alpha}(\overline f^{\epsilon}-f^{\epsilon})\rangle\Big\|_{L^{2}((0, T)\times\Omega)}\|\varphi\|_{L^{2}((0, T)\times\Omega)}.
\end{aligned}
\end{equation}
The conclusion then follows from the fact that $\frac{1}{\epsilon}\langle\vec wD_{x}^{\alpha}(\overline f^{\epsilon}-f^{\epsilon})\rangle$ is bounded in $L^{2}((0, T)\times\Omega)$ and $\overline f\in L^{\infty}(0, T; W_{2}^{6}(\Omega))$.
\hfill$\Box$
\end{lem}

According to (iv) and (v), we can also obtain that $B(\theta^{\epsilon})\rightarrow B(\theta)$ strongly in $C^{0}([0, T];\\ W_{2}^{5}(\Omega))$ up to the same subsequence from the Sobolev embedding Theorem. Therefore, passing to the limit in (\ref{moment equation}) yields that

\begin{equation}\label{daidingfangcheng1}
\partial_{t}\overline f+\mathrm{div}_{x}\vec{J}=B(\theta)-\overline f.
\end{equation}

Dividing $\epsilon$ on both sides of $(\ref{research equations})_{1}$, multiplying the resulting equations with $\varphi$ for any periodic function $\varphi\in C^{\infty}([0, T]\times\Omega\times \mathbb{S}^{2})$ and integrating over $(0, T)\times\Omega\times \mathbb{S}^{2}$, we have
\begin{equation}\nonumber
-\int_{0}^{T}\int_{\Omega}\int_{\mathbb{S}^{2}}\vec{w}\cdot\nabla_{x}\overline f\varphi \mathrm{d}\vec{w}\mathrm{d}\vec{x}dt=\int_{0}^{T}\int_{\Omega}\int_{\mathbb{S}^{2}}\hat g\varphi \mathrm{d}\vec{w}\mathrm{d}\vec{x}dt
\end{equation}
as $\epsilon\rightarrow 0$ according to (i), (iii) in Lemma \ref{main result-6} and Lemma \ref{main result-3-bu-bu}. It follows that $\hat g=-\vec{w}\cdot\nabla_{x}\overline f$ a.e. in $[0, T]\times\Omega\times \mathbb{S}^{2}$. We conclude that
\begin{equation}\label{daidingfangcheng2}
\vec{J}(t, x)=-\langle \vec{w}\otimes\vec{w}\rangle\nabla_{x}\overline f=-\frac{1}{3}\nabla_{x}\overline f.
\end{equation}
Substituting (\ref{daidingfangcheng2}) into (\ref{daidingfangcheng1}), we have
\begin{equation}\label{couple 1-1}
\partial_{t}\overline f-\frac{1}{3}\Delta_{x}\overline f=B(\theta)-\overline f.
\end{equation}

Next, we pass to the limit in the equation for $V^{\epsilon}$ and get
\begin{equation}\label{couple 2}
A_{0}(V)\frac{\partial V}{\partial t}+\sum_{j=1}^{3}A_{j}(V)\frac{\partial V}{\partial x_{j}}=F(V, \overline f).
\end{equation}

The system is complemented with initial data $\overline f(0, \vec x)=\overline h(\vec x)$ and $V(0, \vec{x})=V^{0}(\vec x).$

 Finally, we will show the existence and uniqueness of system (\ref{equations about limit of overline f epsilon}). Moreover, we will prove that the strong convergence about $f^{\epsilon}$ and $V^{\epsilon}$ can be verified for any subsequence $\epsilon\rightarrow 0$. Up to a subsequence, $f^{\epsilon}$ and $V^{\epsilon}$ converge to $\overline f$ and $V$ strongly in $L^{2}(0, T; L^{2}(\mathbb{S}^{2}; W_{2}^{5}(\Omega)))$ and $C^{0}([0, T]; W_{2}^{5}(\Omega))$ respectively. Moreover, $f^{\epsilon}$ and $V^{\epsilon}$ converge in $L^{\infty}(0, T; L^{2}(\mathbb{S}^{2}; W_{2}^{6}(\Omega)))$ and $L^{\infty}(0, T; W_{2}^{6}(\Omega))$ respectively in the weak-$\ast$ sense. Note that $\overline f\in L^{2}(0, T; W_{2}^{7}(\Omega))$ according to Lemma \ref{main result-3-bu-bu}. Therefore, we obtain the existence of system (\ref{equations about limit of overline f epsilon}) as $\epsilon\rightarrow 0$ in (\ref{research equations}). Finally, we can obtain the uniqueness by the $L^{2}$ energy method combining with the Gronwall's inequality. Therefore, the above convergence holds without restricting to a subsequence.\hfill$\Box$

\section{Hilbert Expansion}\label{wu}
In order to get the strong convergence for $f^{\epsilon}$ in $L^{\infty}((0, T)\times\Omega\times \mathbb{S}^{2})$, we consider the Hilbert expansion in this section.

\subsection{Interior expansion}
First, we give the interior expansion of $\theta^{\epsilon}$ as follows:
\begin{equation}\label{expansion of theta epsilon}
\theta^{\epsilon}\sim\theta_{0}+\epsilon\theta_{1}+\epsilon^{2}\theta_{2}+\cdot\cdot\cdot.
\end{equation}
Then, we get
\begin{equation}\label{expansion of B(theta epsilon)}
\begin{split}
B(\theta^{\epsilon})&\sim(\theta_{0}+\epsilon\theta_{1}+\cdot\cdot\cdot)^{4}
\\&=\theta_{0}^{4}+4\theta_{0}^{3}\theta_{1}\epsilon+(4\theta_{2}\theta_{0}^{3}+6\theta_{1}^{2}\theta_{0}^{2})\epsilon^{2}+\cdot\cdot\cdot
\\&=B_{0}(\theta_{0})+\epsilon B_{1}(\theta_{0}; \theta_{1})+\epsilon^{2}B_{2}(\theta_{0}; \theta_{2})+\cdot\cdot\cdot
\end{split}
\end{equation}
where $B_{k}(\theta_{0}; \theta_{k}):=B_{k}(\theta_{0}; \theta_{1}; \theta_{2}; \theta_{3}; \cdot\cdot\cdot; \theta_{k})$ and $B_{0}(\theta_{0})=B(\theta_{0})$ for $k\geq 0$ and $k\in\mathbb{Z}$.
We also expand
\begin{equation}\label{expand rho epsilon}
\rho^{\epsilon}\sim\rho_{0}+\epsilon\rho_{1}+\epsilon^{2}\rho_{2}+\cdot\cdot\cdot
\end{equation}
and
\begin{equation}\label{expand u epsilon}
\vec u^{\epsilon}\sim\vec u_{0}+\epsilon\vec u_{1}+\epsilon^{2}\vec u_{2}+\cdot\cdot\cdot.
\end{equation}
Define the interior expansion of $f^{\epsilon}(t, \vec{x}, w)$ as
\begin{equation}\label{expansion about f epsilon including k term}
f^{\epsilon}(t, \vec x, \vec w)\sim \sum_{k=0}^{\infty}{\epsilon}^{k}f_{k}(t, \vec x, \vec w),
\end{equation}
where $f_{k}$ can be defined by comparing the order of $\epsilon$ via plugging (\ref{expansion of B(theta epsilon)}) and (\ref{expansion about f epsilon including k term}) into $(\ref{research equations})_{4}$.
It follows that
\begin{equation}\label{equality about f0 epsilon}
f_{0}-\overline f_{0}=0,\hspace{6.5cm}
\end{equation}
\begin{equation}\label{equality about f0 f1 epsilon}
f_{1}-\overline f_{1}=-\vec w\cdot\nabla_{x}f_{0},\hspace{4.9cm}
\end{equation}
\begin{equation}\label{equality about f0 f1 f2 epsilon}
f_{2}-\overline f_{2}=-\partial_{t}f_{0}-\vec w\cdot\nabla_{x}f_{1}-f_{0}+B(\theta_{0}),\hspace{1.2cm}
\end{equation}
\begin{equation}
f_{3}-\overline f_{3}=-\partial_{t}f_{1}-\vec w\cdot\nabla_{x}f_{2}-f_{1}+B_{1}(\theta_{0}; \theta_{1}),\hspace{0.5cm}
\end{equation}
\begin{equation}
f_{4}-\overline f_{4}=-\partial_{t}f_{2}-\vec w\cdot\nabla_{x}f_{3}-f_{2}+B_{2}(\theta_{0}; \theta_{2}),\hspace{0.5cm}
\end{equation}

$$\cdot\cdot\cdot$$

\begin{equation}\nonumber
f_{k+2}-\overline f_{k+2}=-\partial_{t}f_{k}-\vec w\cdot\nabla_{x}f_{k+1}-f_{k}+B_{k}(\theta_{0}; \theta_{k}),\hspace{-0.5cm}
\end{equation}
for $k\geq 2.$
Substituting (\ref{equality about f0 epsilon}) into (\ref{equality about f0 f1 epsilon}), we obtain
\begin{equation}\label{equality about f0 f1 epsilon 2}
f_{1}=\overline f_{1}-\vec w\cdot\nabla_{x}\overline f_{0}.
\end{equation}
Substituting (\ref{equality about f0 f1 epsilon 2}) into (\ref{equality about f0 f1 f2 epsilon}), we get
\begin{equation}\label{equality about f0 f1 f2 epsilon 2}
\overline f_{2}-f_{2}=\partial_{t}f_{0}+\vec w\cdot\nabla_{x}\overline f_{1}-(\vec w\cdot\nabla_{x})^{2}\overline f_{0}+f_{0}-B(\theta_{0})
\end{equation}
Integrating (\ref{equality about f0 f1 f2 epsilon 2}) over $\vec w\in \mathbb{S}^{2}$, we have
\begin{equation}\nonumber
\partial_{t}\overline f_{0}-\frac{1}{3}\Delta_{x}\overline f_{0}+\overline f_{0}=B(\theta_{0}).
\end{equation}
Similarly, we can get $f_{1}(t, \vec{x}, \vec{w}), f_{2}(t, \vec{x}, \vec{w})$ and $f_{k}(t, \vec{x}, \vec{w})(k\geq 3)$ as follows
\begin{equation}\nonumber
\partial_{t}\overline f_{1}-\frac{1}{3}\Delta_{x}\overline f_{1}+\overline f_{1}=B_{1}(\theta_{0}; \theta_{1}),
\end{equation}
\begin{equation}\nonumber
\partial_{t}\overline f_{2}-\frac{1}{3}\Delta_{x}\overline f_{2}+\overline f_{2}=B_{2}(\theta_{0}; \theta_{2})+\frac{1}{3}\Delta_{x}B(\theta_{0})+\frac{1}{4\pi}\int_{\mathbb{S}^{2}}(\vec{w}\cdot\nabla_{x})^{4}\overline f_{0}\mathrm{d}\vec w,
\end{equation}
\begin{equation}\nonumber
\begin{split}
\partial_{t}\overline f_{k}-\frac{1}{3}\Delta_{x}\overline f_{k}+\overline f_{k}=&B_{k}(\theta_{0}; \theta_{k})+\frac{1}{4\pi}\int_{\mathbb{S}^{2}}(\vec w\cdot\nabla_{x})(f_{k-1}+\partial_{t}f_{k-1})\mathrm{d}\vec w\\&-\frac{1}{4\pi}\int_{\mathbb{S}^{2}}(\vec w\cdot\nabla_{x})^{2}f_{k-2}\mathrm{d}\vec w-\frac{1}{4\pi}\int_{\mathbb{S}^{2}}(\vec w\cdot\nabla_{x})^{3}f_{k-1}\mathrm{d}\vec w\\&-\frac{1}{4\pi}\int_{\mathbb{S}^{2}}(\vec w\cdot\nabla_{x})^{2}\partial_{t}f_{k-2}\mathrm{d}\vec w+\frac{1}{3}\Delta_{x}(B_{k-2}(\theta_{0}; \theta_{k-2})).
\end{split}
\end{equation}
We canalso  obtain the equation about $\rho_{k}$, $\vec u_{k}$ and $\theta_{k}$ via plugging (\ref{expansion of theta epsilon}), (\ref{expansion of B(theta epsilon)}), (\ref{expand rho epsilon}), (\ref{expand u epsilon}) and (\ref{expansion about f epsilon including k term}) into $(\ref{research equations})_{1}$, $(\ref{research equations})_{2}$ and $(\ref{research equations})_{3}$. Drirect calculations shows that $(\rho_0,\vec{u}_0,\theta_0)$ satisfies
\begin{equation}\nonumber
\partial_{t}\rho_{0}+\mathrm{div}_{x}(\vec{u}_{0}\theta_{0})=0,
\end{equation}
\begin{equation}\nonumber
\partial_{t}(\rho_{0}\vec{u}_{0})+\mathrm{div}_{x}(\rho_{0}\vec{u}_{0}\otimes\vec{u}_{0}+\rho_{0}\theta_{0}),
=-\frac{1}{3}\nabla_{x}\overline f_{0},
\end{equation}
\begin{equation}\nonumber
\partial_{t}(\rho_{0}\theta_{0})+\mathrm{div}_{x}(\rho_{0}\theta_{0}\vec{u}_{0})+\rho_{0}\theta_{0}\mathrm{div}_{x}\vec{u}_{0}
=\overline f_{0}-B(\theta_{0})+\frac{1}{3}\nabla_{x}\overline f_{0}\cdot\vec u_{0}.
\end{equation}
 $(\rho_1,\vec{u}_1,\theta_1)$ satisfies
\begin{equation}\nonumber
\begin{aligned}
\frac{\partial\rho_{1}}{\partial t}+\vec{u}_{0}\cdot\nabla_{x}\rho_{1}+\rho_{0}\mathrm{div}_{x}\vec{u}_{1}=-\mathrm{div}_{x}\vec{u}_{0}\rho_{1}-\vec{u}_{1}\cdot\nabla_{x}\rho_{0},
\end{aligned}
\end{equation}
\begin{equation}\nonumber
\begin{aligned}
\frac{\partial\vec{u}_{1}}{\partial t}+\vec{u}_{0}\cdot\nabla_{x}\vec{u}_{1}+\nabla_{x}\theta_{1}+\frac{\theta_{0}}{\rho_{0}}\nabla_{x}\rho_{1}
=&-\frac{\rho_{1}}{\rho_{0}}\frac{\partial \vec{u}_{0}}{\partial t}-\vec{u}_{1}\cdot\nabla_{x}\vec{u}_{0}-\frac{\rho_{1}}{\rho_{0}}\vec{u}_{0}\cdot\nabla_{x}\vec{u}_{0}
\\&-\frac{\rho_{1}}{\rho_{0}}\nabla_{x}\theta_{0}-\frac{\theta_{1}}{\rho_{0}}\nabla_{x}\rho_{0}-\frac{1}{3}\nabla_{x}\overline f_{1},
\end{aligned}
\end{equation}
\begin{equation}\nonumber
\begin{aligned}
\frac{\partial\theta_{1}}{\partial t}+\vec{u}_{0}\cdot\nabla_{x}\theta_{1}+\theta_{0}\mathrm{div}_{x}\vec{u}_{1}=&\frac{\overline f_{1}-B_{1}(\theta_{0}; \theta_{1})}{\rho_{0}}-\frac{\rho_{1}}{\rho_{0}}\frac{\partial\theta_{0}}{\partial t}-\vec{u}_{1}\cdot\nabla_{x}\theta_{0}-\frac{\rho_{1}}{\rho_{0}}\vec{u}_{0}\cdot\nabla_{x}\theta_{0}
\\&-\theta_{1}\mathrm{div}_{x}\vec{u}_{0}-\frac{\rho_{1}}{\rho_{0}}\theta_{0}\mathrm{div}_{x}\vec{u}_{0}+\frac{1}{3}\nabla_{x}\overline f_{1}\cdot\vec{u}_{0},
\end{aligned}
\end{equation}

and  $(\rho_2,\vec{u}_2,\theta_2)$ satisfies
\begin{equation}\nonumber
\begin{aligned}
&\frac{\partial\rho_{2}}{\partial t}+\vec{u}_{0}\cdot\nabla_{x}\rho_{2}+\rho_{0}\mathrm{div}_{x}\vec{u}_{2}=-\vec{u}_{2}\cdot\nabla_{x}\vec{u}_{0}
-\vec{u}_{1}\cdot\nabla_{x}\rho_{1}-\rho_{2}\mathrm{div}_{x}\vec{u}_{0}-\rho_{1}\mathrm{div}_{x}\vec{u}_{1},
\end{aligned}
\end{equation}
\begin{equation}\nonumber
\begin{aligned}
\frac{\partial\vec{u}_{2}}{\partial t}+\vec{u}_{0}\cdot\nabla_{x}\vec{u}_{2}+\nabla_{x}\theta_{2}+\frac{\theta_{0}}{\rho_{0}}\nabla_{x}\rho_{2}=&-\frac{\rho_{1}}{\rho_{0}}\frac{\partial\vec{u}_{1}}{\partial t}-\frac{\rho_{2}}{\rho_{0}}\frac{\partial\vec{u}_{0}}{\partial t}-\vec{u}_{1}\cdot\nabla_{x}\vec{u}_{1}+\vec{u}_{2}\cdot\nabla_{x}\vec{u}_{0}
\\&-\frac{\rho_{1}}{\rho_{0}}(\vec{u}_{0}\cdot\nabla_{x}\vec{u}_{1}+\vec{u}_{1}\cdot\nabla_{x}\vec{u}_{0})
-\frac{\rho_{2}}{\rho_{0}}\vec{u}_{0}\cdot\nabla_{x}\vec{u}_{0}-\frac{\rho_{1}}{\rho_{0}}\nabla_{x}\theta_{1}\\&-\frac{\rho_{2}}{\rho_{0}}\nabla_{x}\theta_{0}
-\frac{\theta_{1}}{\rho_{0}}\nabla_{x}\rho_{1}-\frac{\theta_{2}}{\rho_{0}}\nabla_{x}\rho_{0}+\frac{1}{3\rho_{0}}(2\nabla_{x}\overline f_{0}-\nabla_{x}\overline f_{2}\\&+\partial_{t}\nabla_{x}\overline f_{0}+\nabla_{x}\overline f_{1}-\nabla_{x}B(\theta_{0}))
\end{aligned}
\end{equation}
and
\begin{equation}\nonumber
\begin{aligned}
\frac{\partial\theta_{2}}{\partial t}+\vec{u}_{0}\cdot\nabla_{x}\theta_{2}+\theta_{0}\mathrm{div}_{x}\vec{u}_{2}=&\frac{\overline f_{2}-B_{2}(\theta_{0}; \theta_{2})}{\rho_{0}}-\frac{\rho_{1}}{\rho_{0}}\frac{\partial\theta_{1}}{\partial t}-\frac{\rho_{2}}{\rho_{0}}\frac{\partial\theta_{0}}{\partial t}
-\frac{\rho_{2}}{\rho_{0}}\vec{u}_{0}\cdot\nabla_{x}\theta_{0}\\&-\frac{\rho_{1}}{\rho_{0}}(\vec{u}_{1}\cdot\nabla_{x}\theta_{0}+\vec{u}_{0}\cdot\nabla_{x}\theta_{1})
-(\vec{u}_{1}\cdot\nabla_{x}\theta_{1}+\vec{u}_{2}\cdot\nabla_{x}\theta_{0})
\\&-\frac{\rho_{2}\theta_{0}}{\rho_{0}}\mathrm{div}\vec{u}_{0}-\frac{\rho_{1}}{\rho_{0}}(\theta_{1}\mathrm{div}\vec{u}_{0}+\theta_{0}\mathrm{div}\vec{u}_{1})
-(\theta_{1}\mathrm{div}\vec{u}_{1}+\theta_{2}\mathrm{div}\vec{u}_{0})\\&-\frac{1}{3\rho_{0}}(2\nabla_{x}\overline f_{0}-\nabla_{x}\overline f_{2}+\partial_{t}\nabla_{x}\overline f_{0}+\nabla_{x}\overline f_{1}-\nabla_{x}B(\theta_{0}))\cdot\vec{u}_{0}
\end{aligned}
\end{equation}
for $k=0, 1, 2$. We don't give the equalities about $k\geq 3$ which can't influence our final result.

\subsection{Initial layer expansion.}
In order to determine the initial condition for $f_{k}$ and $\theta_{k}$, we resort to the initial layer expansion. We define the new variable $\tau$ by making the scaling transform for $(f^{\epsilon}(t, \vec{x}, \vec w), \rho^{\epsilon}(t, \vec x), \vec{u}^{\epsilon}(t, \vec x), \theta^{\epsilon}(t, \vec x))\rightarrow (f^{\epsilon}(\tau, \vec{x}, \vec w), \rho^{\epsilon}(\tau, \vec x), \vec{u}^{\epsilon}(\tau, \vec x), \\\theta^{\epsilon}(\tau, \vec x))$ with $\tau\in [0, \infty)$ as
\begin{equation}\nonumber
\tau=\frac{t}{\epsilon^{2}},
\end{equation}
which implies
\begin{equation}\nonumber
\frac{\partial f^{\epsilon}}{\partial t}=\frac{1}{\epsilon^{2}}\frac{\partial f^{\epsilon}}{\partial\tau}, \frac{\partial \rho^{\epsilon}}{\partial t}=\frac{1}{\epsilon^{2}}\frac{\partial \rho^{\epsilon}}{\partial\tau}, \frac{\partial \vec{u}^{\epsilon}}{\partial t}=\frac{1}{\epsilon^{2}}\frac{\partial \vec{u}^{\epsilon}}{\partial\tau}, \frac{\partial \theta^{\epsilon}}{\partial t}=\frac{1}{\epsilon^{2}}\frac{\partial \theta^{\epsilon}}{\partial\tau}.
\end{equation}
In this new variable, the system $(\ref{research equations})$ can be rewritten as
\begin{equation}\label{research equation in the initial layer}\left\{
\begin{split}
&\frac{1}{\epsilon^{2}}\frac{\partial \rho^{\epsilon}}{\partial\tau}+\vec{u}^{\epsilon}\cdot\nabla_{x}\rho^{\epsilon}+\rho^{\epsilon}\mathrm{div}_{x}\vec{u}^{\epsilon}=0,~~\mathrm{in} \ \ (0,T)\times \Omega,\\
&\frac{1}{\epsilon^{2}}\frac{\partial \vec u^{\epsilon}}{\partial\tau}+\vec{u}^{\epsilon}\cdot\nabla_{x}\vec{u}^{\epsilon}+\nabla_{x}\theta^{\epsilon}
+\frac{\theta^{\epsilon}}{\rho^{\epsilon}}\nabla_{x}\rho^{\epsilon}=\frac{1}{\rho^{\epsilon}}\Big\langle\Big(\frac{1}{\epsilon}+\epsilon\Big)\vec{w}(f^{\epsilon}-\overline f^{\epsilon})\Big\rangle,~~\mathrm{in} \ \ (0,T)\times \Omega,\\
&\frac{1}{\epsilon^{2}}\frac{\partial \theta^{\epsilon}}{\partial\tau}+\vec{u}^{\epsilon}\cdot\nabla_{x}\theta^{\epsilon}+\theta^{\epsilon}\mathrm{div}_{x}\vec{u}^{\epsilon}
=\frac{1}{\rho^{\epsilon}}(\overline f^{\epsilon}-B(\theta^{\epsilon}))-\frac{1}{\rho^{\epsilon}}\Big(\frac{1}{\epsilon}+\epsilon\Big)\langle\vec{w}(f^{\epsilon}-\overline f^{\epsilon})\rangle\cdot\vec u^{\epsilon},\\
&\frac{\partial f^{\epsilon}}{\partial\tau}+\epsilon \vec w\cdot \nabla_{x}f^{\epsilon}+f^{\epsilon}-\overline f^{\epsilon}+{\epsilon}^{2}f^{\epsilon}={\epsilon}^{2}B(\theta^{\epsilon})\ \ \mathrm{in} \ \ (0, \infty)\times \Omega\times \mathbb{S}^{2},\\
&\rho^{\epsilon}(0, \vec{x})=\rho^{0}(\vec x),\;\;\vec{u}^{\epsilon}(0, \vec{x})=\vec{u}^{0}(\vec x),\;\;\theta^{\epsilon}(0, \vec{x})=\theta^{0}(\vec{x})~~\mathrm{in}  \ \ \Omega,\\
&f^{\epsilon}(0, \vec{x}, \vec{w})=h(\vec{x}, \vec{w}),~~\mathrm{in} \ \ \Omega\times \mathbb{S}^{2}.\\
\end{split}\right.
\end{equation}
We define the initial layer expansion as follows:
\begin{equation}\label{expansion about the initial layer}
f_{I}^{\epsilon}\thicksim \sum_{k=0}^{\infty}\epsilon^{k}f_{I, k}(\tau, \vec{x}, \vec{w}),\;\;\rho_{I}^{\epsilon}\thicksim\sum_{k=0}^{\infty}\epsilon^{k}\rho_{I, k}(\tau, \vec{x}),
\end{equation}
\begin{equation}\label{expansion about u epsilon initial layer}
\vec u_{I}^{\epsilon}\thicksim\sum_{k=0}^{\infty}\epsilon^{k}\vec u_{I, k}(\tau, \vec{x}),\;\;
\theta_{I}^{\epsilon}\thicksim\sum_{k=0}^{\infty}\epsilon^{k}\theta_{I, k}(\tau, \vec{x}).
\end{equation}
Then, $\rho_{I, k}, \vec u_{I, k}, \theta_{I, k}$ and $f_{I, k}$ can be determined by comparing the order of $\epsilon$ via plugging $f^{\epsilon}=\sum_{k=0}^{\infty}\epsilon^{k}f_{k}+\epsilon^{k}f_{I, k}$, $\rho^{\epsilon}=\sum_{k=0}^{\infty}\epsilon^{k}\rho_{k}+\epsilon^{k}\rho_{I, k}$, $\vec{u}^{\epsilon}=\sum_{k=0}^{\infty}\epsilon^{k}\vec{u}_{k}+\epsilon^{k}\vec{u}_{I, k}$ and $\theta^{\epsilon}=\sum_{k=0}^{\infty}\epsilon^{k}\theta_{k}+\epsilon^{k}\theta_{I, k}$ into (\ref{research equation in the initial layer}). Thus, we have
\begin{equation}\label{equality about f I0}
\frac{\partial f_{I, 0}}{\partial\tau}+f_{I, 0}-\overline f_{I, 0}=0,
\end{equation}
\begin{equation}\label{equality about f I0 f I1}
\vec w\cdot\nabla_{x}f_{I, 0}+\frac{\partial f_{I, 1}}{\partial\tau}+f_{I, 1}-\overline f_{I, 1}=0,\ \
\end{equation}and
\begin{equation}\label{equality about f I0 f I1 f I2}
f_{I, 0}+\vec w\cdot\nabla_{x}f_{I, 1}+\frac{\partial f_{I, 2}}{\partial\tau}+f_{I, 2}-\overline f_{I, 2}=B(\theta_{I, 0}).
\end{equation}
The following analysis reveals the equation satisfied by $f_{I, 0}$:

Integrating (\ref{equality about f I0}) over $\vec w\in \mathbb{S}^{2}$, we have
\begin{equation}\nonumber
\frac{\partial\overline f_{I, 0}}{\partial\tau}=0,
\end{equation}
which further implies that
\begin{equation}\nonumber
\overline f_{I, 0}(\tau, \vec{x})=\overline f_{I, 0}(0, \vec{x}).
\end{equation}
Therefore we can deduce from (\ref{equality about f I0})

\begin{equation}\nonumber
\begin{split}
f_{I, 0}(\tau, \vec{x}, \vec{w})&=e^{-\tau}f_{I, 0}(0, \vec{x}, \vec{w})+\int_{0}^{\tau}\overline f_{I, 0}(s, \vec{x})e^{s-\tau}\mathrm{d}s\\&=e^{-\tau}f_{I, 0}(0, \vec{x}, \vec{w})+(1-e^{-\tau})\overline f_{I, 0}(0, \vec{x}).
\end{split}
\end{equation}
This means that we have

\begin{equation}\nonumber\label{equation}\left\{
\begin{split}
&\frac{\partial\overline f_{I, 0}}{\partial\tau}=0,\\
&f_{I, 0}(\tau, \vec{x}, \vec{w})=e^{-\tau}f_{I, 0}(0, \vec{x}, \vec{w})+(1-e^{-\tau})\overline f_{I, 0}(0, \vec{x}).\\
\end{split}\right.
\end{equation}
Similarly, we can derive that $f_{I, 1}$ and $f_{I, 2}$ satisfying
\begin{equation}\nonumber\label{equation}\left\{
\begin{split}
&\frac{\partial\overline f_{I, 1}}{\partial\tau}=-\frac{1}{4\pi}\int_{\mathbb{S}^{2}}\vec w\cdot\nabla_{x}f_{I, 0}\mathrm{d}\vec{w},\\
&f_{I, 1}(\tau, \vec{x}, \vec{w})=e^{-\tau}f_{I, 0}(0, \vec{x}, \vec{w})+\int_{0}^{\tau}(\overline f_{I, 1}-\vec w\cdot\nabla_{x}f_{I, 0})(s, \vec{x}, \vec{w})e^{s-\tau}\mathrm{d}s\\
\end{split}\right.
\end{equation}
and
\begin{equation}\nonumber\label{equation}\left\{
\begin{split}
&\frac{\partial\overline f_{I, 2}}{\partial\tau}=\frac{1}{4\pi}\int_{\mathbb{S}^{2}}(B(\theta^{0}+\theta_{I, 0})-B(\theta^{0})-\vec w\cdot\nabla_{x}f_{I, 1}-f_{I, 0})\mathrm{d}\vec{w},\\
&f_{I, 2}(\tau, \vec{x}, \vec{w})=e^{-\tau}f_{I, 2}(0, \vec{x}, \vec{w})+\int_{0}^{\tau}(B(\theta^{0}+\theta_{I, 0})-B(\theta^{0})-\vec w\cdot\nabla_{x}f_{I, 1}-f_{I, 0})(s, \vec{x}, \vec{w})e^{s-\tau}\mathrm{d}s.\\
\end{split}\right.
\end{equation}
We can also get
\begin{equation}\nonumber
\frac{\partial \rho_{I, 0}}{\partial\tau}=0,\;\; \frac{\partial \vec u_{I, 0}}{\partial\tau}=0,\;\; \frac{\partial \theta_{I, 0}}{\partial\tau}=0,
\end{equation}
\begin{equation}\nonumber
\frac{\partial \rho_{I, 1}}{\partial\tau}=0,\;\; \frac{\partial \vec u_{I, 1}}{\partial\tau}=\frac{1}{\rho^{0}}\langle\vec w(f_{I, 0}-\overline f_{I, 0})\rangle,\;\; \frac{\partial \theta_{I, 1}}{\partial\tau}=0,
\end{equation}
\begin{equation}\nonumber
\frac{\partial \rho_{I, 2}}{\partial\tau}=-(\vec{u}^{0}+\vec{u}_{I, 0})\cdot\nabla_{x}(\rho^{0}+\rho_{I, 0})+\vec{u}^{0}\cdot\nabla_{x}\rho^{0}-(\rho^{0}+\rho_{I, 0})\mathrm{div}_{x}(\vec{u}^{0}+\vec{u}_{I, 0})+\rho_{0}\mathrm{div}_{x}\vec{u}^{0},
\end{equation}
\begin{equation}\nonumber
\begin{aligned}
\frac{\partial \vec u_{I, 2}}{\partial\tau}=&-(\vec{u}^{0}+\vec{u}_{I, 0})\cdot\nabla_{x}(\vec{u}^{0}+\vec{u}_{I, 0})+\vec{u}^{0}\cdot\nabla_{x}\vec{u}^{0}-\nabla_{x}\theta_{I, 0}\\&-\frac{1}{\rho^{0}}\Big\{(\theta^{0}+\theta_{I, 0})\nabla_{x}(\rho^{0}+\rho^{I, 0})-\theta^{0}\nabla_{x}\rho^{0}\Big\}+\frac{1}{\rho^{0}}\langle\vec w(f_{I, 0}-\overline f_{I, 0})\rangle,
\end{aligned}
\end{equation} and
\begin{equation}\nonumber
\begin{aligned}
\frac{\partial \theta_{I, 2}}{\partial\tau}=&-(\vec{u}^{0}+\vec{u}_{I, 0})\cdot\nabla_{x}(\theta^{0}+\theta_{I, 0})+\vec{u}^{0}\cdot\nabla_{x}\theta^{0}-(\rho^{0}+\theta_{I, 0})\mathrm{div}_{x}(\vec{u}^{0}+\vec{u}_{I, 0})+\theta^{0}\mathrm{div}_{x}\vec{u}^{0}\\&+\frac{1}{\rho^{0}}(f_{I, 0}-B(\theta_{I, 0}))-\frac{1}{\rho^{0}}(\langle\vec w(f_{1}^{0}+f_{I, 1}-\overline f_{1}^{0}-\overline f_{I, 1})\rangle\cdot(\vec{u}^{0}+\vec{u}_{I, 1})\\&+\frac{1}{\rho^{0}}(\langle\vec w(f_{1}^{0}-\overline f_{1}^{0})\rangle\cdot\vec{u}^{0}-\frac{1}{\rho^{0}}(\langle\vec w(f^{0}_{0}+f_{I, 0}-\overline f_{0}^{0}-\overline f_{I, 0})\rangle\cdot(\vec{u}_{1}+\vec{u}_{I, 0})\\&+\frac{1}{\rho_{0}}(\langle\vec w(f_{0}^{0}-\overline f_{0}^{0})\rangle\cdot\vec{u}_{1},
\end{aligned}
\end{equation}
where we have used the Taylor expansion, i.e.
\begin{equation}\nonumber
\rho_{0}(t, \vec x)=\rho^{0}+\sum_{l=1}^{b}\frac{1}{l!}\partial_{t}^{l}\rho_{0}(t, \vec x)|_{t=0}t^{l}+\frac{t^{b+1}}{(b+1)!}\partial_{t}^{b+1}\tilde{\rho_{0}}.
\end{equation}

\begin{equation}\nonumber
\theta_{0}(t, \vec x)=\theta^{0}+\sum_{l=1}^{b}\frac{1}{l!}\partial_{t}^{l}\theta_{0}(t, \vec x)|_{t=0}t^{l}+\frac{t^{b+1}}{(b+1)!}\partial_{t}^{b+1}\tilde{\theta_{0}}.
\end{equation}

\begin{equation}\nonumber
\vec u_{i}(t, \vec x)=\vec u_{i}^{0}+\sum_{l=1}^{b}\frac{1}{l!}\partial_{t}^{l}\vec u_{i}(t, \vec x)|_{t=0}t^{l}+\frac{t^{b+1}}{(b+1)!}\partial_{t}^{b+1}\tilde{\vec u_{i}}.
\end{equation}

\begin{equation}\nonumber
f_{i}(t, \vec x, \vec w)=f_{i}^{0}+\sum_{l=1}^{b}\frac{1}{l!}\partial_{t}^{l}f_{i}(t, \vec x)|_{t=0}t^{l}+\frac{t^{b+1}}{(b+1)!}\partial_{t}^{b+1}\tilde{f_{i}}.
\end{equation}
Here $\partial_{t}^{b+1}\tilde{\rho_{0}}=\partial_{t}^{b+1}\rho_{0}(t', \vec x), \partial_{t}^{b+1}\tilde{\theta_{0}}=\partial_{t}^{b+1}\theta_{0}(t', \vec x), \partial_{t}^{b+1}\tilde{\vec u_{i}}=\partial_{t}^{b+1}\vec u_{i}(t', \vec x), \partial_{t}^{b+1}\tilde{f_{i}}=\partial_{t}^{b+1}\tilde{f_{i}}\\(t', \vec x), f_{i}^{0}=f_{i}(0, \vec x, \vec w), \vec u_{i}^{0}=\vec u_{i}(0, \vec x)$ where $i=0, 1,$ and $t'\in (0, t).$ Furthermore, it suffices to choose $b=1$ for our problem. We don't give the equalities about $k\geq 3$ which can't influence our final result.

\subsection{Construction of asymptotic expansion.}
The bridge between the interior solution with the initial layer is the initial condition of (\ref{research equations}). First, we have the following relations
\begin{equation}\nonumber
f_{0}(0, \vec{x}, \vec{w})+f_{I, 0}(0, \vec{x}, \vec{w})=h(\vec{x}, \vec{w}),\hspace{0.5 cm}
\end{equation}
\begin{equation}\nonumber
\rho_{0}(0, \vec{x})+\rho_{I, 0}(0, \vec{x})=\rho^{0}(\vec{x}),\hspace{-0.2 cm}
\end{equation}
\begin{equation}\nonumber
\vec{u}_{0}(0, \vec{x})+\vec{u}_{I, 0}(0, \vec{x})=\vec{u}^{0}(\vec{x}),\hspace{-0.2 cm}
\end{equation}
\begin{equation}\nonumber
\theta_{0}(0, \vec{x})+\theta_{I, 0}(0, \vec{x})=\theta^{0}(\vec{x}),\hspace{-0.2 cm}
\end{equation}
\begin{equation}\nonumber
f_{k}(0, \vec{x}, \vec{w})+f_{I, k}(0, \vec{x}, \vec{w})=0,k\geq 1,
\end{equation}
\begin{equation}\nonumber
\rho_{k}(0, \vec{x})+\rho_{I, k}(0, \vec{x})=0,k\geq 1,\hspace{-1 cm}
\end{equation}
\begin{equation}\nonumber
\vec{u}_{k}(0, \vec{x})+\vec{u}_{I, k}(0, \vec{x})=0,k\geq 1,\hspace{-1 cm}
\end{equation}
\begin{equation}\nonumber
\theta_{k}(0, \vec{x})+\theta_{I, k}(0, \vec{x})=0,k\geq 1.\hspace{-1 cm}
\end{equation}
The construction of $f_{k}, f_{I, k}$, $\theta_{k}$ and $\theta_{I, k}$ are as follows:

\noindent\textbf{Step 1.} Construction of zeroth-order terms.

The zeroth-order initial layer $f_{I, 0}$ is defined as
\begin{equation}\label{equations about f I0}\left\{
\begin{split}
&f_{I, 0}(\tau, \vec{x}, \vec{w})=\mathcal{F}_{0}(\tau, \vec{x}, \vec{w})-\mathcal{F}_{0}(\infty, \vec{x}),\\
&\frac{\partial\overline{\mathcal{F}_{0}}}{\partial\tau}=0,\\
&\mathcal{F}_{0}(\tau, \vec{x}, \vec{w})=e^{-\tau}\mathcal{F}_{0}(0, \vec{x}, \vec{w})+(1-e^{-\tau})\overline{ \mathcal{F}_{0}}(0, \vec{x}),\\
&\mathcal{F}_{0}(0, \vec{x}, \vec{w})=h(\vec{x}, \vec{w}),\\
&\lim_{\tau\rightarrow\infty}\mathcal{F}_{0}(\tau, \vec{x}, \vec{w})=\mathcal{F}_{0}(\infty, \vec{x}).
\end{split}\right.
\end{equation}
After a careful calculation, we can obtain $\mathcal{F}_{0}(\tau, \vec{x}, \vec{w})=e^{-\tau}h(\vec x, \vec w)+(1-e^{-\tau})\overline h(\vec x)$ which implies $\mathcal{F}_{0}(\tau, \vec{x}, \vec{w})\in C^{\infty}([0, \infty); L^{\infty}(\mathbb{S}^{2}; W_{2}^{6}(\Omega)))$ and $\mathcal{F}_{0}(\infty, \vec{x})=\overline h(\vec{x})\in W_{2}^{6}(\Omega)$.

We define the equation about $\rho_{I, 0}, \vec{u}_{I, 0}, \theta_{I, 0}$ as follows
\begin{equation}\label{equation about theta I0}\left\{
\begin{split}
&(\rho_{I, 0}, \vec{u}_{I, 0}, \theta_{I, 0})=(\hat{\rho_{0}}-\hat{\rho_{0}}(\infty), \hat{\vec {u}_{0}}-\hat{\vec {u}_{0}}(\infty)), \hat{\theta_{0}}-\hat{\theta_{0}}(\infty)),\\
&\frac{\partial \hat{\rho_{0}}}{\partial\tau}=0,\;\;\frac{\partial \hat{\vec u_{0}}}{\partial\tau}=0,\;\;
\frac{\partial \hat{\theta_{0}}}{\partial\tau}=0,\\
&(\hat{\rho_{0}}(0, \vec x), \hat{\vec u_{0}}(0, \vec x), \hat{\theta_{0}}(0, \vec x))=(\rho^{0}(\vec x), \vec {u}^{0}(\vec x), \theta^{0}(\vec x)),\\
&\lim_{\tau\rightarrow\infty}(\hat{\rho_{0}}(\tau, \vec x), \hat{\vec{u}}_{0}(\tau, \vec x), \hat{\theta_{0}}(\tau, \vec x))=(\hat{\rho_{0}}(\infty, \vec x), \hat{\vec{u}}_{0}(\infty, \vec x), \hat{\theta_{0}}(\infty, \vec x)).\\
\end{split}\right.
\end{equation}
We can obtain $(\hat{\rho_{0}}(t, \vec x), \hat{\vec u_{0}}(t, \vec x), \hat{\theta_{0}}(t, \vec x))\equiv(\rho^{0}(\vec x), \vec {u}^{0}(\vec x), \theta^{0}(\vec x))$ implying $(\rho_{I, 0}, \vec{u}_{I, 0}, \theta_{I, 0})\\=(0, 0, 0)$. Then, we have
\begin{equation}\label{equations about f0 epsilon and theta}\left\{
\begin{split}
&f_{0}(t, \vec{x}, \vec{w})=\overline f_{0}(t, \vec{x}),\\
&\partial_{t}\rho_{0}+\mathrm{div}_{x}(\rho_{0}\vec{u}_{0})=0,~~\mathrm{in} \ \ (0,T)\times \Omega,\\
&\partial_{t}(\rho_{0}\vec{u}_{0})+\mathrm{div}_{x}(\rho_{0}\vec{u}_{0}\otimes\vec{u}_{0}+\rho_{0}\theta_{0}\mathbb{I}_{3\times 3})
=-\frac{1}{3}\nabla_{x}\overline f_{0},~~\mathrm{in} \ \ (0,T)\times \Omega,\\
&\partial_{t}(\rho_{0}\theta_{0})+\mathrm{div}_{x}(\rho_{0}\theta_{0}\vec{u}_{0})+\rho_{0}\theta_{0}\mathrm{div}_{x}\vec{u}_{0}
=\overline f_{0}-B(\theta_{0})+\frac{1}{3}\nabla_{x}\overline f_{0}\cdot\vec{u}_{0},~~\mathrm{in} \ \ (0,T)\times \Omega, \\
&\partial_{t}\overline f_{0}-\frac{1}{3}\Delta_{x}\overline f_{0}+\overline f_{0}=B(\theta_{0})\ \ \mathrm{in}\ \ \ (0,T]\times \Omega,\\
&\rho_{0}(0, \vec{x})=\rho^{0}(\vec x),\;\;\vec{u}_{0}(0, \vec{x})=\vec{u}^{0}(\vec x)~~\mathrm{in}  \ \ \Omega,\\
&\theta_{0}(0, \vec{x})=\theta^{0}(\vec{x}),\;\;\overline f_{0}(0, \vec{x})=\mathcal{F}_{0}(\infty, \vec x)~~\mathrm{in}  \ \ \Omega,
\end{split}\right.
\end{equation}
Since $\mathcal{F}_{0}(\infty, x)=\overline h(\vec x)$, we have $(\overline f_{0}, \rho_{0}, \vec u_{0}, \theta_{0})=(\overline f, \rho, \vec u, \theta)$ by the uniqueness of system (\ref{equations about limit of overline f epsilon}). So, we obtain that $(\overline f_{0}, \rho_{0}, \vec u_{0}, \theta_{0})\in (L^{\infty}(0, T; W_{2}^{6}(\Omega))\cap L^{2}(0, T; W_{2}^{7}(\Omega)))\times L^{\infty}(0, T; W_{2}^{6}(\Omega))\times L^{\infty}(0, T; W_{2}^{6}(\Omega))\times L^{\infty}(0, T; W_{2}^{6}(\Omega))$.

\begin{rem}\label{pinglun 1}
Note that $(\overline f_{0}, \rho_{0}, \vec u_{0}, \theta_{0})=(\overline f, \rho, \vec u, \theta)$. Then, we have $\lim_{\epsilon\rightarrow 0}\|\theta^{\epsilon}-\theta_{0}\|_{C^{0}([0, T]; W_{2}^{5}(\Omega))}\\=0$. This convergence result is very important for us to prove (\ref{jielunbudengshi 1-jia}) in Theorem \ref{definition weak solution-1}.
\end{rem}

\noindent\textbf{Step 2.} Construction of first-order terms.

The first-order initial layer $f_{I, 1}$ is defined as
\begin{equation}\label{equation about f I1}\left\{
\begin{split}
&f_{I, 1}(\tau, \vec{x}, \vec{w})=\mathcal{F}_{1}(\tau, \vec{x}, \vec{w})-\mathcal{F}_{1}(\infty, \vec{x}),\\
&\frac{\partial\overline {\mathcal{F}_{1}}}{\partial\tau}=-\frac{1}{4\pi}\int_{S^{2}}\vec{w}\cdot\nabla_{x}f_{I, 0}\mathrm{d}\vec{w},\\
&\mathcal{F}_{1}(\tau, \vec{x}, \vec{w})=e^{-\tau}\mathcal{F}_{1}(0, \vec{x}, \vec{w})+\int_{0}^{\tau}(\overline {\mathcal{F}_{1}}- \vec{w}\cdot\nabla_{x}f_{I, 0})e^{s-\tau}\mathrm{d}s.\\
&\mathcal{F}_{1}(0, \vec{x}, \vec{w})=\vec{w}\cdot\nabla_{x}f_{0},\\
&\lim_{\tau\rightarrow\infty}\mathcal{F}_{1}(\tau, \vec{x}, \vec{w})=\mathcal{F}_{1}(\infty, \vec{x}).
\end{split}\right.
\end{equation}

Since $\vec{w}\cdot\nabla_{x}f_{I, 0}\in L^{\infty}((0, \infty)\times\Omega\times \mathbb{S}^{2})$ and $\mathcal{F}_{1}(0, \vec{x}, \vec{w})=\vec{w}\cdot\nabla_{x}\overline h\in L^{\infty}(\mathbb{S}^{2}; W_{2}^{5}(\Omega))$, we derive $\overline{\mathcal{F}_{1}}(\tau, x)=O(1+e^{-\tau})$. Since $\mathcal{F}_{1}(0, \vec{x}, \vec{w})\in L^{\infty}(\mathbb{S}^{2}; W_{2}^{6}(\Omega))$, we have $\mathcal{F}_{1}=O(1+e^{-\tau})$ and $\mathcal{F}_{1}(\infty, \vec x)\in W_{2}^{5}(\Omega)$. Therefore,  We can obtain $f_{I, 1}\in C^{\infty}([0, \infty); L^{\infty}(\mathbb{S}^{2}; W_{2}^{5}(\Omega)))$.

We define the equation about $\rho_{I, 1}, \vec{u}_{I, 1}, \theta_{I, 1}$ as follows
\begin{equation}\label{equation about theta I1}\left\{
\begin{split}
&(\rho_{I, 1}, \vec{u}_{I, 1}, \theta_{I, 1})=(\hat{\rho_{1}}-\hat{\rho_{1}}(\infty), \hat{\vec {u}_{1}}-\hat{\vec {u}_{1}}(\infty)), \hat{\theta_{1}}-\hat{\theta_{1}}(\infty)),\\
&\frac{\partial \hat{\rho_{1}}}{\partial\tau}=0,~~\frac{\partial \hat{\vec u_{1}}}{\partial\tau}=\frac{1}{\rho^{0}}\langle\vec w(f_{I, 0}-\overline f_{I, 0})\rangle,~~\frac{\partial \hat{\theta_{1}}}{\partial\tau}=0,\\
&(\hat{\rho_{1}}(0, \vec x), \hat{\vec u_{1}}(0, \vec x), \hat{\theta_{1}}(0, \vec x))=(0, 0, 0),\\
&\lim_{\tau\rightarrow\infty}(\hat{\rho_{1}}(\tau, \vec x), \hat{\vec{u}}_{1}(\tau, \vec x), \hat{\theta_{1}}(\tau, \vec x))=(\hat{\rho_{1}}(\infty, \vec x), \hat{\vec{u}}_{1}(\infty, \vec x), \hat{\theta_{1}}(\infty, \vec x)).\\
\end{split}\right.
\end{equation}

After a careful calculation, we have $(\rho_{I, 1}, \vec{u}_{I, 1}, \theta_{I, 1})=(0, \frac{1}{\rho^{0}}(1-e^{-\tau})\langle\vec w h\rangle, 0)$ which implies $\hat{\rho_{1}}(\infty, \vec x), \hat{\vec{u}}_{1}(\infty, \vec x), \hat{\theta_{1}}(\infty, \vec x)\in W_{2}^{6}(\Omega)$.

We write the equations about $\rho_{1}, \vec{u}_{1}, \theta_{1}$ and $f_{1}$ as follows
\begin{equation}\label{equations about theta 1 overline f epsilon 1}\left\{
\begin{split}
&f_{1}(t, x, w)=\overline f_{1}-\vec w\cdot\nabla_{x}f_{0},\\
&\frac{\partial\rho_{1}}{\partial t}+\vec{u}_{0}\cdot\nabla_{x}\rho_{1}+\rho_{0}\mathrm{div}_{x}\vec{u}_{1}=-\mathrm{div}_{x}\vec{u}_{0}\rho_{1}-\vec{u}_{1}\cdot\nabla_{x}\rho_{0}\ \ \mathrm{in}\ \ (0, T]\times\Omega,\\
&\frac{\partial\vec{u}_{1}}{\partial t}+\vec{u}_{0}\cdot\nabla_{x}\vec{u}_{1}+\nabla_{x}\theta_{1}+\frac{\theta_{0}}{\rho_{0}}\nabla_{x}\rho_{1}
=-\frac{\rho_{1}}{\rho_{0}}\frac{\partial \vec{u}_{0}}{\partial t}-\vec{u}_{1}\cdot\nabla_{x}\vec{u}_{0}-\frac{\rho_{1}}{\rho_{0}}\vec{u}_{0}\cdot\nabla_{x}\vec{u}_{0}
\\&-\frac{\rho_{1}}{\rho_{0}}\nabla_{x}\theta_{0}-\frac{\theta_{1}}{\rho_{0}}\nabla_{x}\rho_{0}-\frac{1}{3}\nabla_{x}\overline f_{1}\ \ \mathrm{in}\ \ (0, T]\times\Omega,\\
&\frac{\partial\theta_{1}}{\partial t}+\vec{u}_{0}\cdot\nabla_{x}\theta_{1}+\theta_{0}\mathrm{div}_{x}\vec{u}_{1}=\frac{\overline f_{1}-B_{1}(\theta_{0}; \theta_{1})}{\rho_{0}}-\frac{\rho_{1}}{\rho_{0}}\frac{\partial\theta_{0}}{\partial t}-\vec{u}_{1}\cdot\nabla_{x}\theta_{0}-\frac{\rho_{1}}{\rho_{0}}\vec{u}_{0}\cdot\nabla_{x}\theta_{0}
\\&-\theta_{1}\mathrm{div}_{x}\vec{u}_{0}-\frac{\rho_{1}}{\rho_{0}}\theta_{0}\mathrm{div}_{x}\vec{u}_{0}+\frac{1}{3}\nabla_{x}\overline f_{1}\cdot\vec{u}_{0}\ \ \mathrm{in}\ \ (0, T]\times\Omega,\\
&\frac{\partial\overline f_{1}}{\partial t}-\frac{1}{3}\Delta_{x}\overline f_{1}+\overline f_{1}=B_{1}(\theta_{0}; \theta_{1})\ \ \mathrm{in}\ \ (0, T]\times\Omega,\\
&\rho_{1}(0, \vec{x})=0,\;\;\vec{u}_{1}(0, \vec{x})=\hat{\vec{u}}_{1}(\infty, \vec x)~~\mathrm{in}  \ \ \Omega,\\
&\theta_{1}(0, \vec{x})=0,\;\;\overline f_{1}(0, \vec{x})=\mathcal{F}_{1}(\infty, \vec x)~~\mathrm{in}  \ \ \Omega,
\end{split}\right.
\end{equation}

Note that $B_{1}(\theta_{0}; \theta_{1})=4\theta_{0}^{3}\theta_{1}$. Similar to the proof of the first part of Theorem \ref{definition weak solution-0}, we can write system (\ref{equations about theta 1 overline f epsilon 1}) as follows
\begin{equation}\label{equations about theta 1 overline f epsilon 1 jia1}\left\{
\begin{split}
&f_{1}(t, x, w)=\overline f_{1}-\vec w\cdot\nabla_{x}f_{0},\\
&A_{0}(V^{0})\frac{\partial V^{1}}{\partial t}+\sum_{j=1}^{3}A_{j}(V^{0})\frac{\partial V^{1}}{\partial x_{j}}+D(V^{0})V^{1}=F(\overline f_{1}, \nabla_{x}\overline f_{1})\ \ \mathrm{in}\ \ (0, T]\times\Omega,\\
&\frac{\partial\overline f_{1}}{\partial t}-\frac{1}{3}\Delta_{x}\overline f_{1}+\overline f_{1}=B_{1}(\theta_{0}; \theta_{1})\ \ \mathrm{in}\ \ (0, T]\times\Omega,\\
&V^{1}(0, \vec x)=(0, \langle\vec w h\rangle, 0)^{t}~~\mathrm{in}  \ \ \Omega,\\
&\overline f_{1}(0, \vec{x})=\mathcal{F}_{1}(\infty, \vec x)~~\mathrm{in} \ \ \Omega,\\
\end{split}\right.
\end{equation}
where $V^{0}=(\rho_{0}, \vec{u}_{0}, \theta_{0}), V^{1}=(\rho_{1}, \vec{u}_{1}, \theta_{1})$,
\begin{equation}\label{xishujuzhen-1-1}
A_{0}(V^{0})
=
\left[
\begin{array}{cccc}
(\rho_{0})^{-1}&0&0\\
0&\frac{\rho_{0}}{\theta_{0}}\mathbb{I}_{3\times3}&0\\
0&0&\frac{\rho_{0}}{(\theta_{0})^{2}}
\end{array}
\right]
\end{equation}
and $A_{j}$ is a symmetric matrix defined similarly as in (\ref{research equations-3-2}).

For a given $\overline f_{1}'\in L^{2}(0, T; W_{2}^{3}(\Omega))\cap C^{0}([0, T]; W_{2}^{2}(\Omega))$, we consider the following problem
\begin{equation}\label{research equation-2-1-1}\left\{
\begin{split}
&A_{0}(V^{0})\frac{\partial V^{1}}{\partial t}+\sum_{j=1}^{3}A_{j}(V^{0})\frac{\partial V^{1}}{\partial x_{j}}+D(V^{0})V^{1}=F(\overline f_{1}', \nabla_{x}\overline f_{1}')\ \ \mathrm{in}\ \ (0, T]\times\Omega,\\
&V^{1}(0, \vec{x})=0\ \ \mathrm{in}\ \ \Omega.\\
\end{split}\right.
\end{equation}
According to the similar version of Theorem I in \cite{cauchy-problem}, we obtain that there exists a unique solution $V^{1}\in C^{0}([0, T]; W_{2}^{2}(\Omega))$ of (\ref{research equation-2-1-1}) through the density argument. Define the map
\begin{equation}\nonumber
S:\overline f_{1}'\in L^{2}(0, T; W_{2}^{3}(\Omega))\cap C^{0}([0, T]; W_{2}^{2}(\Omega))\mapsto V^{1}\in C^{0}([0, T]; W_{2}^{2}(\Omega)).
\end{equation}
Similarly, we consider the parabolic equation
\begin{equation}\label{research equations-2-1-1}\left\{
\begin{split}
&\frac{\partial\overline f_{1}^{s}}{\partial t}-\frac{1}{3}\Delta_{x}\overline f_{1}^{s}+\overline f_{1}^{s}=sB_{1}(\theta_{0}; S(\overline f_{1}'))\ \ \mathrm{in}\ \ (0, T]\times\Omega,\\
&\overline f_{1}^{s}(0, \vec{x})=s\mathcal{F}_{1}(\infty, \vec{x})\ \ \mathrm{in}  \ \ \Omega.\\
\end{split}\right.
\end{equation}
According to the parabolic theory, we can obtain that there exists a unique solution $\overline f_{1}^{s}\in L^{2}(0, T; W_{2}^{4}(\Omega))\cap L^{\infty}(0, T; W_{2}^{3}(\Omega))$ for any given $\overline f_{1}'\in L^{2}(0, T; W_{2}^{2}(\Omega))\cap C^{0}([0, T]; W_{2}^{1}(\Omega)), \\s\in[0, 1].$ Note that $\partial_{t}\overline f_{1}^{s}\in L^{2}(0, T; W_{2}^{2}(\Omega))\cap L^{\infty}(0, T; W_{2}^{1}(\Omega))$. According to the Aubin-Lions-Simon lemma, we can define a compact map
\begin{equation}\nonumber
\begin{aligned}
\Gamma_{s}:\overline f_{1}'\in L^{2}(0, T; W_{2}^{3}(\Omega))\cap C^{0}([0, T]; W_{2}^{2}(\Omega))\mapsto&\overline f_{1}^{s}\in \{f|f\in L^{2}(0, T; W_{2}^{4}(\Omega))\cap L^{\infty}(0, T; \\&W_{2}^{3}(\Omega)),\partial_{t}f\in L^{2}(0, T; W_{2}^{2}(\Omega))\cap L^{\infty}(0, T; \\&W_{2}^{1}(\Omega))\}.
\end{aligned}
\end{equation}

Now, we assume that $\overline f_{1}^{s}$ and $\theta_{1}=S(\overline f^{s}_{1})$ satisfy
\begin{equation}\label{research equations-2-1-2}\left\{
\begin{split}
&\frac{\partial\overline f_{1}^{s}}{\partial t}-\frac{1}{3}\Delta_{x}\overline f_{1}^{s}+\overline f_{1}^{s}=sB_{1}(\theta_{0}; \theta_{1})\ \ \mathrm{in}\ \ (0,T]\times\Omega,\\
&A_{0}(V^{0})\frac{\partial V^{1}}{\partial t}+\sum_{j=1}^{3}A_{j}(V^{0})\frac{\partial V^{1}}{\partial x_{j}}+D(V^{0})V^{1}=F(\overline f_{1}^{s}, \nabla_{x}\overline f_{1}^{s})\ \ \mathrm{in}\ \ (0, T]\times\Omega,\\
&\overline f_{1}^{s}(0, \vec{x})=s\mathcal{F}_{1}(\infty, \vec{x}),\;\;V^{1}(0, \vec{x})=0\ \ \mathrm{in}\ \ \Omega.
\end{split}\right.
\end{equation}
Multiplying $\overline f^{s}_{1}$ and $\theta_{1}$ on both sides of equation $(\ref{research equations-2-1-2})_{1}$ and $(\ref{research equations-2-1-2})_{2}$ respectively, for small $\delta$, we derive
\begin{equation}\nonumber
\begin{split}
&\int_{\Omega}|\overline f_{1}^{s}|^{2}\mathrm{d}\vec{x}(t)+\int_{\Omega}|V^{1}|^{2}\mathrm{d}\vec{x}(t)+\frac{1}{3}\int_{0}^{t}\int_{\Omega}|\nabla_{x}\overline f_{1}^{s}|^{2}\mathrm{d}\vec{x}\mathrm{d}\tau+\int_{0}^{t}\int_{\Omega}|\overline f_{1}^{s}|^{2}\mathrm{d}\vec{x}d\tau\\\leq&\lambda_{\min}^{-1}\Big(s\int_{0}^{t}\int_{\Omega}B_{1}(\theta_{0}; \theta_{1})\overline f_{1}^{s}\mathrm{d}\vec{x}\mathrm{d}\tau+\int_{0}^{t}\int_{\Omega}|\overline f_{1}^{s}|^{2} \mathrm{d}\vec{x}\mathrm{d}\tau \\&+\frac{1}{2}\int_{0}^{t}\int_{\Omega}\sum_{j=1}^{3}\frac{\partial A_{j}(V^{0})}{\partial x_{j}}V^{1}\cdot V^{1}\mathrm{d}\vec x \mathrm{d}\tau+C(\delta)\int_{0}^{t}\int_{\Omega}(|V^{1}|^{2}+1)\mathrm{d}\vec{x}(t)+\delta\int_{0}^{t}\int_{\Omega}|\nabla_{x}\overline f_{1}^{s}|^{2}\mathrm{d}\vec{x}\mathrm{d}\tau\Big)\\&+s^{2}\int_{\Omega}|\mathcal{F}_{1}(\infty, \vec{x})|^{2}\mathrm{d}\vec{x}+\int_{\Omega}|q_{1}(\infty, \vec{x})|^{2}\mathrm{d}\vec{x}.
\end{split}
\end{equation}
Then, we have
\begin{equation}\nonumber
\begin{aligned}
&\int_{\Omega}|\overline f_{1}^{s}|^{2}\mathrm{d}\vec{x}(t)+\int_{\Omega}|V^{1}|^{2}\mathrm{d}\vec{x}(t)+\int_{0}^{t}\int_{\Omega}|\nabla_{x}\overline f_{1}^{s}|^{2}\mathrm{d}\vec{x}\mathrm{d}\tau\\\leq&C(a, r)\Big(\int_{0}^{t}\int_{\Omega}|\overline f_{1}^{s}|^{2}\mathrm{d}\vec{x}\mathrm{d}\tau+\int_{0}^{t}\int_{\Omega}|V^{1}|^{2}\mathrm{d}\vec{x}\mathrm{d}\tau\Big)+C.
\end{aligned}
\end{equation}
According to the Gronwall's inequality, we obtain $\overline f_{1}^{s}\in L^{\infty}(0, T; L^{2}(\Omega))\cap L^{2}(0, T; W_{2}^{1}(\Omega))$ and $V^{1}\in L^{\infty}(0, T; L^{2}(\Omega))$. After a bootstrap, we have $\overline f_{1}^{s}\in L^{\infty}(0, T; W_{2}^{5}(\Omega))\cap L^{2}(0, T; \\W_{2}^{6}(\Omega))$, $V^{1}\in L^{\infty}(0, T; W_{2}^{5}(\Omega))$ and $\partial_{t}\overline f_{1}^{s}\in L^{\infty}(0, T; W_{2}^{3}(\Omega))$, and these estimates are independent of $s$. Therefore, we can obtain that there exists a solution $(\overline{f}_{1}, \theta_{1})$ of problem (\ref{equations about theta 1 overline f epsilon 1}) such that $\overline f_{1}\in L^{\infty}(0, T; W_{2}^{5}(\Omega))\cap L^{2}(0, T; W_{2}^{6}(\Omega))$ and $V^{1}=S(\overline f_{1})\in L^{\infty}(0, T; W_{2}^{5}(\Omega))$ according to the Schaefer fixed point theorem.

\textbf{Step 3.} Construction of second-order terms.

The second-order initial layer is defined as
\begin{equation}\label{equation about f I2}\left\{
\begin{split}
&f_{I, 2}(\tau, \vec{x}, \vec{w})=\mathcal{F}_{2}(\tau, \vec{x}, \vec{w})-\mathcal{F}_{2}(\infty, \vec{x}),\\
&\frac{\partial\overline {\mathcal{F}_{2}}}{\partial\tau}=-\frac{1}{4\pi}\int_{\mathbb{S}^{2}}(B(\theta^{0}+\theta_{I, 0})-B(\theta^{0})-\vec w\cdot\nabla_{x}f_{I, 1}-f_{I, 0})\mathrm{d}\vec{w},\\
&\mathcal{F}_{2}(\tau, \vec{x}, \vec{w})=e^{-\tau}\mathcal{F}_{2}(0, \vec{x}, \vec{w})+\int_{0}^{\tau}(\overline {\mathcal{F}_{2}}+B(\theta^{0}+\theta_{I, 0})-B(\theta^{0})-\vec w\cdot\nabla_{x}f_{I, 1}-f_{I, 0})e^{s-\tau}\mathrm{d}s.\\
&\mathcal{F}_{2}(0, \vec{x}, \vec{w})=(f_{0}+\partial_{t}f_{0}+\vec w\cdot\nabla_{x}f_{1}-B(\theta_{0}))|_{t=0},\\
&\lim_{\tau\rightarrow\infty}\mathcal{F}_{2}(\tau, \vec{x}, \vec{w})=\mathcal{F}_{2}(\infty, \vec{x}).
\end{split}\right.
\end{equation}

Since $\theta_{I, 0}=0$, we have $B(\theta^{0}+\theta_{I, 0})-B(\theta^{0})=0$. Since $-\vec w\cdot\nabla_{x}f_{I, 1}-f_{I, 0}\in C^{\infty}([0, \infty); L^{\infty}(\mathbb{S}^{2}; W_{2}^{4}(\Omega)))$ and $\mathcal{F}_{2}(0, \vec{x}, \vec{w})=(\overline f_{0}+\partial_{t}f_{0}+\vec w\cdot\nabla_{x}f_{1}-B(\theta_{0}))|_{t=0}\in W_{2}^{4}(\Omega)$, we obtain that $\mathcal{F}_{2}(\tau, \vec{x}, \vec{w})$ is well-defined. Then, we get $f_{I, 2}\in C^{\infty}([0, \infty); L^{\infty}(\mathbb{S}^{2}; W_{2}^{4}(\Omega)))$. Finally, we can obtain $\mathcal{F}_{2}(\infty, \vec{x})\in W_{2}^{4}(\Omega).$

Since $\rho_{I, 0}, \vec{u}_{I, 0}, \theta_{I, 0}=0$, we can get the following equations about $\rho_{I, 2}, \vec{u}_{I, 2}, \theta_{I, 2}$
\begin{equation}\label{equation about theta I2-1}\left\{
\begin{split}
&(\rho_{I, 2}, \vec{u}_{I, 2}, \theta_{I, 2})=(\hat{\rho_{2}}-\hat{\rho_{2}}(\infty), \hat{\vec {u}_{2}}-\hat{\vec {u}_{2}}(\infty)), \hat{\theta_{2}}-\hat{\theta_{2}}(\infty)),\\
&\frac{\partial \hat{\rho_{2}}}{\partial\tau}=0,~~\frac{\partial \hat{\vec u_{2}}}{\partial\tau}=\frac{1}{\rho^{0}}\langle\vec w(f_{I, 1}-\overline f_{I, 1})\rangle,\\
&\frac{\partial \hat{\theta_{2}}}{\partial\tau}=-\frac{1}{\rho^{0}}(\langle\vec w(f_{I, 1}-\overline f_{I, 1})\rangle\cdot\vec{u}^{ 0}+\langle\vec w(f_{I, 0}-\overline f_{I, 0})\rangle\cdot\vec{u}_{1}^{0}),\\
&(\hat{\rho_{2}}(0, \vec x), \hat{\vec u_{2}}(0, \vec x), \hat{\theta_{2}}(0, \vec x))=(0, 0, 0),\\
&\lim_{\tau\rightarrow\infty}(\hat{\rho_{2}}(\tau, \vec x), \hat{\vec{u}}_{2}(\tau, \vec x), \hat{\theta_{2}}(\tau, \vec x))=(\hat{\rho_{2}}(\infty, \vec x), \hat{\vec{u}}_{2}(\infty, \vec x), \hat{\theta_{2}}(\infty, \vec x)).\\
\end{split}\right.
\end{equation}
Similar to the analysis of equations (\ref{equation about theta I1}), we have that $\rho_{I, 2}=0, \vec{u}_{I, 2}, \theta_{I, 2}$ are well-defined. Note that $\hat{\rho_{2}}(\infty, \vec x), \hat{\vec {u}_{2}}(\infty, \vec x)), \hat{\theta_{2}}(\infty, \vec x)\in W_{2}^{5}(\Omega)$.

We write the equations about $\rho_{2}, \vec{u}_{2}, \theta_{2}$ and $f_{2}$ as follows
\begin{equation}\label{equations about overline f2 epsilon theta 2}\left\{
\begin{split}
&f_{2}(t, x, w)=\overline f_{2}+B(\theta_{0})-f_{0}-\partial_{t}f_{0}-\vec w\cdot\nabla_{x}f_{1},\\
&\frac{\partial\rho_{2}}{\partial t}+\vec{u}_{0}\cdot\nabla_{x}\rho_{2}+\rho_{0}\mathrm{div}_{x}\vec{u}_{2}=-\vec{u}_{2}\cdot\nabla_{x}\vec{u}_{0}
-\vec{u}_{1}\cdot\nabla_{x}\rho_{1}-\rho_{2}\mathrm{div}_{x}\vec{u}_{0}\\&-\rho_{1}\mathrm{div}_{x}\vec{u}_{1}\ \ \mathrm{in}\ \ (0, T]\times\Omega,\\
&\frac{\partial\vec{u}_{2}}{\partial t}+\vec{u}_{0}\cdot\nabla_{x}\vec{u}_{2}+\nabla_{x}\theta_{2}+\frac{\theta_{0}}{\rho_{0}}\nabla_{x}\rho_{2}=-\frac{\rho_{1}}{\rho_{0}}\frac{\partial\vec{u}_{1}}{\partial t}-\frac{\rho_{2}}{\rho_{0}}\frac{\partial\vec{u}_{0}}{\partial t}-\vec{u}_{1}\cdot\nabla_{x}\vec{u}_{1}+\vec{u}_{2}\cdot\nabla_{x}\vec{u}_{0}
\\&-\frac{\rho_{1}}{\rho_{0}}(\vec{u}_{0}\cdot\nabla_{x}\vec{u}_{1}+\vec{u}_{1}\cdot\nabla_{x}\vec{u}_{0})
-\frac{\rho_{2}}{\rho_{0}}\vec{u}_{0}\cdot\nabla_{x}\vec{u}_{0}-\frac{\rho_{1}}{\rho_{0}}\nabla_{x}\theta_{1}-\frac{\rho_{2}}{\rho_{0}}\nabla_{x}\theta_{0}
-\frac{\theta_{1}}{\rho_{0}}\nabla_{x}\rho_{1}-\frac{\theta_{2}}{\rho_{0}}\nabla_{x}\rho_{0}\\&+\frac{1}{3\rho_{0}}(2\nabla_{x}\overline f_{0}-\nabla_{x}\overline f_{2}+\partial_{t}\nabla_{x}\overline f_{0}+\nabla_{x}\overline f_{1}-\nabla_{x}B(\theta_{0}))\ \ \mathrm{in}\ \ (0, T]\times\Omega,\\
&\frac{\partial\theta_{2}}{\partial t}+\vec{u}_{0}\cdot\nabla_{x}\theta_{2}+\theta_{0}\mathrm{div}_{x}\vec{u}_{2}=\frac{\overline f_{2}-B_{2}(\theta_{0}; \theta_{2})}{\rho_{0}}-\frac{\rho_{1}}{\rho_{0}}\frac{\partial\theta_{1}}{\partial t}-\frac{\rho_{2}}{\rho_{0}}\frac{\partial\theta_{0}}{\partial t}
-\frac{\rho_{2}}{\rho_{0}}\vec{u}_{0}\cdot\nabla_{x}\theta_{0}\\&-\frac{\rho_{1}}{\rho_{0}}(\vec{u}_{1}\cdot\nabla_{x}\theta_{0}+\vec{u}_{0}\cdot\nabla_{x}\theta_{1})
-(\vec{u}_{1}\cdot\nabla_{x}\theta_{1}+\vec{u}_{2}\cdot\nabla_{x}\theta_{0})
-\frac{\rho_{2}\theta_{0}}{\rho_{0}}\mathrm{div}\vec{u}_{0}\\&-\frac{\rho_{1}}{\rho_{0}}(\theta_{1}\mathrm{div}\vec{u}_{0}+\theta_{0}\mathrm{div}\vec{u}_{1})
-(\theta_{1}\mathrm{div}\vec{u}_{1}+\theta_{2}\mathrm{div}\vec{u}_{0})\\&-\frac{1}{3\rho_{0}}(2\nabla_{x}\overline f_{0}-\nabla_{x}\overline f_{2}+\partial_{t}\nabla_{x}\overline f_{0}+\nabla_{x}\overline f_{1}-\nabla_{x}B(\theta_{0}))\cdot\vec{u}_{0}\ \ \mathrm{in}\ \ (0, T]\times\Omega,\\
&\frac{\partial\overline f_{2}}{\partial t}-\frac{1}{3}\Delta_{x}\overline f_{2}+\overline f_{2}=B_{2}(\theta_{0}; \theta_{2})+\frac{1}{3}\Delta_{x}B(\theta_{0})+\frac{1}{4\pi}\int_{S^{2}}(\vec{w}\cdot\nabla_{x})^{4}\overline f_{0}\mathrm{d}\vec w\ \ \mathrm{in}\ \ (0, T]\times\Omega,\\
&\rho_{2}(0, \vec{x})=0,\;\;\vec{u}_{2}(0, \vec{x})=\hat{\vec {u}_{2}}(\infty, \vec x)~~\mathrm{in}  \ \ \Omega,\\
&\theta_{2}(0, \vec{x})=\hat{\theta_{2}}(\infty, \vec x),\;\;\overline f_{2}(0, \vec{x})=\mathcal{F}_{2}(\infty, \vec x)~~\mathrm{in}  \ \ \Omega.
\end{split}\right.
\end{equation}

Note that $B_{2}(\theta_{0}; \theta_{2})=4\theta_{2}\theta_{0}^{3}+6\theta_{1}^{2}\theta_{0}^{2}$. Similar to the analysis of system (\ref{equations about theta 1 overline f epsilon 1}), we can obtain that there exists a solution $(\overline f_{2}, \rho_{2}, \vec{u}_{2}, \theta_{2})\in (L^{\infty}(0, T; W_{2}^{4}(\Omega))\cap L^{2}(0, T; W_{2}^{5}(\Omega)))\times L^{\infty}(0, T; W_{2}^{4}(\Omega))\times L^{\infty}(0, T; W_{2}^{4}(\Omega))\times L^{\infty}(0, T; W_{2}^{4}(\Omega))$ about the system (\ref{equations about overline f2 epsilon theta 2}).

\subsection{Proof of Theorem \ref{definition weak solution-1}.}
We divide the proof into four steps:

\noindent\textbf{Step 1.} Remainder definitions.
We may rewrite the asymptotic expansion as follows:
\begin{equation}\label{expansion of f epsilon}
f^{\epsilon}\sim \sum_{k=0}^{\infty}\epsilon^{k}f_{k}+\sum_{k=0}^{\infty}\epsilon^{k}f_{I, k}.
\end{equation}
The remainder can be defined as
\begin{equation}\nonumber
\begin{split}
R_{N}&=f^{\epsilon}-\sum_{k=0}^{N}\epsilon^{k}f_{k}-\sum_{k=0}^{N}\epsilon^{k}f_{I, k}\\&=f^{\epsilon}-Q_{N}-\mathcal{D}_{I, N},
\end{split}
\end{equation}
where
\begin{equation}\nonumber
Q_{N}=\sum_{k=0}^{N}\epsilon^{k}f_{k}
\end{equation}
and
\begin{equation}\nonumber
\mathcal{D}_{I, N}=\sum_{k=0}^{N}\epsilon^{k}f_{I, k}.
\end{equation}

We note that the equation $(\ref{research equations})_{1}$ is equivalent to the equation $(\ref{research equation in the initial layer})_{1}$, we write $\mathcal{L}$ to denote the transport operator as follows:

\begin{equation}\label{two kinds of equations about f epsilon}
\begin{split}
\mathcal{L}f^{\epsilon}=&\epsilon^{2}\partial_{t}f^{\epsilon}+\epsilon\vec w\cdot\nabla_{x}f^{\epsilon}+f^{\epsilon}-\overline f^{\epsilon}+\epsilon^{2}f^{\epsilon}\\=&\frac{\partial f^{\epsilon}}{\partial\tau}+\epsilon \vec w\cdot \nabla_{x}f^{\epsilon}+(f^{\epsilon}-\overline f^{\epsilon})+\epsilon^{2}f^{\epsilon}.
\end{split}
\end{equation}

\noindent\textbf{Step 2.} Estimates of $\mathcal{L}Q_{2}.$
\begin{equation}\label{equality about Q2}
\begin{split}
\mathcal{L}Q_{2}&=\epsilon^{2}\partial_{t}Q_{2}+\epsilon\vec w\cdot\nabla_{x}Q_{2}+(Q_{2}-\overline Q_{2})+\epsilon^{2}Q_{2}\\&=\epsilon^{2}B(\theta_{0})+\epsilon^{3}\partial_{t}f_{1}+\epsilon^{4}\partial_{t}f_{2}
+\epsilon^{3}\vec{w}\cdot\nabla_{x}f_{2}
+\epsilon^{3}f_{1}+\epsilon^{4}f_{2}.\\
\end{split}
\end{equation}

We have
\begin{equation}
\|\epsilon^{3}\partial_{t}f_{1}+\epsilon^{4}\partial_{t}f_{2}\|_{L^{\infty}((0, T)\times\Omega\times \mathbb{S}^{2})}\leq C\epsilon^{3},
\end{equation}

\begin{equation}
\|\epsilon^{3}\vec{w}\cdot\nabla_{x}f_{2}\|_{L^{\infty}((0, T)\times\Omega\times \mathbb{S}^{2})}\leq C\epsilon^{3}
\end{equation}
and
\begin{equation}\label{inequality about f3-overline f3}
\|\epsilon^{3}f_{1}+\epsilon^{4}f_{2}\|_{L^{\infty}((0, T)\times\Omega\times \mathbb{S}^{2})}\leq C\epsilon^{3}.
\end{equation}

\noindent\textbf{Step 3.} Estimates of $\mathcal{L}D_{I, 2}.$
\begin{equation}\label{equality about mathcal D I2}
\begin{split}
\mathcal{L}\mathcal{D}_{I, 2}=&\frac{\partial \mathcal{D}_{I, 2}}{\partial\tau}+\epsilon \vec w\cdot \nabla_{x}\mathcal{D}_{I, 2}+(\mathcal{D}_{I, 2}-\overline {\mathcal{D}_{I, 2}})+\epsilon^{2}\mathcal{D}_{I, 2}\\=&\frac{\partial f_{I, 0}}{\partial\tau}+\epsilon \vec w\cdot \nabla_{x}f_{I, 0}+(f_{I, 0}-\overline f_{I, 0})+\epsilon^{2}f_{I, 0}\\&+\epsilon\Big(\frac{\partial f_{I, 1}}{\partial\tau}+\epsilon \vec w\cdot \nabla_{x}f_{I, 1}+(f_{I, 1}-\overline f_{I, 1})+\epsilon^{2}f_{I, 1}\Big)\\&+\epsilon^{2}\Big(\frac{\partial f_{I, 2}}{\partial\tau}+\epsilon \vec w\cdot \nabla_{x}f_{I, 2}+(f_{I, 2}-\overline f_{I, 2})+\epsilon^{2}f_{I, 2}\Big).
\end{split}
\end{equation}

Since $B(\theta^{0}+\theta_{I, 0})-B(\theta^{0})=0,$ we have
\begin{equation}
\begin{split}
\mathcal{L}\mathcal{D}_{I, 2}=\epsilon^{3}f_{I, 1}+\epsilon^{3}\vec w\cdot\nabla_{x}f_{I, 2}+\epsilon^{4}f_{I, 2}.
\end{split}
\end{equation}

Furthermore,
\begin{equation}
\|\epsilon^{3}f_{I, 1}^{\epsilon}+\epsilon^{3}\vec w\cdot\nabla_{x}f_{I, 2}\|_{L^{\infty}((0, T)\times\Omega\times \mathbb{S}^{2})}\leq C\epsilon^{3},
\end{equation}
\begin{equation}\label{inequality about f I2}
\|\epsilon^{4}f_{I, 2}\|_{L^{\infty}((0, T)\times\Omega\times \mathbb{S}^{2})}\leq C\epsilon^{4},
\end{equation}
combining with (\ref{equality about f I0}), (\ref{equality about f I0 f I1}) and (\ref{equality about f I0 f I1 f I2}).

\noindent\textbf{Step 4.} Uniform convergence

In summary, since $\mathcal{L}f^{\epsilon}=\epsilon^{2}B(\theta^{\epsilon}),$  we have
\begin{equation}\nonumber
\begin{split}
\|\mathcal{L}\mathcal{R}_{2}\|_{L^{\infty}((0, T)\times\Omega\times \mathbb{S}^{2})}&=\|\mathcal{L}f^{\epsilon}-\mathcal{L}Q_{2}-\mathcal{L}\mathcal{D}_{I, 2}\|_{L^{\infty}((0, T)\times\Omega\times \mathbb{S}^{2})}
\\&\leq\epsilon^{2}\|B(\theta^{\epsilon})-B(\theta_{0})\|_{L^{\infty}((0, T)\times\Omega)}+C\epsilon^{3},
\end{split}
\end{equation}by collecting (\ref{equality about Q2})$\sim$(\ref{inequality about f3-overline f3}) and (\ref{equality about mathcal D I2})$\sim$(\ref{inequality about f I2}).

Now we consider the estimate of $\|B(\theta^{\epsilon})-B(\theta_{0})\|_{L^{\infty}((0, T)\times\Omega)}$. Since $\overline f_{0}=\overline f$ and $\theta_{0}=\hat\theta$ by the uniqueness of system (\ref{equations about limit of overline f epsilon}), we obtain
\begin{equation}\label{convergence L2 theta epsilon}
\lim_{\epsilon\rightarrow 0}\|\theta^{\epsilon}-\theta_{0}\|_{C^{0}([0, T]; W_{2}^{5}(\Omega))}=0.
\end{equation}
Since $\|\theta^{\epsilon}\|_{L^{\infty}((0, T)\times\Omega)}+ \|\theta_{0}\|_{L^{\infty}((0, T)\times\Omega)}\leq C(r)$, we conclude that
\begin{equation}\label{convergence L2 B theta epsilon}
\|B(\theta^{\epsilon})-B(\theta_{0})\|_{L^{\infty}((0, T)\times\Omega)}\leq C(r)\|\theta^{\epsilon}-\theta_{0}\|_{C^{0}([0, T]; W_{2}^{5}(\Omega))}.
\end{equation}

Finally, we consider the asymptotic expansion to $N=2$ where the remainder $\mathcal{R}_{2}$ satisfies the following equation
\begin{equation}\label{remainder equation-2}\left\{
\begin{split}
&{\epsilon}^{2}\partial_{t}\mathcal{R}_{2}+\epsilon \vec w\cdot \nabla_{x}\mathcal{R}_{2}+\mathcal{R}_{2}-\overline {\mathcal{R}_{2}}+{\epsilon}^{2}\mathcal{R}_{2}=\mathcal{L}\mathcal{R}_{2}\ \ \mathrm{in} \ \ (0, T]\times \Omega\times \mathbb{S}^{2},\\
&\mathcal{R}_{2}(0, \vec{x}, \vec{w})=0\ \ \mathrm{in}\ \ \Omega\times \mathbb{S}^{2}.\\
\end{split}\right.
\end{equation}

By Lemma 3.1 we have
\begin{equation}\nonumber
\|\mathcal{R}_{2}\|_{L^{\infty}((0,T)\times \Omega\times\mathbb{S}^{2})}\leq \frac{C}{\epsilon^{2}}\|\mathcal{L}\mathcal{R}_{2}\|_{L^{\infty}((0,T)\times \Omega\times \mathbb{S}^{2})}\leq C(r)\|\theta^{\epsilon}-\theta_{0}\|_{C^{0}([0, T]; W_{2}^{5}(\Omega))}+C\epsilon.
\end{equation}
Since
\begin{equation}\nonumber
\Big\|\sum_{k=1}^{2}\epsilon^{k}f_{k}^{\epsilon}+\sum_{k=1}^{2}\epsilon^{k}f_{I, k}^{\epsilon}\Big\|_{L^{\infty}((0,T)\times \Omega\times S^{2})}=O(\epsilon),
\end{equation}
 we have
\begin{equation}\nonumber
\|f^{\epsilon}-f_{0}-f_{I, 0}\|_{L^{\infty}((0, T)\times\Omega\times \mathbb{S}^{2})}\leq O(\epsilon)+C(r)\|\theta^{\epsilon}-\theta_{0}\|_{C^{0}([0, T]; W_{2}^{5}(\Omega))},
\end{equation}
which implies (\ref{jielunbudengshi 1-jia}) by the fact that $\overline f=\overline f_{0}=f_{0}$. \hfill$\Box$

\vskip 5mm

\noindent{\bf Acknowledgment.} {This research was partially supported by NSFC grants 12271423, 12071044, and 12131007,  and the Fundamental Research Funds for the Central Universities (No. xzy012022005).}

\end{CJK}
\end{document}